\documentclass{article}

\usepackage{arxiv}

\usepackage[utf8]{inputenc} 
\usepackage[T1]{fontenc}    
\usepackage{hyperref}       
\usepackage{url}            
\usepackage{booktabs}       
\usepackage{amsfonts}       
\usepackage{nicefrac}       
\usepackage{microtype}      
\usepackage{lipsum}		
\usepackage{graphicx}
\usepackage{natbib}
\usepackage{doi}
\usepackage{amsmath, amssymb, amsthm}
\usepackage{nccmath}
\usepackage{mathtools}
\usepackage{enumitem}
\usepackage{colortbl}
\usepackage{xcolor}
\usepackage[labelfont=bf, labelsep=period]{caption}
\usepackage{placeins}
\usepackage{algorithm, algorithmic}

\DeclareMathOperator*{\argmin}{arg\,min}
\newcommand*{\E}{\mathrm{E}}
\newcommand*{\Var}{\mathrm{Var}}
\newcommand*{\Cov}{\mathrm{Cov}}
\newcommand*{\COV}{\mathbf{Cov}}
\newcommand*{\tr}{\mathrm{tr}}

\newcommand*{\T}{^{\!\mathsf{T}}}
\newcommand*{\iT}{^{-\!\mathsf{T}}}

\newtheorem{proposition}{Proposition}

\newtheorem{definition}{Definition}
\newtheorem{lemma}{Lemma}

\newcommand{\mathprefix}[1]{${#1}$\nobreakdash\ignorespaces}

\newcolumntype{L}[1]{>{\raggedright\let\newline\\\arraybackslash\hspace{0pt}}m{#1}}
\newcolumntype{C}[1]{>{\centering\let\newline\\\arraybackslash\hspace{0pt}}m{#1}}
\newcolumntype{R}[1]{>{\raggedleft\let\newline\\\arraybackslash\hspace{0pt}}m{#1}}

\hypersetup{
  colorlinks   = true, 
  urlcolor     = black, 
  linkcolor    = black, 
  citecolor   =  black 
}

\title{Optimal subsampling designs}

\author{ \hspace{1mm} Henrik~Imberg \\
    \href{mailto:imbergh@chalmers.se}{\tt imbergh@chalmers.se}
	\And
	\hspace{1mm} Marina~Axelson-Fisk \\
    \href{mailto:marina.axelson-fisk@chalmers.se}{\tt marina.axelson-fisk@chalmers.se}
	\And
	\hspace{1mm}Johan~Jonasson \\
    \href{mailto:jonasson@chalmers.se}{\tt jonasson@chalmers.se}
    \And
    \vspace{-0.5cm}
    \\
    Department of Mathematical Sciences \\
	Chalmers University of Technology and University of Gothenburg \\
	SE-412 96 Gothenburg, Sweden
}

\date{}

\renewcommand{\shorttitle}{Optimal subsampling designs}

\hypersetup{
pdftitle={Optimal subsampling designs},
pdfauthor={Henrik~Imberg, Johan~Jonasson, Marina~Axelson-Fisk},
pdfkeywords={A-optimality, D-optimality, L-optimality, M-estimation, inverse probability weighting, unequal probability sampling.}
}

\begin{document}
\maketitle


\begin{abstract}
Subsampling is commonly used to overcome computational and economical bottlenecks in the analysis of finite populations and massive datasets. Existing methods are often limited in scope and use optimality criteria (e.g., A-optimality) with well-known deficiencies, such as lack of invariance to the measurement-scale of the data and parameterisation of the model. A unified theory of optimal subsampling design is still lacking. We present a theory of optimal design for general data subsampling problems, including finite population inference, parametric density estimation, and regression modelling. Our theory encompasses and generalises most existing methods in the field of optimal subdata selection based on unequal probability sampling and inverse probability weighting. We derive optimality conditions for a general class of optimality criteria, and present corresponding algorithms for finding optimal sampling schemes under Poisson and multinomial sampling designs. We present a novel class of transformation- and parameterisation-invariant linear optimality criteria which enjoy the best of two worlds: the computational tractability of A-optimality and invariance properties similar to D-optimality. The methodology is illustrated on an application in the traffic safety domain. In our experiments, the proposed invariant linear optimality criteria achieve 92--99\% D-efficiency with 90--95\% lower computational demand. In contrast, the A-optimality criterion has only 46\% and 60\% D-efficiency on two of the examples. 
\end{abstract}



\section{Introduction}
\label{sec:introduction}

Consider a $p$-dimensional parameter $\boldsymbol{\theta}_0$ defined by
\begin{align}
    \label{eq:theta0}
    \boldsymbol{\theta}_0 & = \argmin_{\boldsymbol{\theta} \in \boldsymbol{\boldsymbol{\Omega}}} \ell_0(\boldsymbol{\theta}), 
\end{align}
i.e., as the minimiser of some function $\ell_0(\boldsymbol{\theta})$ over some parameter space $\boldsymbol{\boldsymbol{\Omega}} \subset \mathbb{R}^p$. We assume further that $\boldsymbol{\theta}_0$ is unique, and that $\ell_0(\boldsymbol{\theta})$ is twice differentiable and can be written on the form
\begin{align}
    \label{eq:ell0}
    \ell_0(\boldsymbol{\theta}) & = \sum_{i \in \mathcal D}\ell_i(\boldsymbol{\theta}), \quad \ell_i(\boldsymbol{\theta}) = \ell(\boldsymbol{\theta}; \boldsymbol{v}_i),
\end{align}
with summation over some index set $\mathcal D = \{1, \ldots, N\}$, where $\boldsymbol{v}_i$ is a data vector associated with a member $i \in \mathcal D$. Under these assumptions, $\boldsymbol{\theta}_0$ may also be defined as the unique solution to the estimation equation
\begin{equation}
    \label{eq:grads}
    \sum_{i\in \mathcal D} \boldsymbol{\psi}_i(\boldsymbol{\theta}) = \boldsymbol{0}, \quad \boldsymbol{\psi}_i(\boldsymbol{\theta}) = \nabla_{\!\boldsymbol{\theta}} \ell(\boldsymbol{\theta}; \boldsymbol{v}_i). 
\end{equation}
The data is on the form $\boldsymbol{v}_i = \boldsymbol{y}_i \in \mathcal{Y}$ or $\boldsymbol{v}_i = (\boldsymbol{x}_i, \boldsymbol{y}_i) \in \mathcal X \times \mathcal{Y}$, where $\boldsymbol{y}_i$ is a response vector and $\boldsymbol{x}_i$ a vector of explanatory variables. We will generally not distinguish between the case with and without explanatory variables, and throughout we write the data as $(\boldsymbol{x}_i, \boldsymbol{y}_i)$, keeping in mind that the first entry may be null and $\mathcal X$ the empty set. One may interpret \eqref{eq:theta0}--\eqref{eq:ell0} as an empirical risk minimisation problem \citep{Vapnik1991}. Hence, we will refer to $\ell_0(\boldsymbol{\theta})$ as the (full-data) empirical risk and to $\boldsymbol{\theta}_0$ as the (full-data) empirical risk minimiser (ERM).

The setting above covers a broad range of inference problems, models, and estimation methods in statistics, including maximum likelihood estimation, generalised linear models \citep{Nelder1972, McCullagh1989}, quasi-likelihood methods \citep{Wedderburn1974}, and certain types of M-estimation \citep{Stefanski2002}. Some specific examples, which will be considered further in the Application and Examples in Section \ref{sec:application}, include: 
\begin{enumerate}[label = \roman*)]
    \item Finite population inference: consider a finite population of $N$ individuals, where each individual is associated with a non-random vector characteristic $\boldsymbol{y}_i$. The vector of finite population means $\frac{1}{N}\sum_{i=1}^N \boldsymbol{y}_i$ may be written on the form \eqref{eq:theta0}--\eqref{eq:grads} with $\boldsymbol{v}_i = \boldsymbol{y}_i$ and $\ell(\boldsymbol{\theta}; \boldsymbol{v}_i) = ||\boldsymbol{y}_i - \boldsymbol{\theta}||_2^2 = (\boldsymbol{y}_i - \boldsymbol{\theta})\T(\boldsymbol{y}_i - \boldsymbol{\theta})$.
    \label{list:finite_population_inference}
    \item Parametric density estimation: given independent and identically distributed data $y_1, \ldots, y_N$ from a probability distribution with density function $f_{\boldsymbol{\theta}}(y)$, the maximum likelihood estimate of $\boldsymbol{\theta}$ may be written on the form \eqref{eq:theta0}--\eqref{eq:grads} with $\boldsymbol{v}_i = y_i$ and $\ell(\boldsymbol{\theta}; \boldsymbol{v}_i) = -\log f_{\boldsymbol{\theta}}(y_i)$. 
    \label{list:parametric_density_estimation}
    \item Regression modelling: consider a random sample $\{(\boldsymbol{x}_i, y_i)\}_{i=1}^N$, a vector of regression coefficients $\boldsymbol{\theta}$, a (non-linear) model $f_{\boldsymbol{\theta}}(\boldsymbol{x})$ for the conditional mean of $Y$ given $\boldsymbol{x}$, and a differentiable loss-function $l:\mathbb{R}^2 \to \mathbb{R}_+$ such that $l(\hat y, y) = 0$ if and only if $\hat y = y$. With $\boldsymbol{v}_i = (\boldsymbol{x}_i, y_i)$ and $\ell(\boldsymbol{\theta}; \boldsymbol{v}_i) = l(f_{\boldsymbol{\theta}}(\boldsymbol{x}), y_i)$, the equations \eqref{eq:theta0}--\eqref{eq:grads} define an estimate of the vector of regression coefficients $\boldsymbol{\theta}$. 
    \label{list:regression_modelling}
\end{enumerate}

Now consider a situation where inference based on the full data $\{(\boldsymbol{x}_i, \boldsymbol{y}_i)\}_{i \in \mathcal D}$ is prohibited by economic or computational constraints. For instance, the index set may be so large that complete enumeration to observe the full data  $\{(\boldsymbol{x}_i, \boldsymbol{y}_i)\}_{i \in \mathcal D}$ is practically or economically unfeasible. This is the typical situation in finite population inference \citep[][]{Neyman1938, Hansen1943, Horvitz1952}. Some variables may be expensive to measure and hence affordable to observe only for a small number of instances $i \in \mathcal D$, a situation known as a measurement-constrained experiment \citep{Wang2017, Meng2021, Zhang2021, Imberg2022b}. Another example is when the full data  $\{(\boldsymbol{x}_i, \boldsymbol{y}_i)\}_{i \in \mathcal D}$ is available, but the size $N$ of the dataset is so large that estimation of $\boldsymbol{\theta}$ using \eqref{eq:theta0}–\eqref{eq:ell0} is computationally unfeasible \citep{Ma2015, Drovandi2017, Wang2018, Deldossi2022, Dai2022}. In either case, we may search for an approximate solution based on a subset $\mathcal S \subset \mathcal D$ of size $n \ll N$.

In this paper we focus on methods based on data subsampling through unequal probability sampling and inverse probability weighting. Specifically, we consider an estimator of the form
\begin{align}
    \label{eq:thetahat}
    \hat{\boldsymbol{\theta}}_{\!\boldsymbol{\mu}} & = \argmin_{\boldsymbol{\theta} \in \boldsymbol{\Omega}} \hat{\ell}_{\boldsymbol{\mu}}(\boldsymbol{\theta}), \\
    \label{eq:ellhat}
    \hat{\ell}_{\boldsymbol{\mu}}(\boldsymbol{\theta}) & = \sum_{i \in \mathcal S} \frac{S_i}{\mu_i} \ell_i(\boldsymbol{\theta}),
\end{align}
where $S_i$ is the number of times an element $i \in \mathcal D$ is selected by the sampling mechanism, $\mu_i$ the corresponding expected number of selections, and $\mathcal S = \{i \in \mathcal D: S_i > 0\}$ the random set of selected elements. One may recognise \eqref{eq:ellhat} as the Hansen-Hurwitz estimator \citep{Hansen1943} of the full-data empirical risk function \eqref{eq:ell0}. Hence, we refer to $\hat{\boldsymbol{\theta}}_{\!\boldsymbol{\mu}}$ as the Hansen-Hurwitz empirical risk minimiser. For sampling without replacement, \eqref{eq:ellhat} coincides with the also well-known Horvitz-Thompson estimator of $\ell_0(\boldsymbol{\theta})$ \citep{Horvitz1952}. We also note that $\hat{\ell}_{\boldsymbol{\mu}}(\boldsymbol{\theta})$ is an unbiased estimator of $\ell_0(\boldsymbol{\theta})$, provided that $\mu_i > 0$ for all $i\in\mathcal D$, and $\hat{\boldsymbol{\theta}}_{\!\boldsymbol{\mu}}$ a consistent estimator of the full-data parameter $\boldsymbol{\theta}_0$ under general regularity conditions \citep{Binder1983}.

An important question to ask is how the subset $\mathcal S$ used for the approximate solution \eqref{eq:thetahat} to the problem \eqref{eq:theta0}–\eqref{eq:ell0} should be selected for optimal performance. The problem of optimal subsampling has a long standing tradition within the field of survey sampling for inference regarding finite populations; see, e.g., \citet{Neyman1938, Hajek1959, Cassel1976, Brewer1979} and \citet{Bellhouse1984}. Their work, however, is primarily concerned with linear estimators of scalar finite population characteristics. Stimulated by modern technological developments, the question of optimal subdata selection has attained renewed attention during the past few years also for more complex inference problems, as outlined above. Examples include leverage sampling and approximate numerical linear algebra methods for big data regression \citep{Ma2015, Ma2020}, optimal subsampling algorithms for binary and multinomial logistic regression \citep{Wang2018, Yao2019}, generalised linear models \citep{Ai2021_regression, Zhang2021, Yu2022}, quantile regression \citep{Ai2021_quantile, Wang2021}, and active learning \citep{Imberg2020, Kossen2022, Zhan2022}. However, most of these publications have a highly algorithmic perspective, focusing on a restricted class of models and optimality criteria. Moreover, many of the proposed methods use optimality criteria (e.g., A-optimality) with well-known deficiencies, such as lack of invariance to the measurement-scale of the data and parameterisation of the model. A unified theory of optimal subsampling design is still lacking.

We present a theory of optimal design for general data subsampling problems, including finite population inference, parametric density estimation, and regression modelling using quasi-likelihood methods. We derive optimality conditions for a broad class of optimality criteria, including A-, D-, E-, L-, and Kiefer's $\Phi_q$-optimality criterion \citep{Kiefer1974}. Algorithms to find optimal sampling schemes are presented for Poisson sampling and multinomial sampling designs. We also study optimal design from a distance-minimising perspective, and establish equivalence to traditional optimality criteria. This naturally leads us to a novel class of linear optimality criteria with good theoretical and practical properties, including computational tractability and invariance under affine transformations of the data and re-parameterisation of the model. The presented methodology and algorithms are illustrated in an application in the traffic safety domain.

We start with a brief review of some standard methods in unequal probability sampling and optimal design in Section \ref{sec:preliminaries}. A general theory of optimal design for data subsampling problems is presented in Section \ref{sec:optimal_subsampling_designs}, including algorithms for finding optimal sampling schemes. We discuss optimal design from a distance-minimising perspective in Section \ref{sec:dist}, and present optimal designs for some common statistical distance functions. Comments on the implementation of optimal subsampling methods in practice are provided in Section \ref{sec:practical_implementation}. Examples and experiments are presented in Section \ref{sec:application}. We refer to Appendix \ref{sec:proofs} for proofs.

\section{Preliminaries}
\label{sec:preliminaries}

Consider a class of experiments $\Xi$ and corresponding consistent estimators $\hat{\boldsymbol{\theta}}_{\xi}, \xi \in \Xi$, for an unknown parameter $\boldsymbol{\theta}^*$. The aim of optimal design is to find an experiment $\xi \in \Xi$ that minimises some suitable function $\Phi$ of the covariance matrix of the estimator $\hat{\boldsymbol{\theta}}_{\xi}$. For instance, $\Phi$ may be the sum or product of the eigenvalues of its matrix argument, corresponding to A- or D-optimality \citep{Atkinson1992}, or some other measure of "size" of a matrix.

In the context of data subsampling, the experiment is determined by the choice of sampling design and sampling scheme $\boldsymbol{\mu} = (\mu_1, \ldots, \mu_N)$. For the estimation problem outlined in Section \ref{sec:introduction}, we wish to find a sampling scheme $\boldsymbol{\mu}$ that minimises $\Phi(\COV({\hat{\boldsymbol{\theta}}_{\!\boldsymbol{\mu}}}))$ for some suitable family of sampling designs and objective function $\Phi:\mathbb{R}^{p \times p}\rightarrow \mathbb{R}$. Some common unequal probability sampling designs are presented in Section \ref{sec:unequal_probability_sampling}. Expressions for the approximate covariance matrix of the estimator $\hat{\boldsymbol{\theta}}_{\!\boldsymbol{\mu}}$ are provided in Section \ref{sec:approximate_covariance_matrix}, and a brief review of optimal design in Section \ref{sec:optimality_criteria}.

\subsection{Unequal probability sampling designs} 
\label{sec:unequal_probability_sampling}

We consider the situation where individual elements $i \in \mathcal D$ are selected according to an unequal probability sampling design, i.e., by a random mechanism where each member $i \in \mathcal D$ has a strictly positive and possibly unique selection probability. Following the notation in Section \ref{sec:introduction}, we let $S_i$ be the number of times an element $i \in \mathcal D$ is selected by the sampling mechanism, where sampling may be with or without replacement, and $\mu_i$ be the corresponding expected number of selections. We let $n$ denote the expected size of the subsample, and $\mathcal M_n$ the corresponding domain of $\boldsymbol{\mu} = (\mu_1, \ldots, \mu_N)$, i.e., the set of feasible values of the sampling scheme $\boldsymbol{\mu}$ within a specified family of sampling designs of (expected) size $n$. We assume that sampling is conducted according to one of the following families of sampling designs:
\begin{enumerate}[label=\roman*)]
    \item Poisson sampling with replacement (PO-WR): $S_1, \ldots, S_N$ are independent with $S_i \sim \mathrm{Poisson}(\mu_i)$, $\mu_i > 0$. The sample size $\sum_{i\in \mathcal D} S_i$ is random, with expectation $\E[\sum_{i\in \mathcal D} S_i] = \sum_{i\in \mathcal D}\mu_i = n$. The corresponding domain $\mathcal M_n $of $\boldsymbol{\mu}$ is given by $\mathcal M_n = \{\boldsymbol{\mu} \in \mathbb{R}^N: \mu_i > 0 \text{ for all } i \in \mathcal D \text{ and } \sum_{i\in \mathcal D} \mu_i = n\}$. 
    \item Poisson sampling without replacement (PO-WOR): $S_1, \ldots, S_N$ are independent with $S_i \sim \mathrm{Bernoulli}(\mu_i)$, $\mu_i \in (0, 1]$. The sample size $\sum_{i\in \mathcal D} S_i$ is random, with expectation $\E[\sum_{i\in \mathcal D} S_i] = \sum_{i\in \mathcal D}\mu_i = n$. The corresponding domain $\mathcal M_n $of $\boldsymbol{\mu}$ is given by $\mathcal M_n = \{\boldsymbol{\mu} \in \mathbb{R}^N: \mu_i \in (0, 1] \text{ for all } i \in \mathcal D \text{ and } \sum_{i\in \mathcal D} \mu_i = n\}$. 
    \item Multinomial sampling (MULTI): $(S_1, \ldots, S_N) \sim \mathrm{Multinomial}(n, \boldsymbol{\mu}/n)$, $n \in \mathbb{N}, \mu_i>0, \sum_{i \in \mathcal D} \mu_i = n$. Sampling is done with replacement and the sample size is fixed, i.e., $\sum_{i\in \mathcal D}S_i = n$. The corresponding domain $\mathcal M_n $ of $\boldsymbol{\mu}$ is given by $\mathcal M_n = \{\boldsymbol{\mu} \in \mathbb{R}^N: \mu_i > 0 \text{ for all } i \in \mathcal D \text{ and } \sum_{i \in \mathcal D} \mu_i = n\}$. 
\end{enumerate}

For a given size $n$, the Poisson and multinomial sampling designs are uniquely determined by the mean vector $\boldsymbol{\mu}$. We say that such a design, for a given size $n$, is indexed by the sampling scheme $\boldsymbol{\mu}$.

Methods also exist to select a fixed number of elements without replacement and with fixed selection probabilities, for instance using conditional Poisson sampling \citep{Hajek1981, Tille2006}. This method is, however, both computationally and analytically intractable, and will therefore not be considered in this paper. Additional details may be found in, e.g., \citet{Tille2006} and \citet{Fuller2009}.

\subsection{Covariance matrix of the Hansen-Hurwitz empirical risk minimiser} 
\label{sec:approximate_covariance_matrix}

\citet{Binder1983} showed that under suitable regularity conditions the distribution of the estimator \eqref{eq:thetahat} with respect to the sampling mechanism is approximately Gaussian with mean
\begin{equation}
    \label{eq:unbiased}
    \E[\hat{\boldsymbol{\theta}}_{\!\boldsymbol{\mu}}] = \boldsymbol{\theta}_0 + o(n^{-1/2}), 
\end{equation}
and covariance matrix
\begin{align}
    \label{eq:approximate_covariance}
    \COV(\hat{\boldsymbol{\theta}}_{\!\boldsymbol{\mu}} - \boldsymbol{\theta}_0) =  
    \boldsymbol{\Gamma}(\boldsymbol{\mu}; \boldsymbol{\theta}_0) + o(n^{-1}), \quad
    \boldsymbol{\Gamma}(\boldsymbol{\mu};\boldsymbol{\theta}_0) = 
    \mathbf{H}(\boldsymbol \theta_0)^{-1} 
    \mathbf{V}(\boldsymbol{\mu}; \boldsymbol{\theta}_0)
    \mathbf{H}(\boldsymbol \theta_0)^{-1}.
\end{align}  
Here $o(n^{-1/2})$ and $o(n^{-1})$ are interpreted elementwise and
$
\mathbf{H}(\boldsymbol \theta_0) = \frac{\partial^2\ell_0(\boldsymbol{\theta})}{\partial \boldsymbol \theta \partial \boldsymbol \theta^{\mathsf{T}}}\bigr\rvert_{\boldsymbol{\theta} = \boldsymbol{\theta}_0}
$ 
is the Hessian of the full-data empirical risk function \eqref{eq:ell0} at $\boldsymbol{\theta} = \boldsymbol{\theta}_0$.
\begin{equation}
    \label{eq:vmat}
    \mathbf{V}(\boldsymbol{\mu}; \boldsymbol{\theta}_0) = \COV\left(\nabla_{\!\boldsymbol{\theta}}\hat{\ell}_{\boldsymbol{\mu}}(\boldsymbol{\theta})\bigr\rvert_{\boldsymbol{\theta} = \boldsymbol{\theta}_0}\right)
= \sum_{i,j \in \mathcal D} \frac{\Cov(S_i, S_j)}{\mu_i\mu_j}\boldsymbol{\psi}_i(\boldsymbol{\theta}_0)\boldsymbol{\psi}_j(\boldsymbol{\theta}_0)\T
\end{equation}
is the covariance matrix of the gradient $\nabla_{\!\boldsymbol{\theta}}\hat{\ell}_{\boldsymbol{\mu}}(\boldsymbol{\theta})$ with respect to the sampling mechanism, evaluated at $\boldsymbol{\theta} = \boldsymbol{\theta}_0$, and 
$
\boldsymbol{\psi}_i(\boldsymbol{\theta}) = \nabla_{\!\boldsymbol{\theta}}\ell_i(\boldsymbol{\theta})$. 
We refer to \citet{Binder1983} and \citet{Fuller2009} for further details.

It follows from the properties of the sampling designs described in Section \ref{sec:unequal_probability_sampling}, that the matrix $\mathbf{V}(\boldsymbol{\mu}; \boldsymbol{\theta}_0)$ can be simplified to
\begin{equation}
    \label{eq:Vmat}
    \mathbf{V}(\boldsymbol{\mu}; \boldsymbol{\theta}_0) = 
    \begin{cases}
        \sum_{i \in \mathcal D} \mu_i^{-1}\boldsymbol{\psi}_i(\boldsymbol{\theta}_0)\boldsymbol{\psi}_i(\boldsymbol{\theta}_0)\T, & \text{for PO-WR and MULTI designs, and} \\
        \sum_{i \in \mathcal D} (\mu_i^{-1} - 1)\boldsymbol{\psi}_i(\boldsymbol{\theta}_0)\boldsymbol{\psi}_i(\boldsymbol{\theta}_0)\T, & \text{for PO-WOR}.
    \end{cases}    
\end{equation}
See, e.g., \citet{Tille2006}. To obtain the above result for the multinomial sampling design, we have also used \eqref{eq:grads}.

\subsection{Optimal design} 
\label{sec:optimality_criteria}

For an unknown parameter $\boldsymbol{\theta}^*$, consider a class of experiments $\Xi$ and corresponding consistent estimators $\hat{\boldsymbol{\theta}}_{\xi}, \xi \in \Xi$, with unequal covariance matrices $\boldsymbol{\Gamma}_{\xi}$. Ideally, we would like to find an experiment $\xi^* \in \Xi$ such that $\boldsymbol{\Gamma}_{\xi} - \boldsymbol{\Gamma}_{\xi^*}$ is positive semi-definite for all $\xi \in \Xi$. Such universal optimality, however, is not possible to achieve in general. Hence, instead we consider a function $\Phi: \boldsymbol{S}_{\!+}^{p \times p}\to \mathbb{R}$ on the set of real, symmetric, positive semi-definite $p \times p$ matrices, for which a minimiser $\xi^* \in \Xi$ is sought. For $\Phi$ to be a meaningful measure of optimality we require the function to be monotone for Loewner's ordering, i.e., that 
\begin{equation}
    \label{eq:loewner}
    \Phi(\mathbf{U}) \ge \Phi(\mathbf{V}) \text{ for all } \mathbf{U}, \mathbf{V} \in \boldsymbol{S}_{\!+}^{p \times p} \text{ such that } \mathbf{U} \ge \mathbf{V} , 
\end{equation}
with $\mathbf{U} \ge \mathbf{V}$ meaning that $\mathbf{U} - \mathbf{V}$ is positive semi-definite \citep{Pukelsheim1993}.

Some popular optimality criteria are defined and summarised in Table \ref{tab:OED}. These include the D-optimality criterion (minimise the determinant of the covariance matrix), the E-optimality criterion (minimise the largest eigenvalue of the covariance matrix) and the L-optimality criterion (minimise the average variance of a collection of linear combinations $\mathbf{L}\T\hat{\boldsymbol{\theta}}_{\xi}$). Two important special cases of the L-optimality criterion are the A-optimality criterion (minimise the average variance) and c-optimality criterion (minimise the variance of a linear combination $\mathbf{c}\T\hat{\boldsymbol{\theta}}_{\xi}$), obtained with $\mathbf{L} = \mathbf{I}_{p \times p}$ and $\mathbf{L} = \mathbf{c}$ for some $p \times 1$ vector $\mathbf{c}$, respectively \citep{Silvey1980, Atkinson1992}. Included in Table \ref{tab:OED} is also the $\Phi_{q}$- and $\Phi_{q,{\mathbf{A}}}$-optimality criteria, which encompass all other optimality criteria in this table. In particular, $\Phi_q$-optimality coincides with D-optimality when $q = 0$, A-optimality when $q=1$, and E-optimality when $q = \infty$ \citep{Kiefer1974}. Hence, $\Phi_q$-optimality can be used to interpolate between A-, D- and E-optimality.

The A-, D- and E- optimality criteria have a simple geometric interpretation as follows. Consider the random set $\mathcal C(\hat{\boldsymbol{\theta}}_{\xi}) := \{\boldsymbol{\theta} \in \mathbb{R}^p: (\boldsymbol{\theta} - \hat{\boldsymbol{\theta}}_{\xi})\T\boldsymbol{\Gamma}_{\xi}^{-1}(\boldsymbol{\theta} - \hat{\boldsymbol{\theta}}_{\xi}) \le \chi^2_{p,\alpha}\}$, where $\chi^2_{p, \alpha}$ is the $\alpha$-quantile of a $\chi^2$-distribution with $p$ degrees of freedom. For an (approximately) normally distributed estimator $\hat{\boldsymbol{\theta}}_{\xi}$, this defines an (approximate) $100\times(1-\alpha)$\% ellipsoidal confidence set for $\boldsymbol{\theta}^*$ in $\mathbb{R}^p$. D-optimality minimises the volume of this confidence ellipsoid over the class of experiments $\Xi$. E-optimality minimises the length of its longest axis, and A-optimality the length of the diagonal of the minimal bounding box (parallelepiped) around the confidence ellipsoid \citep{Pronzato2013}.

Another popular optimality criterion is the V-optimality criterion, which minimises the average prediction variance with respect to some measure $\nu(\boldsymbol{x})$ on the design space $\mathcal X$ \citep{Welch1984}. This is a linear optimality criterion and hence is covered by the L-optimality criterion for a matrix $\mathbf{L}$ such that $\mathbf{L}\mathbf{L}\T = \int_{\mathcal{X}}\boldsymbol{\varphi}(\boldsymbol{x})\boldsymbol{\varphi}(\boldsymbol{x})\T d\nu(\boldsymbol{x})$ (Table \ref{tab:OED}) \citep{Atkinson1992}. A natural choice for the measure $\nu(\boldsymbol{x})$ in data subsampling problems is the empirical measure on $\{\boldsymbol{x}_i\}_{i \in \mathcal D}$.

A property that is often desirable for an optimal design, is invariance under a non-singular affine transformation of the data and under a re-parameterisation of the model. That is, the optimal design and the statistical properties of the resulting estimator should not depend on the choice of parameterisation, nor on the scaling or coding of the data prior to modelling. The most common example of a transformation- and parameterisation invariant optimality criterion is the D-optimality criterion. In contrast, the A- and E-optimality criteria are sensitive to changes in the parameterisation or data, and hence lack such invariance properties \citep{Atkinson1992}. An L-optimal design may or may not be parameterisation- and transformation-invariant, depending on whether or not the coefficient matrix $\mathbf{L}$ of the L-optimality criterion is adapted to the parameterisation of the problem and scaling of the data. Some examples of transformation- and parameterisation-invariant linear optimality criteria will be discussed in Section \ref{sec:another_look_at_invariance}.

For further details, were refer to \citet{Silvey1980}, \citet{Atkinson1992} and \citet{Pukelsheim1993} and \citet{Pronzato2013}.

\begin{table}[htb!]
    \centering
    \caption{Definition of some common optimality criteria in optimal design. $\Phi$ is a real-valued function on the set of real, symmetric, positive semi-definite $p \times p$ matrices, $\boldsymbol{\Gamma}$ the $p \times p$ covariance matrix of an estimator $\hat{\boldsymbol{\theta}} = (\hat{\theta}_1, \ldots, \hat{\theta}_p)$, $\lambda_{\max}(\boldsymbol{\Gamma})$ the largest eigenvalue of $\boldsymbol{\Gamma}$, 
    $\mathbf{c}$ a non-zero $p \times 1$ vector,
    $\mathbf{A} = \mathbf{L}$ a non-zero $p \times m$ matrix with columns $\mathbf{a}_1, \ldots, \mathbf{a}_m$, and $\mathbf{I}_{p \times p}$ the $p \times p$ identity matrix. $\mathcal X$ is the set of possible values for the predictors $\boldsymbol{x}$ and $\boldsymbol{\varphi}: \mathcal X\to\mathbb{R}^p$ a feature map of the data.}
    \label{tab:OED}
    \resizebox{\textwidth}{!}{\begin{tabular}{m{3.1cm}m{6.0cm}m{6.25cm}}
        \textbf{Optimality criterion} & \textbf{Description} & \textbf{Objective function} $\Phi(\boldsymbol{\Gamma})$ 
        \\ \hline
        
        A-optimality & 
        Minimise average variance, \newline minimise trace of covariance matrix, \newline minimise sum of eigenvalues. & 
        $\frac{1}{p}\sum_{i = 1}^p \Var(\hat{\theta}_i) = $
        $\frac{1}{p}\tr(\boldsymbol{\Gamma})$
        \arrayrulecolor{lightgray} \\ \hline
        
        c-optimality & 
        Minimise variance of a linear combination or contrast $\mathbf{c}\T\hat{\boldsymbol{\theta}}$. & 
        $\Var(\mathbf{c}\T\hat{\boldsymbol{\theta}}) = $
        $\mathbf{c}\T\boldsymbol{\Gamma}\mathbf{c} = \tr(\boldsymbol{\Gamma}\mathbf{c}\mathbf{c}\T)$
        \arrayrulecolor{lightgray} \\ \hline
        
        D-optimality & 
        Minimise generalised variance, \newline minimise determinant of covariance matrix, \newline minimise product of eigenvalues. &
        $\det(\boldsymbol{\Gamma})^{1/p}$ or $\log \det(\boldsymbol{\Gamma})$
        \arrayrulecolor{lightgray} \\ \hline
        
        D$_{\mathrm{A}}$-optimality & 
        Minimise generalised variance for subset of parameters, collection of linear combinations, or contrasts $\mathbf{A}\T\boldsymbol{\theta}$. &
        $\det(\mathbf{A}\T\boldsymbol{\Gamma}\mathbf{A})$
        \arrayrulecolor{lightgray} \\ \hline
        
        E-optimality & 
        Minimise maximal eigenvalue, \newline
        minimise variance along the direction of largest uncertainty. &
        $\lambda_{\max}(\boldsymbol{\Gamma})$
        \arrayrulecolor{lightgray} \\ \hline
        
        L-optimality &
        Minimise average variance of a collection of linear combinations or contrasts $\mathbf{L}\T\boldsymbol{\theta}$. \newline
        $\mathbf{L} = \boldsymbol{c}  \hspace{0.5cm} \Leftrightarrow$ c-optimality \newline
        $\mathbf{L} = \mathbf{I}_{p \times p} \Leftrightarrow$ A-optimality
        &
        $\frac{1}{m}\sum_{i = 1}^m\Var(\mathbf{a}_i\T\hat{\boldsymbol{\theta}}) = $
        $\frac{1}{m}\tr(\boldsymbol{\Gamma}\mathbf{L}\mathbf{L}\T)$ 
        \arrayrulecolor{lightgray} \\ \hline

        V-optimality &
        Minimise average prediction variance with respect to a measure $d\nu(\boldsymbol{x})$ on $\mathcal X$, assuming a linear model $\hat y = \boldsymbol{\varphi}(\boldsymbol{x})\T\hat{\boldsymbol{\theta}}$. &
        $\int_{\mathcal X}\Var(\boldsymbol{\varphi}(\boldsymbol{x})\T\hat{\boldsymbol{\theta}})d\nu(\boldsymbol{x})$ $=$ \vfill $\tr\left(\boldsymbol{\Gamma}\int_{\mathcal{X}}\boldsymbol{\varphi}(\boldsymbol{x})\boldsymbol{\varphi}(\boldsymbol{x})\T d\nu(\boldsymbol{x}) \right)$
        \arrayrulecolor{lightgray} \\ \hline
  
        $\Phi_{q, \mathbf{A}}$-optimality \newline 
        $q \in [0, \infty]$ & 
        $\Phi_0 \hspace{0.30cm} \Leftrightarrow$ D-optimality \newline 
        $\Phi_{0, \mathbf{A}}  \Leftrightarrow$ D$_{\mathrm{A}}$-optimality \newline 
        $\Phi_1 \hspace{0.32cm} \Leftrightarrow$ A-optimality \newline
        $\Phi_{1,\mathbf{L}} \hspace{0.05cm} \Leftrightarrow$ L-optimality \newline
        $\Phi_{\infty} \hspace{0.18cm} \Leftrightarrow$ E-optimality &
        $\Phi_{q,\mathbf{A}}(\boldsymbol{\Gamma}) \hspace{0.13cm} = \frac{1}{m}\tr[(\mathbf{A}\T\boldsymbol{\Gamma}\mathbf{A})^{q}]^{1/q}, q \in (0, \infty)$
        \newline
        $\Phi_{0,\mathbf{A}}(\boldsymbol{\Gamma}) \hspace{0.14cm} = \lim_{q \downarrow 0}\Phi_{q,\mathbf{A}}(\boldsymbol{\Gamma})$ \newline
        $\Phi_{\infty,\mathbf{A}}(\boldsymbol{\Gamma}) = \lim_{q \uparrow \infty}\Phi_{q,\mathbf{A}}(\boldsymbol{\Gamma})$ \newline 
        $\Phi_{q}(\boldsymbol{\Gamma}) \hspace{0.46cm} = \Phi_{q,\mathbf{A}}(\boldsymbol{\Gamma}), \mathbf{A} = \mathbf{I}_{p \times p}$ 
        \arrayrulecolor{black} \\ \hline
    \end{tabular}}
\end{table}

\section{Optimal subsampling designs}
\label{sec:optimal_subsampling_designs}

In this section we present optimal sampling schemes for a general class of optimality criteria, under an assumption of differentiability. $\Phi$-optimality is defined  in Section \ref{sec:phi_optimality}, where we also present three important lemmas. Optimality criteria for Poisson and multinomial sampling designs are presented in Section \ref{sec:optimality_conditions}, and algorithms for finding optimal sampling schemes in Section \ref{sec:optimal_sampling_schemes}.

First we note that the approximate covariance matrix $\boldsymbol{\Gamma}(\boldsymbol{\mu}; \boldsymbol{\theta}_0)$ of the estimator $\hat{\boldsymbol{\theta}}_{\!\boldsymbol{\mu}}$, as given in \eqref{eq:approximate_covariance}, generally depends on the full data $\{(\boldsymbol{x}_i, \boldsymbol{y}_i)\}_{i \in \mathcal D}$ and full-data parameter $\boldsymbol{\theta}_0$. Clearly, subsampling would not be needed if such information were available at the design stage. This is a general problem in optimal design, however, and not specific to our setup, and hence not a major limitation of the theory we present. We will proceed in this section and Section \ref{sec:dist} as if such information is available, keeping in mind that the resulting theoretically optimal designs can generally not be found in practice. We refer to Section \ref{sec:practical_implementation} for a discussion on the implementation of optimal subsampling designs in practice.

Throughout we assume regularity conditions such that \eqref{eq:approximate_covariance} holds, and that $\mathbf{H}(\boldsymbol{\theta}_0)$ is of full rank. All vectors are assumed to be column vectors, unless otherwise stated. We let $||\mathbf{u}||_2^2 = \mathbf{u}\T\mathbf{u}$ denote the Euclidean norm of a vector $\mathbf{u}$. Also recall that $\psi_i(\boldsymbol{\theta}) = \nabla_{\!\boldsymbol{\theta}}\ell_i(\boldsymbol{\theta})$.

\subsection{Optimality criteria}
\label{sec:phi_optimality}

By an optimal sampling scheme $\boldsymbol{\mu}^*$, we mean the following:

\begin{definition}[$\Phi$-optimality] 
    \label{def:phi_optimality} 
    Consider a function $\Phi: \boldsymbol{S}_{\!+}^{p \times p}\to \mathbb{R}$ that is monotone for Loewner's ordering, i.e., such that \eqref{eq:loewner} holds. Also consider a family of unequal probability sampling designs (e.g., PO-WR, PO-WOR or MULTI) indexed by the sampling scheme $\boldsymbol{\mu}$. Let the expected size $\E[\sum_{i\in\mathcal D}S_i] = n$ be fixed, and let $\mathcal M_n$ denote the corresponding domain of $\boldsymbol{\mu}$. We say that a sampling scheme $\boldsymbol{\mu}^*$ is $\Phi$-optimal if
    \begin{equation*}
        \boldsymbol{\mu}^* = \argmin_{\boldsymbol{\mu} \in \mathcal M_n} \Phi(\boldsymbol{\Gamma}(\boldsymbol{\mu}; \boldsymbol{\theta}_0)),    
    \end{equation*}
    where $\boldsymbol{\Gamma}(\boldsymbol{\mu}; \boldsymbol{\theta}_0)$ is the approximate covariance matrix of $\hat{\boldsymbol{\theta}}_{\!\boldsymbol{\mu}}$, as given in \eqref{eq:approximate_covariance}.
\end{definition}

Finding a $\Phi$-optimal sampling scheme reduces to a non-linear, possibly non-convex, restricted optimisation problem over an \mathprefix{(N-1)}-dimensional hyperplane in $\mathbb{R}^{N}$. While this problem may be addressed by numerical optimisation methods when $N$ is small, this is generally not a viable option for large datasets. We therefore need a theory of optimal design that can be used to devise efficient algorithms for finding optimal sampling schemes when $N$ is large. To make the problem tractable, we will restrict ourselves to optimality criteria $\Phi(\boldsymbol{\Gamma}(\boldsymbol{\mu}; \boldsymbol{\theta}_0))$ that are differentiable with respect to $\boldsymbol{\mu}$ in a neighbourhood of its optimum $\boldsymbol{\mu}^*$. Three important lemmas are provided below.

\begin{lemma}[The chain rule]  
    \label{lemma:chain_rule}
    Consider a function $\Phi: \boldsymbol{S}_{\!+}^{p \times p} \to \mathbb{R}^p$, and assume that $\Phi(\boldsymbol{\Gamma}(\boldsymbol{\mu}; \boldsymbol{\theta}_0))$ is differentiable with respect to $\boldsymbol{\mu}$ in a neighbourhood of some point $\boldsymbol{\mu}^*$. The partial derivative of $\Phi(\boldsymbol{\Gamma}(\boldsymbol{\mu}; \boldsymbol{\theta}_0))$ with respect to $\mu_i$ is then given by
    \begin{equation}
        \label{eq:chain_rule}
        \frac{\partial \Phi(\boldsymbol{\Gamma}(\boldsymbol{\mu}; \boldsymbol{\theta}_0))}{\partial \mu_i} = \tr\left(\boldsymbol{\phi}(\boldsymbol{\Gamma}(\boldsymbol{\mu}; \boldsymbol{\theta}_0))\frac{\partial \boldsymbol{\Gamma}(\boldsymbol{\mu}; \boldsymbol{\theta}_0)}{\partial \mu_i}\right),    
    \end{equation}
    where $\boldsymbol{\phi}(\mathbf{U}) = \frac{\partial \Phi(\mathbf{U})}{\partial \mathbf{U}}$ is the $p \times p$ matrix derivative of $\Phi$ with respect to its matrix argument, and $\frac{\partial \boldsymbol{\Gamma}(\boldsymbol{\mu}; \boldsymbol{\theta}_0)}{\partial \mu_i}$ is the elementwise derivative of $\boldsymbol{\Gamma}(\boldsymbol{\mu}; \boldsymbol{\theta}_0)$ with respect to $\mu_i$.

    Assume further that 
    \begin{enumerate}[label = \roman*)]
        \item $\boldsymbol{\Gamma}(\boldsymbol{\mu}; \boldsymbol{\theta}_0)$ decreases monotonically with $\mu_1, \ldots, \mu_N$ in the Loewner order sense, i.e., $\boldsymbol{\Gamma}(\boldsymbol{\mu}_1; \boldsymbol{\theta}_0) - \boldsymbol{\Gamma}(\boldsymbol{\mu}_2; \boldsymbol{\theta}_0)$ is positive semi-definite for every pair of vectors $\boldsymbol{\mu}_1, \boldsymbol{\mu}_2 \in \mathbb{R}^N_{>0}$ such that $\boldsymbol{\mu}_1 \le \boldsymbol{\mu}_2$ (elementwise), and
        \item $\Phi$ is monotone for Loewner's ordering, i.e., that \eqref{eq:loewner} holds.
    \end{enumerate}
    Then the matrix $\boldsymbol{\phi}(\boldsymbol{\Gamma}(\boldsymbol{\mu}; \boldsymbol{\theta}_0))$ is positive semi-definite and there exists a real matrix $\mathbf{L}(\boldsymbol{\mu}; \boldsymbol{\theta}_0)$ such that $\mathbf{L}(\boldsymbol{\mu}; \boldsymbol{\theta}_0)\mathbf{L}(\boldsymbol{\mu}; \boldsymbol{\theta}_0)\T = \boldsymbol{\phi}(\boldsymbol{\Gamma}(\boldsymbol{\mu}; \boldsymbol{\theta}_0))$.
\end{lemma}

The first part of Lemma \ref{lemma:chain_rule} follows by the chain rule in matrix differential calculus and the symmetry of $\boldsymbol{\Gamma}(\boldsymbol{\mu}; \boldsymbol{\theta}_0)$, and the second by the monotonicity assumptions on $\boldsymbol{\Gamma}(\boldsymbol{\mu}; \boldsymbol{\theta}_0)$ and $\Phi$. The matrix $\mathbf{L}(\boldsymbol{\mu}; \boldsymbol{\theta}_0)$ may, e.g.,  be obtained as the matrix square root of $\boldsymbol{\phi}(\boldsymbol{\Gamma}(\boldsymbol{\mu}; \boldsymbol{\theta}_0))$, or by the Cholesky decomposition when $\boldsymbol{\phi}(\boldsymbol{\Gamma}(\boldsymbol{\mu}; \boldsymbol{\theta}_0))$ is of full rank. Some examples are provided in Lemma \ref{lemma:differentiability}.

\begin{lemma}[$\boldsymbol{\mu}$-differentiable $\Phi$-optimality criteria]  
    \label{lemma:differentiability}
    Consider a PO-WR, PO-WOR or MULTI design, and assume that $\mathbf{H}(\boldsymbol{\theta}_0)$ is of full rank. Let $\mathbf{c}$ be a non-zero $p \times 1$ vector, $\mathbf{L}$ a non-zero $p \times m$ matrix, $\lambda_{\max}(\boldsymbol{\Gamma}(\boldsymbol{\mu}; \boldsymbol{\theta}_0))$ the maximal eigenvalue of $\boldsymbol{\Gamma}(\boldsymbol{\mu}; \boldsymbol{\theta}_0)$, and $\boldsymbol{v}_{\!\boldsymbol{\mu}}$ a corresponding eigenvector. Let $\boldsymbol{\phi}(\boldsymbol{\Gamma}(\boldsymbol{\mu}; \boldsymbol{\theta}_0))$ be defined as in Lemma \ref{lemma:chain_rule}. Then the following holds:
    \begin{enumerate}[label = \alph*)]
        \item $\boldsymbol{\Gamma}(\boldsymbol{\mu}; \boldsymbol{\theta}_0)$ is differentiable with respect to $\boldsymbol{\mu}$ and $\frac{\partial \boldsymbol{\Gamma}(\boldsymbol{\mu}; \boldsymbol{\theta}_0)}{\partial \mu_i} = -\mu_i^{-2}\mathbf{H}(\boldsymbol{\theta}_0)^{-1}\boldsymbol{\psi}_i(\boldsymbol{\theta}_0)\boldsymbol{\psi}_i(\boldsymbol{\theta}_0)\T\mathbf{H}(\boldsymbol{\theta}_0)^{-1}$, provided that $\mu_i > 0$.
        \label{prop:differentiability_a}
        \item The D-optimality objective function $\Phi(\boldsymbol{\Gamma}(\boldsymbol{\mu}; \boldsymbol{\theta}_0)) = \log \det(\boldsymbol{\Gamma}(\boldsymbol{\mu}; \boldsymbol{\theta}_0))$ is differentiable with respect to $\boldsymbol{\mu}$ and $\boldsymbol{\phi}(\boldsymbol{\Gamma}(\boldsymbol{\mu}; \boldsymbol{\theta}_0)) = \boldsymbol{\Gamma}(\boldsymbol{\mu}; \boldsymbol{\theta}_0)^{-1}$, 
        provided that $\boldsymbol{\Gamma}(\boldsymbol{\mu}; \boldsymbol{\theta}_0)$ is of full rank. 
        \label{prop:differentiability_b}
        \item The E-optimality objective function $\Phi(\boldsymbol{\Gamma}(\boldsymbol{\mu}; \boldsymbol{\theta}_0)) = \lambda_{\max}(\boldsymbol{\Gamma}(\boldsymbol{\mu}; \boldsymbol{\theta}_0))$ is differentiable with respect to $\boldsymbol{\mu}$ and $\boldsymbol{\phi}(\boldsymbol{\Gamma}(\boldsymbol{\mu}; \boldsymbol{\theta}_0)) = \boldsymbol{v}_{\!\boldsymbol{\mu}}\boldsymbol{v}_{\!\boldsymbol{\mu}}\T$, provided that $\lambda_{\max}(\boldsymbol{\Gamma}(\boldsymbol{\mu}; \boldsymbol{\theta}_0))$ has multiplicity $1$.
        \label{prop:differentiability_c}
        \item The L-optimality objective function $\Phi(\boldsymbol{\Gamma}(\boldsymbol{\mu}; \boldsymbol{\theta}_0)) = \tr(\boldsymbol{\Gamma}(\boldsymbol{\mu}; \boldsymbol{\theta}_0)\mathbf{L}\mathbf{L}\T)$ is differentiable with respect to $\boldsymbol{\mu}$, and $\boldsymbol{\phi}(\boldsymbol{\Gamma}(\boldsymbol{\mu}; \boldsymbol{\theta}_0)) = \mathbf{L}\mathbf{L}\T$. In particular, this holds for A-optimality with $\mathbf{L} = \mathbf{I}_{p \times p}$ and  c-optimality with $\mathbf{L}  = \mathbf{c}$.
        \label{prop:differentiability_d}
        \item The $\Phi_q$-optimality objective function $\Phi(\boldsymbol{\Gamma}(\boldsymbol{\mu}; \boldsymbol{\theta}_0)) = \tr(\boldsymbol{\Gamma}(\boldsymbol{\mu}; \boldsymbol{\theta}_0)^{q})^{1/q}$ is differentiable with respect to $\boldsymbol{\mu}$ for $q \in (0, \infty)$ and $\boldsymbol{\phi}(\boldsymbol{\Gamma}(\boldsymbol{\mu}; \boldsymbol{\theta}_0)) = \tr(\boldsymbol{\Gamma}(\boldsymbol{\mu}; \boldsymbol{\theta}_0)^q)^{1/q-1}\boldsymbol{\Gamma}(\boldsymbol{\mu}; \boldsymbol{\theta}_0)^{q-1}$, provided that $\boldsymbol{\Gamma}(\boldsymbol{\mu}; \boldsymbol{\theta}_0)$ is of full rank.
        \label{prop:differentiability_e}
        \end{enumerate}
\end{lemma}

Combining the results of Lemma \ref{lemma:chain_rule} and \ref{lemma:differentiability}, we obtain the following:

\begin{lemma}[Partial derivatives of $\Phi(\boldsymbol{\Gamma}(\boldsymbol{\mu}; \boldsymbol{\theta}_0))$]  
    \label{lemma:partial_derivatives}
    Consider a PO-WR, PO-WOR or MULTI design. Also consider a function $\Phi: \boldsymbol{S}_{\!+}^{p \times p} \to \mathbb R$ such that $\Phi$ is monotone for Loewner's ordering. Assume that $\mathbf{H}(\boldsymbol{\theta}_0)$ is of full rank, and that $\Phi(\boldsymbol{\Gamma}(\boldsymbol{\mu}; \boldsymbol{\theta}_0))$ is differentiable with respect to $\boldsymbol{\mu}$ in a neighbourhood of some point $\boldsymbol{\mu}^*$. Let $\mathbf{L}(\boldsymbol{\mu}; \boldsymbol{\theta}_0)$ be defined as in Lemma \ref{lemma:chain_rule}. Then 
    \begin{equation*}
        \frac{\partial \Phi(\boldsymbol{\Gamma}(\boldsymbol{\mu}; \boldsymbol{\theta}_0))}{\partial \mu_i} = -\mu_i^{-2}
        \bigr\rvert\bigr\rvert\mathbf{L}(\boldsymbol{\mu}; \boldsymbol{\theta}_0)\T
        \mathbf{H}(\boldsymbol{\theta}_0)^{-1}
        \boldsymbol{\psi}_i(\boldsymbol{\theta}_0)\bigr\rvert\bigr\rvert^2_{2}.
    \end{equation*}
\end{lemma}

\subsection{Optimality conditions}
\label{sec:optimality_conditions}

Using results of Lemma \ref{lemma:chain_rule}–\ref{lemma:partial_derivatives}, in Proposition \ref{prop:optimality_conditions} we present optimality conditions for Poisson and multinomial sampling designs with respect to a $\Phi$-optimality criterion under an assumption of differentiability.

\begin{proposition}[$\Phi$-optimality conditions] 
    \label{prop:optimality_conditions}
    Consider the family of PO-WR, PO-WOR or MULTI designs of (expected) size $n$. Also consider a function $\Phi: \boldsymbol{S}_{\!+}^{p \times p} \to \mathbb R$ such that $\Phi$ is monotone for Loewner's ordering. Assume that $\mathbf{H}(\boldsymbol{\theta}_0)$ is of full rank, and that $\Phi(\boldsymbol{\Gamma}(\boldsymbol{\mu}; \boldsymbol{\theta}_0))$ is differentiable with respect to $\boldsymbol{\mu}$ in a neighbourhood of some point $\boldsymbol{\mu}^*$. Let $\mathbf{L}(\boldsymbol{\mu}; \boldsymbol{\theta}_0)$ be defined according to Lemma \ref{lemma:chain_rule}, and
    \begin{equation}
        \label{eq:ci} 
        c_i = 
        \bigr\rvert \bigr\rvert
        \mathbf{L}(\boldsymbol{\mu}^*; \boldsymbol{\theta}_0)\T\mathbf{H}(\boldsymbol{\theta}_0)^{-1}\boldsymbol{\psi}_i(\boldsymbol{\theta}_0)
        \bigr\rvert \bigr\rvert^2_{2}.
    \end{equation} 
    Then the following holds:
    \begin{enumerate}[label = \alph*)]
        \item $\boldsymbol{\mu}^*$ is a stationary point of $\Phi(\boldsymbol{\Gamma}(\boldsymbol{\mu}; \boldsymbol{\theta}_0))$ for a PO-WR or MULTI design of size $n$ if
        \begin{equation}
            \label{eq:mui}
            \mu^*_i = n \frac{\sqrt{c_i}}{\sum_{j \in \mathcal D} \sqrt{c_j}} \quad \text{for all } i \in \mathcal D.
        \end{equation}
        \label{prop:optimality_conditions_PO-WR}
        \item $\boldsymbol{\mu}^*$ is a stationary point of $\Phi(\boldsymbol{\Gamma}(\boldsymbol{\mu}; \boldsymbol{\theta}_0))$ for a PO-WOR design of size $n$ if 
        \begin{subequations}
            \begin{alignat}{2}
                \label{eq:powor_a}
                \mu_i^* & \le 1 && \quad \text{for all } i \in \mathcal D, \\
                \label{eq:powor_b}
                \mu^*_i & = (n - n_{\mathcal E}) \frac{\sqrt{c_i}}{\sum_{j \in \mathcal D \setminus \mathcal E} \sqrt{c_j}} && \quad \text{for all $i \in \mathcal D \setminus \mathcal E$}, \\
                \label{eq:powor_c}
                \sqrt{c_i} & \ge \sqrt{c_j}/\mu_j^* && \quad  \text{for all $i \in \mathcal E$ and $j \in \mathcal D \setminus \mathcal E$},
            \end{alignat}
        \end{subequations}
        where $\mathcal E = \{i \in \mathcal D: \mu_i^* = 1\}$ and $n_{\mathcal E} = |\mathcal E|$.
        \label{prop:optimality_conditions_PO-WOR}
    \end{enumerate}
    Consequently, if $\boldsymbol{\mu}^*$ satisfies the optimality conditions according to \ref{prop:convex_a} or \ref{prop:convex_b}, and $\Phi(\boldsymbol{\Gamma}(\boldsymbol{\mu}; \boldsymbol{\theta}_0))$ is convex in $\boldsymbol{\mu}$, then $\boldsymbol{\mu}^*$ is the global minimiser of $\Phi(\boldsymbol{\Gamma}(\boldsymbol{\mu}; \boldsymbol{\theta}_0))$.
\end{proposition}

We note that the matrix $\mathbf{L}(\boldsymbol{\mu}^*; \boldsymbol{\theta}_0)$ in Proposition \ref{prop:optimality_conditions} exists by Lemma \ref{lemma:chain_rule} whenever the objective function is differentiable at $\boldsymbol{\mu}^*$. It need not be unique, however, and may depend on both $\boldsymbol{\mu}^*$ and $\boldsymbol{\theta}_0$. Some examples can be found in Lemma \ref{lemma:differentiability}. For linear optimality criteria, the matrix $\mathbf{L}(\boldsymbol{\mu}^*; \boldsymbol{\theta}_0)$ does not depend on $\boldsymbol{\mu}^*$ but may depend on the full-data parameter $\boldsymbol{\theta}_0$; see Section \ref{sec:another_look_at_invariance} for further discussion and examples.

The result of Proposition \ref{prop:optimality_conditions} follows from Lemma \ref{lemma:partial_derivatives} by the Lagrange multiplier method in \ref{prop:optimality_conditions_PO-WR} and the Karush-Kuhn-Tucker conditions in \ref{prop:optimality_conditions_PO-WOR}. We show in Proposition \ref{prop:convex} that the D- and L-optimality criteria are convex in $\boldsymbol{\mu}$ and hence that global optimality can be deduced.

\begin{proposition}[Convexity of the D- and L-optimality criteria] 
    \label{prop:convex}
    Consider the family of PO-WR, PO-WOR or multinomial sampling designs of (expected) size $n$. Assume that $\mathbf{H}(\boldsymbol{\theta}_0)$ is of full rank. Then 
    \begin{enumerate}[label = \alph*)]
        \item the L-optimality criterion is convex in $\boldsymbol{\mu}$.
         \label{prop:convex_a}
   \end{enumerate}
    Assume further that $\mathbf{V}(\boldsymbol{\mu}; \boldsymbol{\theta}_0)$, defined in \eqref{eq:vmat}, is positive definite for every $\boldsymbol{\mu} \in \mathcal M_n$. Then 
    \begin{enumerate}[label = \alph*)]
        \setcounter{enumi}{1}
        \item the D-optimality criterion is (log) convex in $\boldsymbol{\mu}$.
        \label{prop:convex_b}
    \end{enumerate}
\end{proposition}

The first assumption in Proposition \ref{prop:convex} is needed to ensure that the inverse of $\mathbf{H}(\boldsymbol{\theta}_0)$ exists, and that the approximate covariance matrix $\boldsymbol{\Gamma}(\boldsymbol{\mu}; \boldsymbol{\theta}_0)$ is well-defined. For the D-optimality criterion we also need that $\boldsymbol{\Gamma}(\boldsymbol{\mu}; \boldsymbol{\theta}_0)$ is of full rank, which follows if the additional assumption on $\mathbf{V}(\boldsymbol{\mu}; \boldsymbol{\theta}_0)$ is fulfilled. We note that this is rather an assumption on the model and data than on the sampling design. Moreover, both of the assumptions in Proposition \ref{prop:convex} hold in most situations. One example where these assumptions are violated, however, is encountered in (multivariate) regression analysis when the model matrix $\mathbf{X}$ or response matrix $\mathbf{Y}$ (i.e., the matrices with rows $\boldsymbol{x}_i\T$ and $\boldsymbol{y}_i\T$) has linearly dependent columns. Another example is logistic regression with complete separation, i.e., when the outcome is linearly separable by the predictors. It is also possible that $\mathbf{H}(\boldsymbol{\theta}_0)$ is of full rank while $\mathbf{V}(\boldsymbol{\mu}; \boldsymbol{\theta}_0)$ is rank-deficient on $\mathcal M_n$. In this case the L-optimality criterion is still well-defined, whereas the D-optimality criterion is not. There are various solutions to such problems, e.g, removing redundant columns from the data, using a ridge penalty to avoid rank-deficiency of the Hessian matrix \citep{Hastie2020}, or by restricting the D-optimality criterion to a subset of the parameters using so called D$_{\mathrm{A}}$-optimality (Table \ref{tab:OED}) \citep{Sibson1974}. Most of these situations may be avoided by a careful construction of the model, however.

Even with a convex objective function, it is possible that no feasible global optimum exist since the domain $\mathcal M_n$ is not closed. For the L-optimality criterion this happens if $c_i$ in \eqref{eq:ci} is equals zero for some $i \in \mathcal D$. In this case the objective function does not depend on the corresponding $\mu_i$ and the partial derivative with respect to $\mu_i$ is equal to zero. The optimal choice would be to correspondingly set $\mu_i = 0$, but this is an unfeasible solution. For any choice of $\mu_i >0$, it is always possible to improve the value of the objective function by reducing $\mu_i$ and distribute the regained probability mass optimally on the remaining elements in $\mathcal D$. The existence of a feasible global optimum can be ensured by imposing the additional restriction that $\mu_i \ge \mu_{\min}$ for all $i \in \mathcal D$, and some $\mu_{\min} \in (0, n/N)$. An alternative solution that does not require explicit specification of a lower bound $\mu_{\min}$, but that still ensures a feasible solution with $\mu_i> 0$, is proposed in Section \ref{sec:practical_implementation}.

\subsection{Optimal sampling schemes}
\label{sec:optimal_sampling_schemes}

In this subsection we present algorithms for finding optimal sampling schemes. First consider a linear optimality criterion with respect to a $p \times m$ matrix $\mathbf{L}$. In this case a closed solution for the optimal sampling scheme is available for the PO-WR and MULTI designs, and given by \eqref{eq:ci}–\eqref{eq:mui} with $\mathbf{L}(\boldsymbol{\mu}^*; \boldsymbol{\theta}_0) = \mathbf{L}$, provided that the corresponding $c_i>0$ for all $i \in \mathcal D$. In particular, A-optimality is obtained with $\mathbf{L}(\boldsymbol{\mu}^*; \boldsymbol{\theta}_0) = \mathbf{I}_{p \times p}$, and c-optimality with $\mathbf{L}(\boldsymbol{\mu}^*; \boldsymbol{\theta}_0) = \mathbf{c}$. For PO-WOR, a simple adjustment may be needed to ensure that a feasible solution with $\mu_i \le 1$ is obtained (Algorithm \ref{alg:lopt}).

\begin{algorithm}[!htb]
    \caption{L-optimal sampling schemes for Poisson and multinomial sampling designs.}
    {\scshape Input}: Index set $\mathcal D$, (expected) sample size $n$, non-zero $p \times m$ matrix $\mathbf{L}$, Hessian matrix $\mathbf{H}(\boldsymbol{\theta}_0)$, gradients $\{\boldsymbol{\psi}_i(\boldsymbol{\theta}_0)\}_{i \in \mathcal D}$, family of sampling designs (PO-WR, PO-WOR or MULTI).
    \begin{algorithmic}[1]
    \label{alg:lopt}
    \STATE Let $c_i = ||\mathbf{L}\T\mathbf{H}(\boldsymbol{\theta}_0)^{-1}\boldsymbol{\psi}_i(\boldsymbol{\theta}_0)||_{2}^2$ for all $i \in \mathcal D$.
    \IF{any $c_i$ = 0}
        \label{step:begin_check_ci}
        \STATE {\scshape Stop}. Feasible solution does not exist. 
    \ELSE 
    \label{step:end_check_ci}
    \STATE Let $\mu_i^* = n \frac{\sqrt{c_i}}{\sum_{j \in \mathcal D} \sqrt{c_j}} \text{ for all } i \in \mathcal D$.
    \label{step:mui}
    \IF{PO-WOR}
    	\WHILE{any $\mu_i^* > 1$} 
        \STATE Let $\mathcal E = \{i\in \mathcal D: \mu_i^* \ge 1\}$ and $n_{\mathcal E} = |\mathcal E|$.
        \label{alg:iter_start}
        \STATE Let 
        \begin{fleqn} 
        $
            \mu_i^* = \begin{cases}
                1 & \text{if } i \in \mathcal E, \\
                (n - n_{\mathcal E}) \frac{\sqrt{c_i}}{\sum_{j \in \mathcal D \setminus \mathcal E} \sqrt{c_j}} & \text{if } i \in \mathcal D \setminus \mathcal E.
            \end{cases}
        $
        \end{fleqn}        \label{alg:iter_end}
        \ENDWHILE
    \ENDIF
    \STATE {\scshape Return} optimal sampling scheme $\boldsymbol{\mu}^* = (\mu_1^*, \ldots, \mu_N^*)$.
    \label{step:output}
    \ENDIF
    \end{algorithmic}
\end{algorithm}

Using the result of Proposition \ref{prop:optimality_conditions} and Algorithm \ref{alg:lopt}, in Algorithm \ref{alg:linearisation} we present an iterative algorithm to find optimal sampling schemes for non-linear optimality criteria. The algorithm takes an initial sampling scheme as input and solves a series of convex optimisation problems by a local approximation of the objective function as linear optimality criterion.  The algorithm is terminated for convergence when the relative improvement of the objective function between two consecutive iterations is less than some pre-specified tolerance level $\epsilon$ (e.g., $\epsilon = 10^{-3}$). The algorithm may also be terminated for divergence if the value of the objective function increases between the iterations. If the algorithm converges, it converges to a fixed-point of the function $\boldsymbol{h}(\boldsymbol{u}): \mathbb{R}^{N} \to \mathbb{R}^{N}$ defined by Algorithm \ref{alg:lopt} with $\mathbf{L} = \mathbf{L}(\boldsymbol{\mu}; \boldsymbol{\theta}_0)$, which by Proposition \ref{prop:optimality_conditions} is a stationary point of $\Phi(\boldsymbol{\Gamma}(\boldsymbol{\mu}; \boldsymbol{\theta}_0))$. For L-optimality, the method is exact and terminates within a single iteration. Beyond L-optimality, the algorithm need not converge, and even if it does, it need not converge to a global optimum unless the problem is convex. The performance of this algorithm for non-linear optimality criteria will be evaluated in Section \ref{sec:application}.

\begin{algorithm}[!htb]
    \caption{Fixed-point iteration.}
    {\scshape Input}: Index set $\mathcal D$, (expected) sample size $n$, optimality criterion $\Phi$, Hessian matrix $\mathbf{H}(\boldsymbol{\theta}_0)$, gradients $\{\boldsymbol{\psi}_i(\boldsymbol{\theta}_0)\}_{i \in \mathcal D}$, initial sampling scheme $\boldsymbol{\mu}_0$, family of sampling designs (PO-WR, PO-WOR or MULTI), maximal number of iterations $T$, tolerance parameter $\epsilon > 0$.
    \begin{algorithmic}[1] 
    \label{alg:linearisation}
	\FOR{t = 1, \ldots, T} 
	    \STATE Let $\mathbf{L}_t$ be a matrix such that $\mathbf{L}_t\mathbf{L}_t\T = \boldsymbol{\phi}(\boldsymbol{\Gamma}(\boldsymbol{\mu}_{t-1}; \boldsymbol{\theta}_0))$.
        \STATE Let $c_i = ||\mathbf{L}_t\T\mathbf{H}(\boldsymbol{\theta}_0)^{-1}\boldsymbol{\psi}_i(\boldsymbol{\theta}_0)||_{2}^2$ for all $i \in \mathcal D$.
        \IF{any $c_i = 0$}
             \STATE {\scshape Stop}. Unfeasible solution encountered during iteration.
        \ELSE 
            \STATE Find L-optimal sampling scheme $\boldsymbol{\mu}_t$ with respect to $\mathbf{L} = \mathbf{L}_t$ according to Algorithm \ref{alg:lopt}. 
            \IF{value of objective function increased} 
                \STATE {\scshape Stop}. Algorithm diverged.
            \ELSIF{relative improvement of the objective function $<\epsilon$} 
                \STATE Algorithm converged. {\scshape Return} $\boldsymbol{\mu}^* = \boldsymbol{\mu}_t$.
            \ENDIF   
        \ENDIF
    \ENDFOR
    \end{algorithmic}
\end{algorithm}

\section{A distance-minimising perspective on optimal subsampling designs}
\label{sec:dist}

Recall the overall aim of data subsampling as introduced in Section \ref{sec:introduction}; to find an approximate solution to the originally intractable problem \eqref{eq:theta0}–\eqref{eq:ell0}. 
A natural target for optimal design in this context is therefore to minimise the expected distance $\E[d(\hat{\boldsymbol{\theta}}_{\!\boldsymbol{\mu}})]$ of the estimator $\hat{\boldsymbol{\theta}}_{\!\boldsymbol{\mu}}$ from the full-data parameter $\boldsymbol{\theta}_0$, for some suitable statistical distance function $d : \boldsymbol{\Omega} \to \mathbb{R}_+$. In Section \ref{sec:d_star_optimality} we define a class of optimality criteria for minimising the expected distance, and discuss their relation to traditional optimality criteria. Some specific examples are presented in Section 
\ref{sec:d_star_examples}, and invariance properties discussed in Section \ref{sec:another_look_at_invariance}.

\subsection{\textit{d}-optimality}
\label{sec:d_star_optimality}

Consider a statistical distance function $d(\boldsymbol{\theta})$ such that $d(\boldsymbol{\theta}) \ge 0$ for all $\boldsymbol{\theta} \in \boldsymbol{\Omega}$, with equality only for $\boldsymbol{\theta} = \boldsymbol{\theta}_0$. For analytical and computational tractability we also require the distance function to be twice differentiable, and let $\mathbf{H}_d(\boldsymbol{\theta}) = \frac{\partial^2 d(\boldsymbol{\theta})}
{\partial\boldsymbol{\theta}\partial \boldsymbol{\theta}^{\mathsf{T}}}$ denote the Hessian matrix of $d(\boldsymbol{\theta})$. We have the following result:

\begin{lemma}[Taylor expansion of $d(\boldsymbol{\theta})$] 
    \label{lemma:d_taylor_expansion}
    Let $\boldsymbol{\theta}_0$ and $\hat{\boldsymbol{\theta}}_{\!\boldsymbol{\mu}}$ be defined according to \eqref{eq:theta0}–\eqref{eq:ell0} and \eqref{eq:thetahat}–\eqref{eq:ellhat}. Assume that \eqref{eq:unbiased}--\eqref{eq:approximate_covariance} hold, and that $\hat{\boldsymbol{\theta}}_{\!\boldsymbol{\mu}}$ has bounded $2+\delta$ moments for some $\delta > 0$. Consider a function $d: \boldsymbol{\Omega} \to \mathbb{R}_+$ such that $d(\boldsymbol{\theta}) = 0$ if and only if $\boldsymbol{\theta} = \boldsymbol{\theta}_0$. Assume that $d(\boldsymbol{\theta})$ is twice differentiable in a neighbourhood of $\boldsymbol{\theta}_0$, and that $\mathbf{H}_d(\boldsymbol{\theta}_0)$ is non-zero. Then 
    \[
    \E[d(\hat{\boldsymbol{\theta}}_{\!\boldsymbol{\mu}})] = 
    \frac{1}{2}\tr\left(
    \boldsymbol{\Gamma}(\boldsymbol{\mu}; \boldsymbol{\theta}_0)
    \mathbf{H}_d(\boldsymbol{\theta}_0)
    \right)
    + o(n^{-1}).
    \]
\end{lemma}

The result of Lemma \ref{lemma:d_taylor_expansion} follows from a Taylor expansion of $d(\boldsymbol{\theta})$ at $\boldsymbol{\theta} = \boldsymbol{\theta}_0$ and properties of quadratic forms. Based on this result, we define a class of expected-distance-minimising optimality criteria as follows:

\begin{definition}[\textit{d}-optimality] 
    \label{def:d_optimality} 
    Consider a function $d:\boldsymbol{\Omega} \to \mathbb{R}_+$ satisfying the conditions of Lemma \ref{lemma:d_taylor_expansion}. Also consider a family of unequal probability sampling designs (e.g., PO-WR, PO-WOR or MULTI) indexed by the sampling scheme $\boldsymbol{\mu}$.  Let the expected size $\E[\sum_{i\in\mathcal D}S_i] = n$ be fixed, and let $\mathcal M_n$ denote the corresponding domain of $\boldsymbol{\mu}$. We say that a sampling scheme $\boldsymbol{\mu}^*$ is \textit{d}-optimal with respect to the statistical distance function $d(\boldsymbol{\theta})$ if 
    \[
    \boldsymbol{\mu}^* = \argmin_{\boldsymbol{\mu} \in \mathcal M_n} \tr
    \left(
    \boldsymbol{\Gamma}(\boldsymbol{\mu}; \boldsymbol{\theta}_0)
    \mathbf{H}_d(\boldsymbol{\theta}_0)
    \right).
    \]
\end{definition}

We denote this optimality criterion as \textit{d}-optimality for \textit{distance}, which should not be confused with the D-optimality criterion introduced in Section \ref{sec:optimality_criteria}. We recognise the \textit{d}-optimality criterion as a linear optimality criterion with $\mathbf{L}\mathbf{L}\T = \mathbf{H}_d(\boldsymbol{\theta}_0)$. Indeed, we have the following equivalence result:

\begin{proposition}[Equivalence between d- and $\Phi$-optimality] \hspace{1cm} 
    \label{prop:equivalence}
    \begin{enumerate}[label = \alph*)]
        \item Consider a function $d: \boldsymbol{\Omega} \to \mathbb{R}_+$ satisfying the conditions of Lemma \ref{lemma:d_taylor_expansion} and denote by $\mathbf{H}_d(\boldsymbol{\theta})$ the Hessian of $d(\boldsymbol{\theta})$. Assume that the sampling scheme $\boldsymbol{\mu}^*$ is \textit{d}-optimal with respect to the distance function $d(\boldsymbol{\theta})$. Then there exists a real matrix $\mathbf{L}$ such that $\mathbf{L}\mathbf{L}\T = \mathbf{H}_d(\boldsymbol{\theta}_0)$ and $\boldsymbol{\mu}^*$ is L-optimal with respect to $\mathbf{L}$.
        \label{prop:equivalence_a}
        \item Let $\Phi:\boldsymbol{S}_{\!+}^{p \times p} \to \mathbb{R}$ and $\mathbf{L}(\boldsymbol{\mu}; \boldsymbol{\theta}_0)$ be defined as in Lemma \ref{lemma:chain_rule} and assume that $\Phi(\boldsymbol{\Gamma}(\boldsymbol{\mu}; \boldsymbol{\theta}_0))$ is differentiable with respect to $\boldsymbol{\mu}$ in a neighbourhood of its optimum argument $\boldsymbol{\mu}^*$. Then $\boldsymbol{\mu}^*$ is \textit{d}-optimal with respect to the distance function $d(\boldsymbol{\theta}) = ||\mathbf{L}(\boldsymbol{\mu}^*; \boldsymbol{\theta}_0)
        \T(\boldsymbol{\theta} - \boldsymbol{\theta_0})||_2^2$.
        \label{prop:equivalence_b}
        \end{enumerate}
\end{proposition}

Proposition \ref{prop:equivalence} follows immediately by the definitions and the optimality conditions of Proposition \ref{prop:optimality_conditions}. By this result, any $\Phi$-optimality criterion may be viewed as minimising the expected distance of the estimator $\hat{\boldsymbol{\theta}}_{\!\boldsymbol{\mu}}$ from the full-data parameter $\boldsymbol{\theta}_0$ for a particular choice of distance function. For instance, A-optimality is equivalent to \textit{d}-optimality with $d(\boldsymbol{\theta}) = ||\boldsymbol{\theta} - \boldsymbol{\theta}_0||_2^2$. Beyond linear optimality criteria, the induced distance function may be implicit and depend on the $\Phi$-optimal sampling scheme $\boldsymbol{\mu}^*$. As an example, E-optimality is equivalent to \textit{d}-optimality with $d(\boldsymbol{\theta}) = ||\boldsymbol{v}_{\boldsymbol{\mu}^*}\T(\boldsymbol{\theta} - \boldsymbol{\theta}_0)||_2^2$, where $\boldsymbol{v}_{\boldsymbol{\mu}^*}$ is an eigenvector pertaining to the largest eigenvalue of $\boldsymbol{\Gamma}(\boldsymbol{\mu}^*; \boldsymbol{\theta}_0)$ and $\boldsymbol{\mu}^*$ the corresponding E-optimal sampling scheme. In this case the distance function for the \textit{d}-optimality criterion can only be evaluated if the E-optimal sampling scheme is known.

\subsection{Some distance-minimising designs}
\label{sec:d_star_examples}

Next we show how \textit{d}-optimality may be used to derive a novel class of linear optimality criteria with good theoretical properties, including transformation- and parameterisation invariance. Consider the following statistical distance functions naturally arising in data subsampling applications and commonly encountered in statistics:

\begin{enumerate}[label = \roman*)]
    \item \textbf{Empirical risk distance}: Since $\boldsymbol{\theta}_0$ is defined as the minimiser of the full-data empirical risk \eqref{eq:ell0}, we may measure of the distance of a parameter value $\boldsymbol{\theta}$ from the full-data parameter $\boldsymbol{\theta}_0$ through the attained value of the empirical risk. We define the empirical risk distance of $\boldsymbol{\theta}$ from $\boldsymbol{\theta}_0$ as $d_{\mathrm{ER}}(\boldsymbol{\theta}) = \ell_0(\boldsymbol{\theta}) - \ell_0(\boldsymbol{\theta}_0)$. 
    \label{list:empirical_risk_distance}

    \item \textbf{Kullback-Leibler divergence}: Consider a random vector $\boldsymbol{Y}$ with probability density function $f_{\boldsymbol{\theta}}(\boldsymbol{y})$ and cumulative distribution function $F_{\boldsymbol{\theta}}(\boldsymbol{y})$. Let $\mathcal Y$ denote the domain of $\boldsymbol{Y}$. The Kullback-Leibler divergence of $f_{\boldsymbol{\theta}}$ from $f_{\boldsymbol{\theta}_0}$ is defined as $\mathrm{KL}\left(f_{\boldsymbol{\theta}_0} || f_{\boldsymbol{\theta}}\right) = \int_{\mathcal Y} \log \frac{f_{\boldsymbol{\theta}_0}\!(\boldsymbol{y})}{f_{\boldsymbol{\theta}}(\boldsymbol{y})} dF_{\boldsymbol{\theta}_0}\!(\boldsymbol{y})$. To allow for covariates, we define the Kullback-Leibler distance of $\boldsymbol{\theta}$ from $\boldsymbol{\theta}_0$ as $d_{\mathrm{KL}}(\boldsymbol{\theta}) = \sum_{i \in \mathcal D} \int_{\mathcal Y}  \log \frac{f_{\boldsymbol{\theta}_0}\!(\boldsymbol{y}|\boldsymbol{x}_i)}{f_{\boldsymbol{\theta}}(\boldsymbol{y}|\boldsymbol{x}_i)}dF_{\boldsymbol{\theta}_0}\!(\boldsymbol{y}|\boldsymbol{x}_i)$.
    \label{list:Kullback_Leibler}    
    
    \item \textbf{Mahalanobis distance}: Consider a probability distribution on $\mathbb{R}^p$ with mean vector $\boldsymbol{\gamma}$ and covariance matrix $\boldsymbol{\Sigma}$. The Mahalanobis distance of a point $\boldsymbol{\theta} \in \mathbb{R}^p$ from the mean $\boldsymbol{\gamma}$ is then given by $\sqrt{(\boldsymbol{\theta} -\boldsymbol{\gamma})\T \boldsymbol{\Sigma}^{-1}(\boldsymbol{\theta} -\boldsymbol{\gamma})}$. We define the squared Mahalanobis distance of $\boldsymbol{\theta}$ from $\boldsymbol{\theta}_0$ with respect to a real, symmetric, positive definite dispersion matrix $\boldsymbol{\Sigma}$ as $d_{\boldsymbol{\Sigma}}(\boldsymbol{\theta}) = (\boldsymbol{\theta} -\boldsymbol{\theta}_0)\T \boldsymbol{\Sigma}^{-1}(\boldsymbol{\theta} -\boldsymbol{\theta}_0)$. 
    \label{list:Mahalanobis_distance}
\end{enumerate}

Four natural choices of the dispersion matrix $\boldsymbol{\Sigma}$ for the Mahalanobis distance are: 
\begin{enumerate}[label = iii.\alph*)]
    \item $\boldsymbol{\Sigma} = \boldsymbol{\Gamma}(\boldsymbol{\mu}; \boldsymbol{\theta}_0)$, the approximate covariance matrix of $\hat{\boldsymbol{\theta}}_{\!\boldsymbol{\mu}}$.
    \label{list:approximate_covariance}
    \item $\boldsymbol{\Sigma} = \mathbf{H}(\boldsymbol{\theta}_0)^{-1}$, which for a parametric model is an estimate of the covariance matrix of $\boldsymbol{\theta}_0$, seen as an estimator of some underlying super-population parameter $\boldsymbol{\theta}^*$. In this case, $\mathbf{H}(\boldsymbol{\theta}_0)$ is also known as the observed Fisher information matrix, often denoted as $\mathbf{I}(\boldsymbol{\theta}_0)$ \citep{Efron1978}. 
    \label{list:observed_fisher_information}
    
    \item $\boldsymbol{\Sigma} = \widetilde{\mathbf{H}}(\boldsymbol{\theta}_0)^{-1}$, where $\widetilde{\mathbf{H}}(\boldsymbol{\theta}_0)$ is defined for a parametric model $f_{\boldsymbol{\theta}}(\boldsymbol{y}|\boldsymbol{x})$ as $\widetilde{\mathbf{H}}(\boldsymbol{\theta}) = \E_{\boldsymbol{y}\sim f_{\boldsymbol{\theta}_0}\!(\boldsymbol{y}|\boldsymbol{x})}[\mathbf{H}(\boldsymbol{\theta}_0)]$. In this case, $\widetilde{\mathbf{H}}(\boldsymbol{\theta}_0)$ is also known as the expected Fisher information matrix, often denoted as $\mathcal{I}(\boldsymbol{\theta}_0)$ \citep{Efron1978}. 
    \label{list:expected_fisher_information}
    
    \item $\boldsymbol{\Sigma} = \mathbf{H}(\boldsymbol{\theta}_0)^{-1}\mathbf{V}(\boldsymbol{\theta}_0)\mathbf{H}(\boldsymbol{\theta}_0)^{-1}$, with 
    \begin{equation}
        \label{eq:vmat_nomu}
        \mathbf{V}(\boldsymbol{\theta}_0) = \sum_{i \in \mathcal D}\boldsymbol{\psi}_i(\boldsymbol{\theta}_0)\boldsymbol{\psi}_i(\boldsymbol{\theta}_0)\T, \quad \boldsymbol{\psi}_i(\boldsymbol{\theta}) = \nabla_{\!\boldsymbol{\theta}}\ell_i(\boldsymbol{\theta}).
    \end{equation}
    This choice of the matrix $\boldsymbol{\Sigma}$ corresponds to the "robust estimator" or "sandwich estimator" of the covariance matrix of $\boldsymbol{\theta}_0$, seen as an estimator of some underlying super-population parameter $\boldsymbol{\theta}^*$ under a semi-parametric or presumably misspecified parametric model \citep{Stefanski2002}.
    \label{list:robust_covariance_matrix_estimator}
\end{enumerate}

We define $d_{\mathrm{ER}}$-, $d_{\mathrm{KL}}$- and $d_{\boldsymbol{\Sigma}}$-optimality accordingly, i.e., as \textit{d}-optimality with the distance function taken as indicated by the subscript. We also define $d_{\mathrm{I}}$-, $d_{\mathcal{I}}$- and $d_{\mathrm{S}}$-optimality as $d_{\boldsymbol{\Sigma}}$-optimality with dispersion matrix $\boldsymbol{\Sigma}$ taken as in \ref{list:observed_fisher_information} (the inverse of the observed information matrix), \ref{list:expected_fisher_information} (the inverse of the expected information matrix) and \ref{list:robust_covariance_matrix_estimator} (the sandwich variance estimator), respectively.

Note that $d_{\mathrm{KL}}$- and $d_{\mathcal{I}}$-optimality are defined for parametric models only, whereas $d_{\mathrm{ER}}$-, $d_{\mathrm{I}}$- and $d_{\mathrm{S}}$-optimality are appropriate also for semi-parametric and distribution-free methods, including estimation of finite population characteristics. For regression problems, the \textit{d}-optimality criterion with the empirical risk distance (i.e., $d_{\mathrm{ER}}$-optimality) is closely related to the V-optimality criterion (Table \ref{tab:OED}, Section \ref{sec:optimality_criteria}). Indeed, these two optimality criteria are equivalent for ordinary least squares regression when $\nu(\boldsymbol{x})$ is the empirical measure on $\{\boldsymbol{x}_i\}_{i \in \mathcal D}$. 

The Mahalanobis distance with $\boldsymbol{\Sigma} = \boldsymbol{\Gamma}(\boldsymbol{\mu}; \boldsymbol{\theta}_0)$ arises by considering the uncertainty of $\hat{\boldsymbol{\theta}}_{\!\boldsymbol{\mu}}$ as an estimator of the full-data parameter $\boldsymbol{\theta}_0$. In contrast, our motivation for the dispersion matrices in \ref{list:observed_fisher_information}–\ref{list:robust_covariance_matrix_estimator} above comes from a super-population viewpoint where $\boldsymbol{\theta}_0$ is seen as an estimator of some underlying parameter $\boldsymbol{\theta}^*$ \citep[cf.][]{Hartley1975}. The different choices of dispersion matrix $\boldsymbol{\Sigma}$ then arise naturally trough different measures of uncertainty associated with the full-data parameter $\boldsymbol{\theta}_0$ \citep[cf.][]{Stefanski2002}. We emphasise, however, that the super-population perspective adopted here is purely rhetorical. The resulting distance functions are equally valid even without any intentions of super-population inference. The significance of these particular choices of distance functions and dispersion matrices are highlighted in Proposition \ref{prop:MD_ER_KL} below and further in Section \ref{sec:another_look_at_invariance}.

\begin{proposition}[$d_{\mathrm{ER}}$- $d_{\mathrm{KL}}$-,  $d_{\boldsymbol{\Sigma}}$-optimality and equivalence with L-optimality] \hspace{1cm} 
    \label{prop:MD_ER_KL}
    \begin{enumerate}[label = \alph*)]
        \item \textit{d}-optimality with respect to the empirical risk distance is equivalent to L-optimality with respect to a $p \times p$ matrix $\mathbf{L}$ such that $\mathbf{L}\mathbf{L}\T = \mathbf{H}(\boldsymbol{\theta}_0)$.
        \label{prop:MD_ER_KL_a}
        
        \item \textit{d}-optimality with respect to the Mahalanobis distance is equivalent to L-optimality with respect to a $p \times p$ matrix $\mathbf{L}$ such that $\mathbf{L}\mathbf{L}\T = \boldsymbol{\Sigma}^{-1}$. 
        \label{prop:MD_ER_KL_b}
        
        \item Consider a parametric statistical model with density function $f_{\boldsymbol{\theta}}(\boldsymbol{y}|\boldsymbol{x})$ and cumulative distribution function $F_{\boldsymbol{\theta}}(\boldsymbol{y}|\boldsymbol{x})$. Let $\boldsymbol{\theta}_0$ be defined by \eqref{eq:theta0}–\eqref{eq:ell0} with $\ell_i(\boldsymbol{\theta}) = -\log f_{\boldsymbol{\theta}}(\boldsymbol{y}_i|\boldsymbol{x}_i)$. Assume that the following holds for all $i \in \mathcal D$ and all parameter values $\boldsymbol{\theta}$ in a neighbourhood or $\boldsymbol{\theta}_0$: $d_{\mathrm{KL}}(\boldsymbol{\theta})$ is finite, $\ell_i(\boldsymbol{\theta})$ is two times continuously differentiable with respect to $\boldsymbol{\theta}$, and all first- and second-order derivatives of $\log f_{\boldsymbol{\theta}}(\boldsymbol{y}|\boldsymbol{x}_i)$ are bounded in $L_1$ with respect to the measure $dF_{\boldsymbol{\theta}}(\boldsymbol{y}|\boldsymbol{x}_i)$.   
        Then \textit{d}-optimality with respect to the Kullback-Leibler distance is equivalent to L-optimality with respect to a $p \times p$ matrix $\mathbf{L}$ such that $\mathbf{L}\mathbf{L}\T = \widetilde{\mathbf{H}}(\boldsymbol{\theta}_0)$.
        \label{prop:MD_ER_KL_c}
    \end{enumerate}
\end{proposition}

The result of Proposition \ref{prop:MD_ER_KL} follows immediately from Proposition \ref{prop:equivalence}. Note that for \ref{prop:MD_ER_KL_c} we need conditions on the model $f_{\boldsymbol{\theta}}(\boldsymbol{y}|\boldsymbol{x})$ that allow us to change the order of integration and differentiation.

By Proposition \ref{prop:MD_ER_KL}\ref{prop:MD_ER_KL_a} and \ref{prop:MD_ER_KL_c} we observe that $d_{\mathrm{ER}}$- and $d_{\mathrm{I}}$-optimality are equivalent (take $\boldsymbol{\Sigma} = \mathbf{H}(\boldsymbol{\theta}_0)^{-1}$). The same also holds for $d_{\mathrm{KL}}$- and $d_{\mathcal{I}}$-optimality (take $\boldsymbol{\Sigma} = \widetilde{\mathbf{H}}(\boldsymbol{\theta}_0)^{-1}$). We also note that for many models, including exponential families and generalised linear models with a canonical link function, the observed information matrix $\mathbf{I}(\boldsymbol{\theta}_0) = \mathbf{H}(\boldsymbol{\theta}_0)$ and expected information matrix $\mathcal{I}(\boldsymbol{\theta}_0) = \widetilde{\mathbf{H}}(\boldsymbol{\theta}_0)$ are equal, and that these four optimality criteria hence are equivalent \citep[see, e.g.][]{McCullagh1989}. For a correctly specified parametric model, they are also asymptotically equivalent to $d_{\mathrm{S}}$-optimality (as $N \rightarrow \infty$), since in this case $N\mathbf{H}(\boldsymbol{\theta}_0)^{-1}$, $N\widetilde{\mathbf{H}}(\boldsymbol{\theta}_0)^{-1}$ and $N\mathbf{H}(\boldsymbol{\theta}_0)^{-1}\mathbf{V}(\boldsymbol{\theta}_0)\mathbf{H}(\boldsymbol{\theta}_0)^{-1}$ all converge to the same limit \citep[see, e.g.][]{Stefanski2002}.

The above-mentioned optimality criteria are also related to A-optimality after an appropriate change of variables. Consider, e.g., a linear regression model, and assume that the model matrix $\mathbf{X}$ (i.e., the matrix with rows $\boldsymbol{x}_i\T$) has orthogonal columns. Then the $d_{\mathrm{ER}}$- and $d_{\mathrm{KL}}$-optimality criteria are equivalent to A-optimality, since in this case $\mathbf{H}(\boldsymbol{\theta}_0) = \widetilde{\mathbf{H}}(\boldsymbol{\theta}_0) \propto  \mathbf{X}\T\mathbf{X} =  \mathbf{I}_{p \times p}$. In the non-orthogonal case, the $d_{\mathrm{ER}}$- and $d_{\mathrm{KL}}$-optimality criteria depend on the parameterisation of the model and on the scaling of the data and correlations between the variables, through the Hessian $\mathbf{H}(\boldsymbol{\theta}_0)$. As a consequence, invariance under non-singular affine transformations of the data and under a re-parameterisation of the model is achieved (see Section \ref{sec:another_look_at_invariance}). Geometrically, the A-optimality criterion minimises the expected Euclidean distance of the estimator $\hat{\boldsymbol{\theta}}_{\!\boldsymbol{\mu}}$ from the full-data parameter $\boldsymbol{\theta}_0$ (Proposition \ref{prop:equivalence}, Section \ref{sec:d_star_optimality}). The $d_{\mathrm{ER}}$- and $d_{\mathrm{KL}}$-optimality criteria minimise the expected distance with respect to the natural geometry of the model space.

Finally we consider the relation between \textit{d}-optimality and D-optimality. These two criteria coincide if the distance function is taken as the squared Mahalanobis distance $d_{\boldsymbol{\Sigma}}(\boldsymbol{\theta})$ with dispersion matrix $\boldsymbol{\Sigma} = \boldsymbol{\Gamma}(\boldsymbol{\mu^*}; \boldsymbol{\theta}_0)$, where $\boldsymbol{\mu^*}$ is the D-optimal sampling scheme (see Proposition \ref{prop:equivalence} and Proposition \ref{prop:MD_ER_KL}\ref{prop:MD_ER_KL_b}). In particular, D-optimality is equivalent to L-optimality with $\mathbf{L} = \mathbf{H}(\boldsymbol{\theta}_0)\mathbf{V}(\boldsymbol{\mu}^*; \boldsymbol{\theta}_0)^{-1/2}$, and with $\mathbf{V}(\boldsymbol{\mu}; \boldsymbol{\theta}_0)$ defined as in \eqref{eq:vmat}. This result is not very practical, however, since the coefficient matrix of the L-optimality criterion depends on the D-optimal sampling scheme $\boldsymbol{\mu}^*$. An optimality criterion closely related to D-optimality is L-optimality with $\mathbf{L} = \mathbf{H}(\boldsymbol{\theta}_0)\mathbf{V}(\boldsymbol{\theta}_0)^{-1/2}$, where $\mathbf{V}(\boldsymbol{\theta}_0)$ given by \eqref{eq:vmat_nomu} does not depend on $\boldsymbol{\mu}$. By Proposition \ref{prop:MD_ER_KL}\ref{prop:MD_ER_KL_b}, this is equivalent to $d_{\mathrm{S}}$-optimality.

We point out that having the coefficient matrix $\mathbf{L}$ depending on the full-data Hessian $\mathbf{H}(\boldsymbol{\theta}_0)$ and parameter $\boldsymbol{\theta}_0$ is not restrictive, since all optimal designs anyway depend on unknown full-data characteristics. Methods to handle this issue will be addressed in Section \ref{sec:practical_implementation}.

\subsection{Invariance properties}
\label{sec:another_look_at_invariance}

In addition to their appealing geometric and statistical interpretation, the expected-distance-minimising optimality criteria introduced in the previous section have two desirable properties: computational tractability and parameterisation invariance. Indeed, belonging to the class of linear optimality criteria, the $d_{\mathrm{ER}}$-, $d_{\mathrm{KL}}$ and $d_{\mathrm{S}}$-optimality criteria have simple solutions for the optimal sampling schemes according to Algorithm \ref{alg:lopt}. The invariance properties of these optimality criteria and their corresponding optimal sampling schemes are established below.

Consider a re-parameterisation $\boldsymbol{g}: \boldsymbol{\theta} \mapsto \boldsymbol{\eta}$, where $\boldsymbol{g}$ is a one-to-one differentiable mapping on the parameter space. Under such a transformation the full-data empirical risk minimiser $\boldsymbol{\theta}_0$ is equivariant in the sense that the minimiser of the induced empirical risk 
$
\ell^*_{0}(\boldsymbol{\eta}) := \sum_{i\in \mathcal D} \ell_i(\boldsymbol{g}^{-1}(\boldsymbol{\eta}))
$
is given by $\boldsymbol{\eta}_0 = \boldsymbol{g}(\boldsymbol{\theta}_0)$ \citep[see, e.g.,][]{Casella2002}. By similar arguments, the Hansen-Hurwitz empirical risk minimiser for $\boldsymbol{\eta}_0$ is given by $\hat{\boldsymbol{\eta}}_{\!\boldsymbol{\mu}} = \boldsymbol{g}(\hat{\boldsymbol{\theta}}_{\!\boldsymbol{\mu}})$. Evaluating the derivatives of the induced empirical risk $\ell^*_{0}(\boldsymbol{\eta})$, by \eqref{eq:approximate_covariance} we obtain the covariance matrix of $\hat{\boldsymbol{\eta}}_{\!\boldsymbol{\mu}}$ as 
\begin{align}
    \label{eq:cov_eta}
    \COV(\hat{\boldsymbol{\eta}}_{\!\boldsymbol{\mu}} - \boldsymbol{\eta}_0) = \boldsymbol{\Gamma}_{\!\boldsymbol{g}}(\boldsymbol{\mu}; \boldsymbol{\theta}_0) + o(n^{-1}), \quad
    \boldsymbol{\Gamma}_{\!\boldsymbol{g}}(\boldsymbol{\mu}; \boldsymbol{\theta}_0) = 
    \mathbf{J}_{\boldsymbol{g}}(\boldsymbol{\theta}_0) \boldsymbol{\Gamma}(\boldsymbol{\mu}; \boldsymbol{\theta}_0) \mathbf{J}_{\boldsymbol{g}}(\boldsymbol{\theta}_0)\T,    
\end{align}
where $\mathbf{J}_{\boldsymbol{g}}(\boldsymbol{\theta})$ is the Jacobian of $\boldsymbol{g}$, i.e,. the matrix with rows $\nabla_{\boldsymbol{\theta}} g_i(\boldsymbol{\theta})\T$. We say that an optimality criterion is invariant under a re-parameterisation $\boldsymbol{g}: \boldsymbol{\theta} \mapsto \boldsymbol{\eta}$ if the optimal sampling schemes for $\boldsymbol{\Gamma}(\boldsymbol{\mu}; \boldsymbol{\theta}_0)$ and $\boldsymbol{\Gamma}_{\!\boldsymbol{g}}(\boldsymbol{\mu}; \boldsymbol{\theta}_0)$ are equal. Invariance of the  $d_{\mathrm{ER}}$-, $d_{\mathrm{KL}}$- and $d_{\mathrm{S}}$-optimality criteria is established in Proposition \ref{prop:invariance}.

\begin{proposition}[Parameterisation invariance] 
    \label{prop:invariance}
    Let $\mathbf{V}(\boldsymbol{\theta}_0)$ be defined as in \eqref{eq:vmat_nomu}, and assume that $\mathbf{H}(\boldsymbol{\theta}_0)$ and $\mathbf{V}(\boldsymbol{\theta}_0)$ are of full rank. Then the $d_{\mathrm{ER}}$- and $d_{\mathrm{S}}$-optimality criteria are invariant under a re-parameterisation $\boldsymbol{g}: \boldsymbol{\theta} \mapsto \boldsymbol{\eta}$, where $\boldsymbol{g}$ is a one-to-one differentiable mapping on the parameter space. Under the assumptions of Proposition \ref{prop:MD_ER_KL}\ref{prop:MD_ER_KL_c}, the same also holds for the $d_{\mathrm{KL}}$-optimality criterion.
\end{proposition}

Similar results may also be obtained for invariance under non-singular affine transformations of the data. Indeed, in many cases a transformation of the data induces a transformation on the parameter space that satisfies the conditions on the transformation $\boldsymbol{g}$ in Proposition \ref{prop:invariance}. Care needs to be taken, however, to make sure that the empirical risk function is still defined after applying the transformation, and that the transformation produces a mathematically equivalent model. Under such circumstances, the notions of transformation- and parameterisation-invariance are interchangeable in most practical situations. Exceptions exist, however, where a transformation of the data renders the Hessian $\mathbf{H}(\boldsymbol{\theta}_0)$ unchanged. In such a case, the $d_{\mathrm{ER}}$- and $d_{\mathrm{KL}}$-optimality criteria are no longer invariant under affine transformations of the data. We provide such an example in Section \ref{sec:finite_population_inference}. We note that even in such cases the D- and $d_{\mathrm{S}}$-optimality criteria remain invariant under affine transformations of the data.

\section{Practical implementation}
\label{sec:practical_implementation}

Thus far, we have assumed the full data $\{(\boldsymbol{x}_i, \boldsymbol{y}_i)\}_{i \in \mathcal D}$ and full-data parameter $\boldsymbol{\theta}_0$ to be known. However, if such information were available at the design stage, subsampling would not be needed in the first place. In this section we describe a practical approach to optimal subsampling. In Section \ref{sec:anticipated_variance} we introduce the anticipated covariance matrix \citep[cf.][]{Isaki1982} to be used in the optimisation as a surrogate for the unknown covariance matrix $\boldsymbol{\Gamma}(\boldsymbol{\mu}; \boldsymbol{\theta}_0)$. Sequential optimal design and multi-stage sampling procedures, where the information needed for the optimisation is acquired gradually during the sampling process, are discussed in Section \ref{sec:multistage_sampling}.

\subsection{Auxiliary-variable-assisted subsampling designs}
\label{sec:anticipated_variance}

In addition to the data $\{(\boldsymbol{x}_i, \boldsymbol{y}_i)\}_{i \in \mathcal D}$, we now assume the existence of a collection of auxiliary variables $\{\boldsymbol{z}_i\}_{i \in \mathcal D}$, which are available \textit{a priori} for all members $i \in \mathcal D$. Depending on context, the auxiliary variables may include some of the variables in $\boldsymbol{x}_i$ and/or some of the variables in $\boldsymbol{y}_i$. For instance, consider a case-control study to investigate the effect of some exposure variables on a known binary outcome. In this case the auxiliary variables contain the (scalar) outcome $y_i$, and possibly some of the explanatory variables or some proxies for those \citep[cf.][]{Imberg2022a}. The opposite situation is encountered in active learning \citep{Settles2012}. In this case all predictor vectors $\boldsymbol{x}_i$ are known but the outcomes $y_i$ can be observed only for a subset $\mathcal S \subset \mathcal D$, hence $\boldsymbol{z}_i = \boldsymbol{x}_i$ \citep[cf.][]{Bach2007, Wang2017, Meng2021, Zhang2021, Imberg2022b}. In the extreme case, one may even have access to the full-data $\{(\boldsymbol{x}_i, \boldsymbol{y}_i)\}_{i \in \mathcal D}$, but using this information to calculate $\boldsymbol{\theta}_0$ may be too computationally demanding to be feasible \citep[see, e.g.][]{Ma2015, Drovandi2017, Wang2018, Deldossi2022}. Any case in between those extremes may be encountered in practice. The auxiliary variables may be weakly, strongly, or even perfectly correlated with the unobserved study variables. The stronger the correlation, the greater the potential benefits of optimal sampling.

The algorithms presented in Section \ref{sec:optimal_sampling_schemes} for finding optimal sampling schemes require information about the full-data Hessian matrix $\mathbf{H}(\boldsymbol{\theta})$ and gradients $\boldsymbol{\psi}_i(\boldsymbol{\theta})$, evaluated at the full-data parameter $\boldsymbol{\theta}_0$. Moreover, the Hessian depends on the explanatory variables $\boldsymbol{x}_i$, if such are included in the model, and sometimes also on the outcomes $\boldsymbol{y}_i$. Similarly, the gradients depend on both the outcomes and the explanatory variables. To handle this we introduce a collection of random variables $\{(\mathbf{X}_i, \mathbf{Y}_i)\}_{i \in \mathcal D}$ to describe our uncertainty in the unknown values of the data $\{\boldsymbol{x}_i,\boldsymbol{y}_i\}_{i \in \mathcal D}$. For any variable also included in $\boldsymbol{z}_i$, we may associate a degenerate (deterministic) distribution with the corresponding component of $(\mathbf{X}_i, \mathbf{Y}_i)$ conditioned on $\boldsymbol{z}_i$. We also assume that we have a preliminary estimate $\tilde{\boldsymbol{\theta}}_0$ of the full-data parameter $\boldsymbol{\theta}_0$, and an auxiliary model $f(\boldsymbol{x}, \boldsymbol{y}|\boldsymbol{z})$ for the conditional distribution of the random variables $(\mathbf{X}_i, \mathbf{Y}_i)$ given auxiliary variables $\boldsymbol{z}_i$. Such information may be available from domain knowledge, previous studies, a pilot sample, or a combination of those. In Section \ref{sec:multistage_sampling} we will discuss how such information can be acquired gradually during the subsampling process. Below we define the anticipated covariance matrix as the target of optimisation under an assisting auxiliary model for the unknowns.

\begin{definition}[Anticipated covariance] 
    \label{def:anticipated_covariance}
    Consider a data triplet $\{(\boldsymbol{X}_i, \boldsymbol{Y}_i, \boldsymbol{z}_i)\}_{i \in \mathcal D}$, where $(\mathbf{X}_i, \mathbf{Y}_i)$ is a random vector and $\boldsymbol{z}_i$ are known for all $i \in \mathcal D$. Also consider a preliminary estimate $\tilde{\boldsymbol{\theta}}_0$ of the full-data parameter $\boldsymbol{\theta}_0$, and a model $f(\boldsymbol{x}, \boldsymbol{y}|\boldsymbol{z})$ for the conditional distribution of $(\mathbf{X}_i, \mathbf{Y}_i)$ given auxiliary variables $\boldsymbol{z}_i$. The anticipated covariance matrix of $\hat{\boldsymbol{\theta}}_{\!\boldsymbol{\mu}}$ is defined as
    \[
    \widetilde{\boldsymbol{\Gamma}}(\boldsymbol{\mu}; \tilde{\boldsymbol{\theta}}_0) = \E_{(\boldsymbol{x},\boldsymbol{y})\sim f(\boldsymbol{x}, \boldsymbol{y}|\boldsymbol{z})}\!\left[\boldsymbol{\Gamma}(\boldsymbol{\mu}; \boldsymbol{\theta}_0)\right]_{\boldsymbol{\theta}_0 = \tilde{\boldsymbol{\theta}}_0}.
    \]
\end{definition}

The anticipated covariance matrix in Definition \ref{def:anticipated_covariance} is our prediction of the actual unknown covariance matrix $\boldsymbol{\Gamma}(\boldsymbol{\mu}; \boldsymbol{\theta}_0)$, given the available auxiliary information. We use the term \textit{anticipated} rather than \textit{expected}, as adopted from \citet{Isaki1982}, to emphasise that the expectation involved in the above definition is a hypothetical construct and generally differs from the expectation under the data generating mechanism.

All results in Section \ref{sec:optimal_subsampling_designs} and \ref{sec:dist} may now be restated for $\Phi$-optimality with respect to the anticipated covariance matrix $\widetilde{\boldsymbol{\Gamma}}(\boldsymbol{\mu}; \tilde{\boldsymbol{\theta}}_0)$ instead of the approximate covariance matrix $\boldsymbol{\Gamma}(\boldsymbol{\mu}; \boldsymbol{\theta}_0)$. Under weak assumptions on the model $f(\boldsymbol{x}, \boldsymbol{y}|\boldsymbol{z})$ that allow us to replace the order of integration and differentiation, all that changes is that the coefficients $c_i$ in Algorithm \ref{alg:linearisation} are replaced by their corresponding expectations
\begin{equation}
    \label{eq:Ci_random}
    \tilde{c}_i := 
    \E_{(\boldsymbol{x},\boldsymbol{y})\sim f(\boldsymbol{x}, \boldsymbol{y}|\boldsymbol{z})}\left[C_i|\{\boldsymbol{z}_i\}_{i \in \mathcal D}\right], \quad  
    C_i = \bigr\rvert \bigr\rvert
    \mathbf{L}_t\T\mathbf{H}(\boldsymbol{\theta}_0)^{-1}\boldsymbol{\psi}_i(\boldsymbol{\theta}_0)
    \bigr\rvert \bigr\rvert^2_{2, \boldsymbol{\theta}_0 = \tilde{\boldsymbol{\theta}}_0},   
\end{equation}
where $C_i$ is a function of the random variables $\{(\mathbf{X}_i, \mathbf{Y}_i)\}_{i \in \mathcal D}$, and $\mathbf{L}_t$ a matrix such that $\mathbf{L}_t\mathbf{L}_t\T = \boldsymbol{\phi}(\widetilde{\boldsymbol{\Gamma}}(\boldsymbol{\mu}_{t-1}; \tilde{\boldsymbol{\theta}}_0))$.

We note that $C_i$ in \eqref{eq:Ci_random} is a positive random variable, which implies that $\tilde{c}_i > 0$ as long as $C_i > 0$ with positive probability. This is fulfilled whenever the covariance matrices for the components of $(\mathbf{X}_i, \mathbf{Y}_i)$ not included in $\boldsymbol{z}_i$ are of full rank for all $i$. Hence, considering the anticipated covariance under an auxiliary distribution that properly acknowledge the uncertainty in the unknowns, we effectively avoid the situation where the presented algorithms (Algorithm \ref{alg:lopt} and \ref{alg:linearisation}) converge to an unfeasible solution.

\subsection{Sequential optimal design}
\label{sec:multistage_sampling}

The anticipated covariance introduced in the previous section takes us one step closer to a practical framework for
optimal subsampling. With this notion, optimal sampling schemes may be found using the methods of Section \ref{sec:optimal_sampling_schemes}, with the unknown values of the coefficients $c_i$ replaced by their expectations \eqref{eq:Ci_random} under an assisting auxiliary model $f(\boldsymbol{x}, \boldsymbol{y}|\boldsymbol{z})$ and a preliminary parameter estimate $\tilde{\boldsymbol{\theta}}_0$. In most cases, however, even this information is unavailable before any data is observed. This problem may be approached using sequential optimal design. Hence, subsampling is performed in multiple stages, where the information acquired from previous sampling stages may be utilised to devise optimal sampling schemes in succeeding stages. We acknowledge that many algorithms and methods in this spirit have already been presented \citep[see, e.g.][]{Bach2007, Wang2018, Imberg2020, Ai2021_regression}. A general procedure is presented in Algorithm \ref{alg:sequential_optimal_designs}.

\begin{algorithm}
    \caption{K-stage subsampling procedure.}
    \label{alg:sequential_optimal_designs}
    \textbf{Input}: Index set $\mathcal D$, optimality criterion $\Phi$, family of sampling designs (PO-WR, PO-WOR or MULTI), number of sampling stages $K$, batch sizes $\{n_k\}_{k=1}^K$. 
    \begin{algorithmic}[1]
        \FOR{$k = 1, 2, \ldots, K$}
            \STATE Calculate (optimal) sampling scheme.
            \STATE Select a random subsample of size $n_k$.
            \STATE Estimate the target parameter $\boldsymbol{\theta}_0$.
            \STATE Update the auxiliary model $f(\boldsymbol{x}, \boldsymbol{y}|\boldsymbol{z})$.
            \label{step:optimal_sampling_scheme}
            \STATE Evaluate performance/precision.
            \STATE {\scshape Stop} if sufficient precision in reached. {\scshape Else} continue.
        \ENDFOR
    \end{algorithmic}   
    \footnotesize
\end{algorithm}

The number of sampling stages $K$ in Algorithm \ref{alg:sequential_optimal_designs} may range from a single stage with $n$ observations, to $n$ stages with a single observation in each subsample. In linear regression, for instance, there is no need for sequential subsampling if the explanatory variables $\boldsymbol{x}_i$ are known. This holds since in this case \eqref{eq:Ci_random} is a function of the predictors $\boldsymbol{x}_i$ (which are known), the Hessian $\mathbf{H}(\boldsymbol{\theta}_0)$ (which only depends on the predictors $\boldsymbol{x}_i$), and the second moments of the residuals. See \citet{Ma2020} for various optimality criteria and corresponding optimal sampling schemes in this context. At the other extreme, active learning methods utilise a large number of sampling stages, often with a single observation per stage to gain maximal flexibility in the sampling process \citep{Bach2007, Imberg2020, Kossen2022, Zhan2022}. Subsampling methods in big data often rely on two sampling stages: an initial simple random sample followed by an optimal unequal probability sample \citep{Wang2018, Ai2021_regression, Wang2021}.

An estimator for $\boldsymbol{\theta}_0$ after $k$ sampling stages may be defined as
\begin{gather}
    \label{eq:thetahat_seq}    
    \hat{\boldsymbol{\theta}}_{\boldsymbol{\mu}}^{(k)} 
    = \argmin_{\boldsymbol{\theta} \in \boldsymbol{\Omega}} \hat{\ell}_{\boldsymbol{\mu}}^{(k)}(\boldsymbol{\theta}), \\ 
    \nonumber
    \hat{\ell}_{\!\boldsymbol{\mu}}^{(k)}(\boldsymbol{\theta}) 
    = m_k^{-1}\sum_{j=1}^k n_j\hat{\ell}_{\boldsymbol{\mu},j}(\boldsymbol{\theta}), \quad
    \hat{\ell}_{\boldsymbol{\mu},j}(\boldsymbol{\theta})
    = \sum_{i \in \mathcal D} S_{ji}w_{ji}\ell_i(\boldsymbol{\theta}),
\end{gather}
where $S_{ji}$ is the number of times an instance $i \in \mathcal D$ is selected by the sampling mechanism at stage $j$, $\mu_{ji}$ the corresponding expected number of selections, $m_k = n_1 + \ldots + n_k$ the cumulative sample size after $k$ stages, and $w_{ji} = 1/\mu_{ji}$. Here $ \hat{\ell}_{\boldsymbol{\mu},j}(\boldsymbol{\theta})$ is an unbiased Hansen-Hurwitz estimator of the full-data empirical risk $\ell_0(\boldsymbol{\theta})$ from the sample obtained at stage $j$, and $\hat{\ell}_{\!\boldsymbol{\mu}}^{(k)}(\boldsymbol{\theta})$ a pooled estimator calculated from the first $k$ subsamples.

The properties of the resulting estimator \eqref{eq:thetahat_seq}, have been studied in some specific cases, where it has been proven that under suitable regularity conditions the estimator $\hat{\boldsymbol{\theta}}_{\boldsymbol{\mu}}^{(k)}$ is asymptotically normally distributed and consistent for $\boldsymbol{\theta}_0$. See, e.g., \citet{Ai2021_regression} and \citet{Yu2022} for results on generalised linear models and quasi-likelihood methods when the number of sampling stages $K=2$. \citet{Imberg2022b} established the asymptotic properties of estimators for finite population vector characteristics when the subsample sizes $n_k$ are bounded and the number of sampling stages $K\rightarrow \infty$. Combining martingale limit theory \citep{Hall1980} with the asymptotics of estimating equation estimators in survey sampling \citep{Binder1983}, consistency and asymptotic normality of \eqref{eq:thetahat_seq} when the batch sizes $n_k$ are bounded and $K \rightarrow \infty$ may also be deduced \citep[cf.][]{Zhang2021}. We conjecture that a similar result holds also in the case when the number of sample stages $K$ is bounded and the subsample sizes $n_k$ tend to infinity, along with $N \rightarrow \infty$ and $n_k /N \rightarrow \gamma_k \in (0,1)$. A thorough treatment of this issue, however, is a topic for future research.

\section{Application and Examples}
\label{sec:application}

There is already an extensive amount of publications demonstrating the benefits of optimal subsampling; see, e.g., the references in Section \ref{sec:introduction}. We will not provide further evidence for these already convincing results. Instead, in this section we illustrate the presented methodology through examples, and compare different optimality criteria for data subsampling in terms of computation aspects and estimator efficiency.

We consider an application in scenario generation for virtual safety assessment of an advanced driver assistance system. A brief background to the application, description of the data and problem formulation is provided in Section \ref{sec:background}. Examples, illustrations and results for parametric density estimation are presented in Section \ref{sec:parametric_density_estimation}, regression modelling in Section \ref{sec:regression_modelling}, and finite population inference in Section \ref{sec:finite_population_inference}.

\subsection{Materials and methods}
\label{sec:background}

\paragraph{Background} Road traffic injuries is a major cause of death worldwide \citep{WHO2018}. Countermeasures, such as advanced driver assistance systems, are constantly developed to mitigate these risks. One way to evaluate such systems before they enter the market is through virtual simulations \citep{Anderson2013, Seyedi2021}. Since such evaluations are performed in a virtual rather than physical test environment, they are more cost-efficient than traditional test beds. This, however, comes at the cost of a huge computational load. Computation demands can be substantially reduced through subsampling \citep{Mullins2018, Imberg2022b, Sun2022}.

\paragraph{Dataset} Our dataset consists of 44,220 observations generated through variations of 44 reconstructed real rear-end crashes. The variations were generated by altering the driver behaviour of the ensuing vehicle in terms of glance behaviour (off-road glance duration after a specific anchoring point in time) and braking profile (maximal deceleration during braking). For each such variation, a corresponding scenario was setup in a virtual environment and simulation software, through which the entire course of events could be simulated. The outcomes of such a simulation include whether a collision occurred or not, and the impact speed if there was a collision. Thus, each observation in the dataset represents a synthetic event that describes what could have happened in the original crash event under certain variations of the conditions. Each scenario was further run under two 'treatment conditions': a scenario with an advanced emergency braking (AEB) system, and a baseline manual driving scenario without the AEB.

The following variables are included in the dataset:
\begin{itemize}
    \item Input variables: case identifier (categorical with 44 levels corresponding to the 44 original rear-end crashes) off-road glance duration (67 levels, 0–6.6 s), and maximal deceleration during braking (15 levels, 3.3–10.3 m/s$^2$).
    \item Direct outcomes: crash indicator (1 if there was a collision and 0 otherwise) and impact speed with the AEB system and under the baseline manual driving scenario.
    \item Calculated outcomes: injury risk with the AEB system and under the baseline manual driving scenario, impact speed reduction, injury risk reduction, and crash avoidance indicator with the AEB system compared to baseline manual driving. 
\end{itemize}
Associated with each observation is also an observation weight $w_i > 0$, describing the probability of the specific input parameter configuration (i.e., off-road glance duration and maximal deceleration during braking) occurring in real life. Additional details may be found in \citet{Imberg2022b}.

\paragraph{Target characteristics} We are interested in the following:
\begin{enumerate}[label = \roman*)]
    \item The impact speed distribution under the baseline scenario, restricted to the subset of input values that produce a crash.
    \label{list:impact_speed_distribution}
    \item The impact speed response surface under the baseline scenario, as a function of the off-road glance duration and maximal deceleration.
    \label{list:response_surface_modelling}
    \item The mean impact speed reduction, mean injury risk reduction, and crash avoidance rate with the AEB compared to baseline manual driving, restricted to the subset of variations for which there is a crash in the baseline scenario.
    \label{list:finite_populaiton_inference}
\end{enumerate}

Characteristics of the dataset, including the baseline impact speed distribution, impact speed response surface, and safety benefit distribution of the AEB compared to baseline manual driving, are presented graphically in Figure \ref{fig:descriptives} and \ref{fig:response_surface} in Appendix \ref{sec:supplementary_figures}.

As often is the case in practice, we assume that running all simulations of interest is practically unfeasible and subsampling inevitable. In such a case, the input variables (i.e., case identifier, off-road glance duration, and maximal deceleration during braking) and scenario probabilities are available \textit{a priori} for all instances in the dataset. Hence, these are our auxiliary variables. The remaining variables can only be observed for a subset on which inference will be based. For simplicity, we restrict our consideration in problem \ref{list:impact_speed_distribution} (Section \ref{sec:parametric_density_estimation}) and \ref{list:finite_populaiton_inference} (Section  \ref{sec:finite_population_inference}) to simulations that produce a crash in the baseline scenario. Thus, the 4299 observations that did not result in a crash are excluded from the corresponding evaluations.

\paragraph{Performance evaluation}

We evaluate the performance of the proposed optimal subsampling methods in terms of computation time and statistical efficiency on the application and inference problems described above. Also, for non-linear optimality criteria, we evaluate the number of iterations needed for convergence of the fixed-point iteration algorithm (Algorithm \ref{alg:linearisation}, Section \ref{sec:optimal_sampling_schemes}), i.e., the time it takes to find the optimal sampling scheme. For a sampling scheme $\boldsymbol{\mu}$, the statistical efficiency of the estimator $\hat{\boldsymbol{\theta}}_{\!\boldsymbol{\mu}}$ with respect to a criterion $\Phi$ is measured by the relative $\Phi$-efficiency
\begin{equation*}
    \Phi\text{-eff}(\boldsymbol{\mu}) = \Phi(\boldsymbol{\Gamma}(\boldsymbol{\mu}^*; \boldsymbol{\theta}_0)) / \Phi(\boldsymbol{\Gamma}(\boldsymbol{\mu}; \boldsymbol{\theta}_0)),
\end{equation*}
where $\boldsymbol{\mu}^*$ is the $\Phi$-optimal sampling scheme \citep{Atkinson1992, Pukelsheim1993}. The relative $\Phi$-efficiency measures the extent to which the sampling scheme $\boldsymbol{\mu}$ exhausts the maximum information for $\boldsymbol{\theta}_0$ with respect to the criterion $\Phi$. Its inverse is the relative increase in the sample size needed to reach the same level of performance as the optimal design with respect to the $\Phi$-optimality criterion. The relative efficiencies are evaluated analytically using the expression \eqref{eq:approximate_covariance} for the approximate covariance matrix.

The following optimality criteria are considered: A-, c-, D-, and E-optimality, $\Phi_q$-optimality with $q=0.5, q=5$, and $q=10$, $d_{\mathrm{ER}}$-optimality (i.e., L-optimality with $\mathbf{L} = \mathbf{H}(\boldsymbol{\theta}_0)^{1/2}$, which for all models in this evaluation also is equivalent to $d_{\mathrm{KL}}$-optimality), and $d_{\mathrm{S}}$-optimality (i.e., L-optimality with $\mathbf{L}= \mathbf{H}(\boldsymbol{\theta}_0)\mathbf{V}(\boldsymbol{\theta}_0)^{-1/2}$). See Section \ref{sec:optimality_criteria} and \ref{sec:d_star_examples} for additional details and definitions. For the D-optimality criterion, the non-logarithmic version of the objective function ($\det(\boldsymbol{\Gamma})^{1/p}$) is used (Table \ref{tab:OED}, Section \ref{sec:optimality_criteria}).

All algorithms and evaluations are implemented using the \texttt{R} language and environment for statistical computing, version 4.2.3 \citep{RCoreTeam}. Computations are carried out using a single core on a desktop running Windows 11 with an 2.1 GHz Intel i7 processor. The subsample size is set to 1\% of the full-data size. Sampling schemes for linear optimality criteria are calculated according Algorithm \ref{alg:lopt}, and sampling schemes for non-linear optimality criteria are calculated according to Algorithm \ref{alg:linearisation} with tolerance parameter $\epsilon = 0.001$. The full data $\{(\boldsymbol{x}_i, \boldsymbol{y}_i)\}_{i\in \mathcal D}$ and full-data parameter $\boldsymbol{\theta}_0$ are assumed to be known, so that the theoretically optimal sampling schemes can be found. The dataset and \texttt{R} code is available online at \href{https://github.com/imbhe/OSD}{\nolinkurl{https://github.com/imbhe/OSD}}.

Results are presented for PO-WR and multinomial sampling designs, which produce identical analytical results. By similar means, analogous results may be obtained for PO-WOR.

\subsection{Parametric density estimation}
\label{sec:parametric_density_estimation}

First we consider the distribution of the impact speed under the baseline scenario, illustrated in Figure \ref{fig:descriptives} in Appendix \ref{sec:supplementary_figures}.

\paragraph{Model} 

The impact speed is assumed to follow a log-normal distribution with parameter $(\eta, \sigma)$ for the mean and standard deviation of the log impact speed. The full-data parameter $\boldsymbol{\theta}_0 = (\eta_0, \sigma_0)\T$ is defined as
\begin{equation}
    \label{eq:lognormal}
    \boldsymbol{\theta}_0 = \argmin_{\eta \in \mathbb{R}, \sigma \in \mathbb{R}_{>0}} \frac{1}{2}\sum_{i=1}^N w_i \left(\frac{(\log y_i - \eta)^2}{\sigma^2} + \log \sigma^2 \right),
\end{equation}
where $w_i$ is an observation weight known \textit{a priori}, and $y_i$ is the impact speed in scenario $i \in \mathcal D$. Without loss of generality, we assume that the observation weights have been normalised so that $\sum_{i=1}^N w_i = 1$.

\paragraph{Optimal sampling schemes} 

As an illustrative example we consider the c-optimality criterion with $\mathbf{c} = (1, 0)\T$, i.e., minimising the variance of estimating the location parameter $\eta_0$. Since the optimality criterion is linear here, the optimal sampling scheme can be found according to Algorithm \ref{alg:lopt} with 
\begin{gather*}   
    \mathbf{L} = (1, 0)\T, \quad  
    \mathbf{H}(\boldsymbol{\theta}_0) =
    \frac{1}{\sigma_0^2}
    \begin{pmatrix}
    1 & 0 \\
    0 & 2
    \end{pmatrix}
    , \quad 
    \boldsymbol{\psi}_i(\boldsymbol{\theta}_0) = - w_i\left(\frac{\log y_i - \eta_0}{\sigma_0^2}, \frac{(\log y_i - \eta_0)^2}{\sigma_0^3} - \frac{1}{\sigma_0}\right)\T,
\end{gather*}
and
\begin{align}
\label{eq:copt_ci}
c_i = \bigr\rvert\bigr\rvert
\mathbf{L}\T
\mathbf{H}(\boldsymbol{\theta}_0)^{-1}
\boldsymbol{\psi}_i(\boldsymbol{\theta}_0)
\bigr\rvert\bigr\rvert^2_{2} \propto w_i^2(\log y_i - \eta_0)^2.
\end{align}

To find an optimal sampling scheme with respect to the anticipated variance of our estimator for $\eta_0$, we replace $y_i$ by a random variable $Y_i$ and evaluate the corresponding expectation of \eqref{eq:copt_ci} under an assumed model for $Y_i$. If we assume that $\log Y_i$ has mean $\hat y_i$ and variance $\sigma_i^2$, we obtain 
\begin{equation*}
    \tilde c_i \propto \sqrt{\E_{Y_i}[w_i^2(\log Y_i - \eta_0)^2]} = w_i \sqrt{(\hat y_i - \eta_0)^2 + \sigma_i^2},
\end{equation*}
which in practice may be evaluated at a preliminary estimate $\tilde{\eta}_0$ of $\eta_0$. The predictions $\hat y_i$ and dispersion parameters $\sigma_i^2$ may be modelled as functions of the observed auxiliary variables (i.e., the case identifier, off-road glance duration, and maximal deceleration during braking), and estimated from a pilot sample or using sequential subsampling methods (Algorithm \ref{alg:sequential_optimal_designs}). The resulting sampling scheme is guaranteed to produce strictly positive sampling probabilities as long as all $w_i, \sigma_i>0$.

\paragraph{Results} 

The computation time, number of iterations needed for convergence, and relative efficiencies for various optimality criteria are presented in Table \ref{tab:baseline_impact_speed}. The D-optimal sampling scheme was found in four fixed-point iterations with Algorithm \ref{alg:linearisation}. The $\Phi_{0.5}$-, $\Phi_{5}$- optimal sampling schemes were found in three and 25 iterations, respectively. An E-optimal sampling scheme could not be found, due to the non-convexity of the objective function. The computation time for finding an optimal sampling scheme ranged from 0.10 seconds for the linear optimality criteria to 0.96 s for the D-optimality criterion and 4.95 s for the $\Phi_5$-optimality criterion. The optimal sampling schemes of the $d_{\mathrm{ER}}$- and $d_{\mathrm{S}}$-optimality criteria reached 97--99\% A-efficiency, 96--99\% D-efficiency, and 92--94\% $\Phi_5$-efficiency. The A-optimal sampling scheme had a similar performance.

\begin{table}[htb!]
\centering
\caption{Performance measures for estimating the log-normal model \eqref{eq:lognormal} by optimal subsampling with various optimality criteria. The columns show the number of fixed-point iterations and execution time to find the optimal sampling scheme, and relative efficiencies with respect to other optimality criteria.} 
\label{tab:baseline_impact_speed}
\begin{tabular}{lrrrrrrrr}
 Optimality criterion & No. iterations & Time (s) & A-eff & c$_{(1,0)}$-eff & c$_{(0,1)}$-eff & D-eff & $d_{\mathrm{ER}}$-eff & $\Phi_5$-eff \\ 
  \hline
A &  & 0.11 & \textbf{1.00} & 0.90 & 0.76 & 0.99 & 0.99 & 0.96 \\ 
  c, $\mathbf{c} = (1,0)\T$ &  & 0.10 & 0.14 & \textbf{1.00} & 0.04 & 0.24 & 0.10 & 0.09 \\ 
  c, $\mathbf{c} = (0,1)\T$ &  & 0.10 & 0.27 & 0.16 & \textbf{1.00} & 0.47 & 0.34 & 0.19 \\ 
  D & 4 & 0.96 & 0.98 & 0.86 & 0.79 & \textbf{1.00} & 0.99 & 0.91 \\ 
  $d_{\mathrm{ER}}$ &  & 0.10 & 0.99 & 0.84 & 0.83 & 0.99 & \textbf{1.00} & 0.92 \\ 
  $d_{\mathrm{S}}$ &  & 0.10 & 0.97 & 0.84 & 0.79 & 0.96 & 0.98 & 0.94 \\ 
  E & Diverged & - &  &  &  &  &  &  \\ 
  $\Phi_{0.5}$ & 3 & 0.78 & $>$0.99 & 0.88 & 0.78 & $>$0.99 & 0.99 & 0.94 \\ 
  $\Phi_5$ & 25 & 4.95 & 0.95 & 0.91 & 0.65 & 0.90 & 0.91 & \textbf{1.00} \\ 
  $\Phi_{10}$ & Diverged & - &  &  &  &  &  &  \\ 
   \hline
\end{tabular}
\end{table}

\subsection{Regression modelling}
\label{sec:regression_modelling}

Next we consider the distribution of the baseline impact speed as a function of the input variables to the scenario generation, i.e., the off-road glance duration and maximal deceleration during braking.

We first note that the impact speed increases monotonically with increased levels of the off-road glance duration and decreased levels of deceleration. Hence, variations generated from the same original rear-end crash have an upper bound on their impact speed, attained for the variation having the off-road glance duration at its maximum and the deceleration level at its minimum. We assume that this maximal impact speed is known, e.g., observed by running the corresponding virtual simulation. The impact speed may then be expressed relative to the maximal impact speed for that specific case, with values in the common range $[0, 1]$. Note that in this case the explanatory variables are known \textit{a priori}, whereas the outcome (i.e., relative impact speed) can only be observed after running the corresponding virtual simulation.

\paragraph{Model}

A simple model for a response variable on the unit interval is a quasi-binomial logistic regression model, for which the full-data parameter $\boldsymbol{\theta}_0$ is defined as 
\begin{equation}
    \label{eq:qblr}
    \boldsymbol{\theta}_0 = \argmin_{\boldsymbol{\theta} \in \mathbb{R}^p} -\sum_{i=1}^N y_i \log p_i(\boldsymbol{\theta}) + (1 - y_i)\log(1 -p_i(\boldsymbol{\theta})), \quad p_i(\boldsymbol{\theta}) = (1 + \exp(-\boldsymbol{x}_i\T\boldsymbol{\theta}))^{-1},    
\end{equation}
where $\boldsymbol{x}_i$ is a feature vector pertaining to instance $i$, and $\boldsymbol{\theta}$ a vector of regression coefficients. As explanatory variables we include the case identifier of the original rear-end crash event (categorical with 44 levels, dummy coded into 44 binary variables), the off-road glance duration, the maximal deceleration during braking, and all three-way interactions. For each of the 44 cases, the impact speed response surface is then described by 4 parameters: an intercept parameter and three slope parameters corresponding to the off-road glance duration, deceleration level, and the interaction between those. The joint parameter vector $\boldsymbol{\theta}$ is of dimension $44 \times 4 = 176$. Note that in this case we do not include the observation weights $w_i$ in the empirical risk function, since these are functions of the explanatory variables and hence ignorable in this context. Illustrations of the observed and predicted impact speed response surfaces for three of the cases are presented in Figure \ref{fig:response_surface} in Appendix \ref{sec:supplementary_figures}.

\paragraph{Optimal sampling schemes} 

For illustrative purposes, we consider the $d_{\mathrm{ER}}$-optimality criterion. Since this is a linear optimality criterion, the optimal sampling scheme can be found according to Algorithm \ref{alg:lopt} with 
\begin{align*}
\mathbf{L} = \mathbf{H}(\boldsymbol{\theta}_0)^{1/2}, \quad  
\mathbf{H}(\boldsymbol{\theta}_0) =
\mathbf{X}\T\mathbf{W}(\boldsymbol{\theta}_0)\mathbf{X}, \quad
\boldsymbol{\psi}_i(\boldsymbol{\theta}_0)
= -(y_i - p_i(\boldsymbol{\theta}_0))\boldsymbol{x}_i,
\end{align*}
where $\mathbf{W}(\boldsymbol{\theta})$ is the diagonal matrix with entries $p_i(\boldsymbol{\theta})(1-p_i(\boldsymbol{\theta}))$ and $\mathbf{X}$ the matrix with rows $\boldsymbol{x}_i\T$, and 
\begin{align}
    \label{eq:ci_response}
    c_i = ||\mathbf{L}\T\mathbf{H}(\boldsymbol{\theta}_0)^{-1}\boldsymbol{\psi}_i(\boldsymbol{\theta}_0)||_2^2 = (y_i - p_i(\boldsymbol{\theta}_0))^2 \boldsymbol{x}_i\T(\mathbf{X}\T\mathbf{W}(\boldsymbol{\theta}_0)\mathbf{X})^{-1}\boldsymbol{x}_i.
\end{align}

To find a d$_{ER}$-optimal sampling scheme with respect to the anticipated covariance matrix, we replace $y_i$ in \eqref{eq:ci_response} by a random variable $Y_i$ and evaluate the corresponding expectation under a model for $Y_i$ given the known explanatory variables $\boldsymbol{x}_i$. For instance, we may assume that $Y_i$ has mean $p(\boldsymbol{\theta}_0)$ and variance $p_i(\boldsymbol{\theta}_0)(1-p_i(\boldsymbol{\theta}_0))$. We then obtain 
\begin{equation}
    \label{eq:leverage_sampling}
    \tilde c_i = \sqrt{\E_{Y_i}[(Y_i - p_i(\boldsymbol{\theta}_0))^2 \boldsymbol{x}_i\T (\mathbf{X}\T\mathbf{W}(\boldsymbol{\theta}_0)\mathbf{X})^{-1}\boldsymbol{x}_i]} = \sqrt{h_{ii}(\boldsymbol{\theta}_0)}, 
\end{equation}
where $h_{ii}(\boldsymbol{\theta}_0)$ is the $i^{\mathrm{th}}$ diagonal element of the 'hat matrix', or projection matrix 
\[
\mathbf{W}(\boldsymbol{\theta}_0)^{1/2}
\mathbf{X}
(\mathbf{X}\T\mathbf{W}(\boldsymbol{\theta}_0)\mathbf{X})^{-1}
\mathbf{X}\T
\mathbf{W}(\boldsymbol{\theta}_0)^{1/2}
\]
\citep[see][]{Hoaglin1978, Pregibon1981}. To account for the influence a data point $(\boldsymbol{x}_i, y_i)$ exerts on its own prediction, it is appropriate to deflate the variance of $Y_i$ by a factor $1-h_{ii}(\boldsymbol{\theta}_0)$, resulting in
\begin{equation*}
    \label{eq:lev2}
    \tilde c_i = \sqrt{h_{ii}(\boldsymbol{\theta}_0)(1-h_{ii}(\boldsymbol{\theta}_0))}
\end{equation*}
instead of \eqref{eq:leverage_sampling} \citep[cf.][]{Ma2020}. 
In practice we may evaluate $\tilde{c}_i$ at a preliminary estimate $\tilde{\boldsymbol{\theta}}_0$ obtained from a pilot sample or estimated using sequential subsampling methods (Algorithm \ref{alg:sequential_optimal_designs}). The resulting sampling scheme is guaranteed to produce strictly positive sampling probabilities as long as the predictions $p_i(\tilde{\boldsymbol{\theta}}_0)$ are bounded away from $0$ and $1$.

\paragraph{Results} 

Table \ref{tab:impact_speed_response_surface} shows the computation time, relative efficiencies, and number of iterations needed to find an optimal sampling scheme for various optimality criteria. Optimal sampling schemes were found in five fixed-point iterations for D-optimality, four iterations for $\Phi_{0.5}$-optimality, and could not be found for the $\Phi_5$-, $\Phi_{10}$- and E-optimality criteria. Finding an L-optimal sampling scheme required 95\% less computation time than for the non-linear D-optimality criterion. The $d_{\mathrm{ER}}$- and $d_{\mathrm{S}}$-optimal schemes attained 40--47\% A-efficiency and 92--96\% D-efficiency. The A-optimal sampling scheme had only 60\% D-efficiency. The $\Phi_{0.5}$-optimal sampling scheme, which interpolates between A- and D-optimality, achieved 92\% A-efficiency and 80\% D-efficiency.

\begin{table}[htb!]
\centering
\caption{Performance measures for estimating the quasi-binomial logistic regression model \eqref{eq:qblr} by optimal subsampling with various optimality criteria. The columns show the number of fixed-point iterations and execution time to find the optimal sampling scheme, and relative efficiencies with respect to other optimality criteria. The computation time for fitting the model to the full dataset was 8.46 seconds.} 
\label{tab:impact_speed_response_surface}
\begin{tabular}{lrrrrrrr}
 Optimality criterion & No. iterations & Time (s) & A-eff & D-eff & $d_{\mathrm{ER}}$-eff & $d_{\mathrm{S}}$-eff & $\Phi_{0.5}$-eff \\ 
  \hline
A &  & 1.16 & \textbf{1.00} & 0.60 & 0.47 & 0.42 & 0.93 \\ 
  D & 5 & 27.69 & 0.49 & \textbf{1.00} & 0.89 & 0.94 & 0.77 \\ 
  $d_{\mathrm{ER}}$ &  & 1.11 & 0.47 & 0.92 & \textbf{1.00} & 0.91 & 0.71 \\ 
  $d_{\mathrm{S}}$ &  & 1.12 & 0.40 & 0.96 & 0.92 & \textbf{1.00} & 0.67 \\ 
  E & Diverged & - &  &  &  &  &  \\ 
  $\Phi_{0.5}$ & 4 & 20.37 & 0.92 & 0.80 & 0.68 & 0.65 & \textbf{1.00} \\ 
  $\Phi_5$ & Diverged & - &  &  &  &  &  \\ 
  $\Phi_{10}$ & Diverged & - &  &  &  &  &  \\ 
   \hline
\end{tabular}
\end{table}

\subsection{Finite population inference}
\label{sec:finite_population_inference}

We finally consider the potential safety benefit of the AEB system compared to a baseline manual driving scenario. For a scenario $i \in \mathcal D$, let $\boldsymbol{y}_i = (y_{i1}, y_{i2}, y_{i3})\T$, where $y_{1i}$ is the impact speed reduction, $y_{i2}$ the injury risk reduction, and $y_{i3}$ the binary crash avoidance indicator with the AEB system compared to baseline manual driving. The distributions of these characteristics are illustrated in Figure \ref{fig:descriptives} in Appendix \ref{sec:supplementary_figures}.

\paragraph{Model} We are interested in the mean impact speed reduction, mean injury risk reduction and crash avoidance rate, given by the vector total
\begin{equation*}
    \boldsymbol{t}_{\boldsymbol{y}}  = \sum_{i=1}^N w_i \boldsymbol{y}_i,
\end{equation*}
where the observation weights $w_i$ are normalised so that $\sum_{i=1}^N w_i = 1$. This can also be expressed as
\begin{equation}
    \label{eq:theta0_fps}
    \boldsymbol{t}_{\boldsymbol{y}} = \boldsymbol{\theta}_0 = \argmin_{\boldsymbol{\theta} \in \mathbb{R}^3} \frac{1}{2}\sum_{i=1}^N w_i ||\boldsymbol{y}_i - \boldsymbol{\theta}||_2^2.   
\end{equation}

\paragraph{Optimal sampling schemes} 

As an example, consider the $d_{\mathrm{S}}$-optimality criterion. Since this is a linear optimality criterion, the optimal sampling scheme can be found according to Algorithm \ref{alg:lopt} with 
\[
\mathbf{L} = \mathbf{H}(\boldsymbol{\theta}_0)\mathbf{V}(\boldsymbol{\theta}_0)^{-1/2}, \quad 
\mathbf{H}(\boldsymbol{\theta}_0) = \mathbf{I}_{3 \times 3}, \quad \mathbf{V}(\boldsymbol{\theta}_0) = \sum_{i \in \mathcal D} \boldsymbol{\psi}_i(\boldsymbol{\theta}_0)\boldsymbol{\psi}_i(\boldsymbol{\theta}_0)^T, \quad
\boldsymbol{\psi}_i(\boldsymbol{\theta}_0) = -w_i(\boldsymbol{y}_i - \boldsymbol{\theta}_0),
\]
and
\begin{alignat}{2}
    \label{eq:ci_fps}
    c_i & = 
    ||\mathbf{L}\T\mathbf{H}(\boldsymbol{\theta}_0)^{-1}\boldsymbol{\psi}_i(\boldsymbol{\theta}_0)||_2^2 = 
    w_i^2(\boldsymbol{y}_i - \boldsymbol{\theta}_0)\T\mathbf{V}(\boldsymbol{\theta}_0)^{-1}(\boldsymbol{y}_i - \boldsymbol{\theta}_0).
\end{alignat}

In order to find the L-optimal sampling scheme with respect to the anticipated covariance matrix, we introduce a random vector $\mathbf{Y}_i$, substitute $\mathbf{Y}_i$ for $\boldsymbol{y}_i$ in \eqref{eq:ci_fps}, and evaluate the expectation. Let therefore $\hat{\boldsymbol{y}}_i$ and $\hat{\boldsymbol{\Sigma}}_i$ denote the mean vector and covariance matrix of $\mathbf{Y}_i$, respectively. By properties of quadratic forms \citep{Mathai1992}, we obtain
\begin{align}
    \label{eq:citilde_fps}
    \tilde c_i 
    & = \sqrt{\E\left[w_i^2(\boldsymbol{y}_i - \boldsymbol{\theta}_0)\T\mathbf{V}(\boldsymbol{\theta}_0)^{-1}(\boldsymbol{y}_i - \boldsymbol{\theta}_0)\right]} = w_i\sqrt{\left[(\hat{\boldsymbol{y}}_i - \boldsymbol{\theta}_0)\T\mathbf{V}(\boldsymbol{\theta}_0)^{-1} (\hat{\boldsymbol{y}}_i - \boldsymbol{\theta}_0) + \tr(\mathbf{V}(\boldsymbol{\theta}_0)^{-1}\hat{\boldsymbol{\Sigma}}_i)\right]}.
\end{align}
To implement optimal sampling in practice, we evaluate \eqref{eq:citilde_fps} at a preliminary estimate $\tilde{\boldsymbol{\theta}}_0$ obtained from a pilot sample. The predictions $\hat{\boldsymbol{y}}_i$ and dispersion matrices $\hat{\boldsymbol{\Sigma}}_i$ may be modelled as functions of the observed auxiliary variables (i.e., the case identifier, off-road glance duration, and maximal deceleration during braking) and iteratively updated using sequential subsampling methods (Algorithm \ref{alg:sequential_optimal_designs}). The sampling scheme derived from \eqref{eq:citilde_fps} is guaranteed to produce strictly positive sampling probabilities as long as all $\hat{\boldsymbol{\Sigma}}_i$ are full-rank.

\paragraph{Results}

Results in terms of computation time, number of iterations needed for convergence, and relative efficiencies of various optimality criteria are presented in Table \ref{tab:finite_population_inference}. The optimal sampling scheme was found in four fixed-point iterations for the D-optimality criterion, and in two iterations for the other non-linear optimality criteria. The computation time ranged from 0.10 for the linear optimality criteria, to 0.94 s for the D-optimality criterion. The $d_{\mathrm{ER}}$-optimal sampling scheme had 100\% A-efficiency, 46\% D-efficiency and $>$99\% E-efficiency. In fact, in this case the $d_{\mathrm{ER}}$-optimality criterion is identical to A-optimality. In contrast, the $d_{\mathrm{S}}$--optimal sampling scheme had 73\% A-efficiency, 98\% D-efficiency, and 72\% E-efficiency. The $\Phi_{0.5}$-criterion had 99\% A-efficiency, 58\% D-efficiency and 99\% E-efficiency. The A- and E-optimality criteria were largely driven by the mean impact speed reduction, as this was measured on a scale that was orders of magnitude larger than the measurement-scale for the injury risk reduction and crash avoidance (Figure \ref{fig:descriptives}, Appendix \ref{sec:supplementary_figures}).

\begin{table}[htb!]
\centering
\caption{Performance measures for estimating the vector of finite population means \eqref{eq:theta0_fps} by optimal subsampling with various optimality criteria. The columns show the number of fixed-point iterations and execution time to find the optimal sampling scheme, and relative efficiencies with respect to other optimality criteria.} 
\label{tab:finite_population_inference}
\resizebox{\textwidth}{!}{\begin{tabular}{lrrrrrrrr}
 Optimality criterion & No. iterations & Time (s) & A-eff & c$_{(1,0,0)}$-eff & c$_{(0,1,0)}$-eff & c$_{(0,0,1)}$-eff & D-eff & E-eff \\ 
  \hline
A &  & 0.10 & \textbf{1.00} & $>$0.99 & 0.36 & 0.25 & 0.46 & $>$0.99 \\ 
  c, $\mathbf{c} = (1,0,0)\T$ &  & 0.10 & 0.98 & \textbf{1.00} & 0.05 & 0.04 & 0.20 & $>$0.99 \\ 
  c, $\mathbf{c} = (0,1,0)\T$ &  & 0.10 & 0.12 & 0.12 & \textbf{1.00} & 0.11 & 0.22 & 0.12 \\ 
  c, $\mathbf{c} = (0,0,1)\T$ &  & 0.10 & 0.41 & 0.41 & 0.50 & \textbf{1.00} & 0.70 & 0.41 \\ 
  D & 4 & 0.94 & 0.65 & 0.65 & 0.76 & 0.82 & \textbf{1.00} & 0.65 \\ 
  $d_{\mathrm{ER}}$ &  & 0.10 & \textbf{1.00} & $>$0.99 & 0.36 & 0.25 & 0.46 & $>$0.99 \\ 
  $d_{\mathrm{S}}$ &  & 0.10 & 0.73 & 0.72 & 0.77 & 0.77 & 0.98 & 0.72 \\ 
  E & 2 & 0.70 & 0.99 & $>$0.99 & 0.07 & 0.06 & 0.22 & \textbf{1.00} \\ 
  $\Phi_{0.5}$ & 2 & 0.57 & 0.99 & 0.99 & 0.46 & 0.38 & 0.58 & 0.99 \\ 
  $\Phi_5$ & 2 & 0.58 & 0.99 & $>$0.99 & 0.07 & 0.06 & 0.22 & \textbf{1.00} \\ 
  $\Phi_{10}$ & 2 & 0.57 & 0.99 & $>$0.99 & 0.07 & 0.06 & 0.22 & \textbf{1.00} \\ 
   \hline
\end{tabular}}
\end{table}

\section{Discussion}
\label{sec:discussion}

We have presented a theory of optimal subsampling design for a general class of estimators, sampling designs, and optimality criteria. Although the presented optimality conditions are valid for any differentiable objective function, the algorithms for finding optimal sampling schemes are most appropriate for convex functions. Further research could include development of methods to handle non-convex optimality criteria, such as E- and G-optimality \citep{Kiefer1960, Kiefer1974}.

From an applied perspective, we believe that the proposed invariant linear optimality criteria (i.e., $d_{\mathrm{ER}}$-, $d_{\mathrm{KL}}$- and $d_{\mathrm{S}}$-optimality) offer a good compromise between computational and statistical efficiency. Non-linear optimality criteria require iterative procedures and computationally expensive covariance matrix evaluations, which limits their usability in problems and applications where computational complexity is a major concern. Further studies evaluating the performance of these methods in practice and in other applications are encouraged.

Sequential subsampling is a viable approach to implement optimal subsampling methods in practice. The theoretical properties of the estimators derived from such sequential subsampling methods have so far only been studied rigorously in limited settings. Further research in this direction is requested.

\section*{Acknowledgement}

We would like to thank Malin Svärd and Simon Lundell at Volvo Car Corporation for allowing us to use their data in our experiments.

\bibliographystyle{apalike}
\bibliography{references}

\begin{thebibliography}{}

\bibitem[Ai et~al., 2021a]{Ai2021_quantile}
Ai, M., Wang, F., Yu, J., and Zhang, H. (2021a).
\newblock Optimal subsampling for large-scale quantile regression.
\newblock {\em Journal of Complexity}, 62:101512.

\bibitem[Ai et~al., 2021b]{Ai2021_regression}
Ai, M., Yu, J., Zhang, H., and Wang, H. (2021b).
\newblock Optimal subsampling algorithms for big data regressions.
\newblock {\em Statistica Sinica}.

\bibitem[Anderson et~al., 2013]{Anderson2013}
Anderson, R., Doecke, S., Mackenzie, J., and Ponte, G. (2013).
\newblock Potential benefits of autonomous emergency braking based on in-depth
  crash reconstruction and simulation.
\newblock In {\em Proceedings of the 23rd International Conference on Enhanced
  Safety of Vehicles}.

\bibitem[Atkinson and Donev, 1992]{Atkinson1992}
Atkinson, A.~C. and Donev, A.~N. (1992).
\newblock {\em {Optimum Experimental Designs}}.
\newblock Clarendon Press, Oxford.

\bibitem[Bach, 2007]{Bach2007}
Bach, F.~R. (2007).
\newblock {Active learning for misspecified generalized linear models}.
\newblock In {\em {Advances in Neural Information Processing Systems 19}}.

\bibitem[Bellhouse, 1984]{Bellhouse1984}
Bellhouse, D.~R. (1984).
\newblock A review of optimal designs in survey sampling.
\newblock {\em Canadian Journal of Statistics}, 12(1):53--65.

\bibitem[Binder, 1983]{Binder1983}
Binder, D.~A. (1983).
\newblock On the variances of asymptotically normal estimators from complex
  surveys.
\newblock {\em {International Statistical Review}}, 51(3):279--292.

\bibitem[Boyd and Vandenberghe, 2004]{Boyd2004}
Boyd, S. and Vandenberghe, L. (2004).
\newblock {\em Convex Optimization}.
\newblock Cambridge University Press, Cambridge.

\bibitem[Brewer, 1979]{Brewer1979}
Brewer, K. R.~W. (1979).
\newblock A class of robust sampling designs for large-scale surveys.
\newblock {\em Journal of the American Statistical Association},
  74(368):911--915.

\bibitem[Casella and Berger, 2001]{Casella2002}
Casella, G. and Berger, R. (2001).
\newblock {\em Statistical Inference}.
\newblock {Duxbury}, Pacific Grove.

\bibitem[Cassel et~al., 1976]{Cassel1976}
Cassel, C.~M., Särndal, C.~E., and Wretman, J.~H. (1976).
\newblock {Some results on generalized difference estimation and generalized
  regression estimation for finite populations}.
\newblock {\em Biometrika}, 63(3):615--620.

\bibitem[Dai et~al., 2022]{Dai2022}
Dai, W., Song, Y., and Wang, D. (2022).
\newblock A subsampling method for regression problems based on minimum energy
  criterion.
\newblock {\em Technometrics}.
\newblock Advance online publication.
  \url{https://doi.org/10.1080/00401706.2022.2127915}.

\bibitem[Deldossi and Tommasi, 2022]{Deldossi2022}
Deldossi, L. and Tommasi, C. (2022).
\newblock Optimal design subsampling from big datasets.
\newblock {\em Journal of Quality Technology}, 54(1):93--101.

\bibitem[Drovandi et~al., 2017]{Drovandi2017}
Drovandi, C.~C., Holmes, C.~C., McGree, J.~M., Mengersen, K., Richardson, S.,
  and Ryan, E.~G. (2017).
\newblock {Principles of Experimental Design for Big Data Analysis}.
\newblock {\em Statistical Science}, 32(3):385--404.

\bibitem[Efron and Hinkley, 1978]{Efron1978}
Efron, B. and Hinkley, D.~V. (1978).
\newblock Assessing the accuracy of the maximum likelihood estimator: Observed
  versus expected fisher information.
\newblock {\em Biometrika}, 65(3):457--482.

\bibitem[Fuller, 2009]{Fuller2009}
Fuller, W.~A. (2009).
\newblock {\em {Sampling Statistics}}.
\newblock Wiley, Hoboken.

\bibitem[H{\'a}jek, 1981]{Hajek1981}
H{\'a}jek, J. (1981).
\newblock {\em Sampling from a Finite Population}.
\newblock Marcel Dekker, New York.

\bibitem[Hall and Heyde, 1980]{Hall1980}
Hall, P. and Heyde, C. (1980).
\newblock {\em Martingale Limit Theory and Its Application}.
\newblock Academic Press, New York.

\bibitem[Hansen and Hurwitz, 1943]{Hansen1943}
Hansen, M.~H. and Hurwitz, W.~N. (1943).
\newblock On the theory of sampling from finite populations.
\newblock {\em The Annals of Mathematical Statistics}, 14(4):333--362.

\bibitem[Hartley and Sielken, 1975]{Hartley1975}
Hartley, H.~O. and Sielken, R.~L. (1975).
\newblock A "super-population viewpoint" for finite population sampling.
\newblock {\em Biometrics}, 31(2):411--422.

\bibitem[Hastie, 2020]{Hastie2020}
Hastie, T. (2020).
\newblock Ridge regularization: An essential concept in data science.
\newblock {\em Technometrics}, 62(4):426--433.

\bibitem[Hoaglin and Welsch, 1978]{Hoaglin1978}
Hoaglin, D.~C. and Welsch, R.~E. (1978).
\newblock The hat matrix in regression and {ANOVA}.
\newblock {\em The American Statistician}, 32(1):17--22.

\bibitem[Horn and Johnson, 1990]{Horn1990}
Horn, R. and Johnson, C. (1990).
\newblock {\em Matrix Analysis}.
\newblock Cambridge University Press, Cambridge.

\bibitem[Horvitz and Thompson, 1952]{Horvitz1952}
Horvitz, D.~G. and Thompson, D.~J. (1952).
\newblock A generalization of sampling without replacement from a finite
  universe.
\newblock {\em Journal of the American Statistical Association},
  47(260):663--685.

\bibitem[Hájek, 1959]{Hajek1959}
Hájek, J. (1959).
\newblock Optimal strategy and other problems in probability sampling.
\newblock {\em Časopis pro pěstování matematiky}, 84(4):387--423.

\bibitem[Imberg et~al., 2020]{Imberg2020}
Imberg, H., Jonasson, J., and Axelson-Fisk, M. (2020).
\newblock Optimal sampling in unbiased active learning.
\newblock In {\em Proceedings of the 23rd International Conference on
  Artificial Intelligence and Statistics}.

\bibitem[Imberg et~al., 2022a]{Imberg2022a}
Imberg, H., Lisovskaja, V., Selpi, and Nerman, O. (2022a).
\newblock Optimization of two-phase sampling designs with application to
  naturalistic driving studies.
\newblock {\em IEEE Transactions on Intelligent Transportation Systems},
  23(4):3575--3588.

\bibitem[Imberg et~al., 2022b]{Imberg2022b}
Imberg, H., Yang, X., Flannagan, C., and Bärgman, J. (2022b).
\newblock Active sampling: A machine-learning-assisted framework for finite
  population inference with optimal subsamples.
\newblock {\tt arXiv:2212.10024 [stat.ME]}.

\bibitem[Isaki and Fuller, 1982]{Isaki1982}
Isaki, C.~T. and Fuller, W.~A. (1982).
\newblock Survey design under the regression superpopulation model.
\newblock {\em Journal of the American Statistical Association},
  77(377):89--96.

\bibitem[Kiefer, 1974]{Kiefer1974}
Kiefer, J. (1974).
\newblock {General Equivalence Theory for Optimum Designs (Approximate
  Theory)}.
\newblock {\em The Annals of Statistics}, 2(5):849--879.

\bibitem[Kiefer and Wolfowitz, 1960]{Kiefer1960}
Kiefer, J. and Wolfowitz, J. (1960).
\newblock The equivalence of two extremum problems.
\newblock {\em Canadian Journal of Mathematics}, 12:363--366.

\bibitem[Kossen et~al., 2022]{Kossen2022}
Kossen, J., Farquhar, S., Gal, Y., and Rainforth, T. (2022).
\newblock Active surrogate estimators: An active learning approach to
  label-efficient model evaluation.
\newblock In {\em Advances in Neural Information Processing Systems}.

\bibitem[Kullback and Leibler, 1951]{Kullback1951}
Kullback, S. and Leibler, R.~A. (1951).
\newblock {On Information and Sufficiency}.
\newblock {\em The Annals of Mathematical Statistics}, 22(1):79--86.

\bibitem[Ma et~al., 2015]{Ma2015}
Ma, P., Mahoney, M.~W., and Yu, B. (2015).
\newblock A statistical perspective on algorithmic leveraging.
\newblock {\em Journal of Machine Learning Research}, 16:861--911.

\bibitem[Ma et~al., 2020]{Ma2020}
Ma, P., Zhang, X., Xing, X., Ma, J., and Mahoney, M.~W. (2020).
\newblock Asymptotic analysis of sampling estimators for randomized numerical
  linear algebra algorithms.
\newblock In {\em Proceedings of the 23rd International Conference on
  Artificial Intelligence and Statistics}.

\bibitem[Mathai and Provost, 1992]{Mathai1992}
Mathai, A. and Provost, S. (1992).
\newblock {\em Quadratic Forms in Random Variables}.
\newblock Marcel Dekker, New York.

\bibitem[McCullagh and Nelder, 1989]{McCullagh1989}
McCullagh, P. and Nelder, J.~A. (1989).
\newblock {\em {Generalized Linear Models}}.
\newblock CRC Press, Boca Raton.

\bibitem[Meng et~al., 2021]{Meng2021}
Meng, C., Xie, R., Mandal, A., Zhang, X., Zhong, W., and Ma, P. (2021).
\newblock {LowCon: A design-based subsampling approach in a misspecified linear
  model}.
\newblock {\em Journal of Computational and Graphical Statistics},
  30(3):694--708.

\bibitem[Mullins et~al., 2018]{Mullins2018}
Mullins, G.~E., Stankiewicz, P.~G., Hawthorne, R.~C., and Gupta, S.~K. (2018).
\newblock Adaptive generation of challenging scenarios for testing and
  evaluation of autonomous vehicles.
\newblock {\em Journal of Systems and Software}, 137:197--215.

\bibitem[Nelder and Wedderburn, 1972]{Nelder1972}
Nelder, J.~A. and Wedderburn, R. W.~M. (1972).
\newblock Generalized linear models.
\newblock {\em Journal of the Royal Statistical Society. Series A (General)},
  135(3):370--384.

\bibitem[Neyman, 1938]{Neyman1938}
Neyman, J. (1938).
\newblock Contribution to the theory of sampling human populations.
\newblock {\em Journal of the American Statistical Association},
  33(201):101--116.

\bibitem[Petersen and Pedersen, 2012]{Petersen2012}
Petersen, K.~B. and Pedersen, M.~S. (2012).
\newblock \textit{The Matrix Cookbook}.

\bibitem[Pregibon, 1981]{Pregibon1981}
Pregibon, D. (1981).
\newblock {Logistic regression diagnostics}.
\newblock {\em The Annals of Statistics}, 9(4):705--724.

\bibitem[Pronzato and P{\'a}zman, 2013]{Pronzato2013}
Pronzato, L. and P{\'a}zman, A. (2013).
\newblock {\em Design of Experiments in Nonlinear Models}.
\newblock Springer, New York.

\bibitem[Pukelsheim, 1993]{Pukelsheim1993}
Pukelsheim, F. (1993).
\newblock {\em Optimal Design of Experiments}.
\newblock Wiley, New York.

\bibitem[{R Core Team}, 2023]{RCoreTeam}
{R Core Team} (2023).
\newblock {\em R: A Language and Environment for Statistical Computing}.
\newblock R Foundation for Statistical Computing, Vienna.

\bibitem[Settles, 2012]{Settles2012}
Settles, B. (2012).
\newblock Active learning.
\newblock {\em Synthesis Lectures on Artificial Intelligence and Machine
  Learning}, 6(1):1--114.

\bibitem[Seyedi et~al., 2021]{Seyedi2021}
Seyedi, M., Koloushani, M., Jung, S., and Vanli, A. (2021).
\newblock Safety assessment and a parametric study of forward
  collision-avoidance assist based on real-world crash simulations.
\newblock {\em Journal of Advanced Transportation}.
\newblock Advance online publication.
  \url{https://doi.org/10.1155/2021/4430730}.

\bibitem[Sibson, 1974]{Sibson1974}
Sibson, R. (1974).
\newblock {D$_A$}-optimality and duality.
\newblock In {\em Progress of Statistics, Volume 2: Proceedings of the 9th
  European Meeting of Statisticians}.

\bibitem[Silvey, 1980]{Silvey1980}
Silvey, S. (1980).
\newblock {\em Optimal Design}.
\newblock Chapman \& Hall, London.

\bibitem[Stefanski and Boos, 2002]{Stefanski2002}
Stefanski, L.~A. and Boos, D.~D. (2002).
\newblock The calculus of {M}-estimation.
\newblock {\em The American Statistician}, 56(1):29--38.

\bibitem[Sun et~al., 2022]{Sun2022}
Sun, J., Zhou, H., Xi, H., Zhang, H., and Tian, Y. (2022).
\newblock Adaptive design of experiments for safety evaluation of automated
  vehicles.
\newblock {\em IEEE Transactions on Intelligent Transportation Systems},
  23(9):14497--14508.

\bibitem[Till{\'{e}}, 2006]{Tille2006}
Till{\'{e}}, Y. (2006).
\newblock {\em {Sampling Algorithms}}.
\newblock Springer, New York.

\bibitem[Vapnik, 1991]{Vapnik1991}
Vapnik, V. (1991).
\newblock Principles of risk minimization for learning theory.
\newblock In {\em Proceedings of the 4th International Conference on Neural
  Information Processing Systems}.

\bibitem[Wang and Ma, 2021]{Wang2021}
Wang, H. and Ma, Y. (2021).
\newblock {Optimal subsampling for quantile regression in big data}.
\newblock {\em Biometrika}, 108(1):99--112.

\bibitem[Wang et~al., 2018]{Wang2018}
Wang, H., Zhu, R., and Ma, P. (2018).
\newblock Optimal subsampling for large sample logistic regression.
\newblock {\em Journal of the American Statistical Association},
  113(522):829--844.

\bibitem[Wang et~al., 2017]{Wang2017}
Wang, Y., Yu, A.~W., and Singh, A. (2017).
\newblock On computationally tractable selection of experiments in
  measurement-constrained regression models.
\newblock {\em Journal of Machine Learning Research}, 18(143):1--41.

\bibitem[Wedderburn, 1974]{Wedderburn1974}
Wedderburn, R. W.~M. (1974).
\newblock {Quasi-likelihood functions, generalized linear models, and the
  Gauss-Newton method}.
\newblock {\em Biometrika}, 61(3):439--447.

\bibitem[Welch, 1984]{Welch1984}
Welch, W.~J. (1984).
\newblock Computer-aided design of experiments for response estimation.
\newblock {\em Technometrics}, 26(3):217--224.

\bibitem[{World Health Organization}, 2018]{WHO2018}
{World Health Organization} (2018).
\newblock {\em Global status report on road safety 2018}.
\newblock \url{https://www.who.int/publications/i/item/9789241565684}.

\bibitem[Yao and Wang, 2019]{Yao2019}
Yao, Y. and Wang, H. (2019).
\newblock Optimal subsampling for softmax regression.
\newblock {\em Statistical Papers}, 60:585--599.

\bibitem[Yu et~al., 2022]{Yu2022}
Yu, J., Wang, H., Ai, M., and Zhang, H. (2022).
\newblock Optimal distributed subsampling for maximum quasi-likelihood
  estimators with massive data.
\newblock {\em Journal of the American Statistical Association},
  117(537):265--276.

\bibitem[Zhan et~al., 2022]{Zhan2022}
Zhan, X., Wang, Y., and Chan, A.~B. (2022).
\newblock Asymptotic optimality for active learning processes.
\newblock In {\em Proceedings of the Thirty-Eighth Conference on Uncertainty in
  Artificial Intelligence}.

\bibitem[Zhang et~al., 2021]{Zhang2021}
Zhang, T., Ning, Y., and Ruppert, D. (2021).
\newblock Optimal sampling for generalized linear models under measurement
  constraints.
\newblock {\em Journal of Computational and Graphical Statistics},
  30(1):106--114.

\end{thebibliography}


\appendix

\setcounter{figure}{0}

\section{Proofs}
\label{sec:proofs}

\subsection{Proof of Lemma \ref{lemma:chain_rule}}

According to the chain rule in matrix differential calculus \citep{Petersen2012} we have that
\[
\frac{\partial \Phi(\boldsymbol{\Gamma}(\boldsymbol{\mu}; \boldsymbol{\theta}_0))}{\partial \mu_i} = \tr\left(
\boldsymbol{\phi}(\boldsymbol{\Gamma}(\boldsymbol{\mu}; \boldsymbol{\theta}_0))\T
\frac{\partial\boldsymbol{\Gamma}(\boldsymbol{\mu}; \boldsymbol{\theta}_0)}{\partial\mu_i}
\right),
\]
where $\boldsymbol{\phi}(\mathbf{U}) = \frac{\partial \Phi(\mathbf{U})}{\partial \mathbf{U}}$ is the $p \times p$ matrix derivative of $\Phi$ with respect to its matrix argument, and $\frac{\partial \boldsymbol{\Gamma}(\boldsymbol{\mu}; \boldsymbol{\theta}_0)}{\partial \mu_i}$ the elementwise derivative of $\boldsymbol{\Gamma}(\boldsymbol{\mu}; \boldsymbol{\theta}_0)$ with respect to $\mu_i$. Since $\Phi(\boldsymbol{\Gamma}(\boldsymbol{\mu}; \boldsymbol{\theta}_0))$ is symmetric, $\boldsymbol{\phi}(\boldsymbol{\Gamma}(\boldsymbol{\mu}; \boldsymbol{\theta}_0))$ must also be symmetric, which proves \eqref{eq:chain_rule}.

By the assumptions, $\boldsymbol{\Gamma}(\boldsymbol{\mu}; \boldsymbol{\theta}_0)$ decreases monotonically with $\mu_i$ in the Loewner order sense, which implies that $\frac{\partial\boldsymbol{\Gamma}(\boldsymbol{\mu}; \boldsymbol{\theta}_0)}{\partial\mu_i}$ is negative semi-definite. Therefore, there exists a real matrix $\mathbf{U}$ such that $\mathbf{U}\mathbf{U}\T = -\frac{\partial\boldsymbol{\Gamma}(\boldsymbol{\mu}; \boldsymbol{\theta}_0)}{\partial\mu_i}$. Moreover, $\Phi(\boldsymbol{\Gamma}(\boldsymbol{\mu}; \boldsymbol{\theta}_0))$ is monotone for Loewner's ordering and hence a monotone decreasing function of $\mu_i$, so we must have 
\[
\frac{\partial \Phi(\boldsymbol{\Gamma}(\boldsymbol{\mu}; \boldsymbol{\theta}_0))}{\partial \mu_i} 
\le 0,
\]
which by the above is equivalent to 
\[
\tr\left(
\mathbf{U}\T
\boldsymbol{\phi}(\boldsymbol{\Gamma}(\boldsymbol{\mu}; \boldsymbol{\theta}_0))
\mathbf{U}
\right) \ge 0, \quad \mathbf{U}\mathbf{U}\T = -\frac{\partial\boldsymbol{\Gamma}(\boldsymbol{\mu}; \boldsymbol{\theta}_0)}{\partial\mu_i}.
\]
This inequality holds true for every $\mathbf{U}$, and hence for every possible value of the matrix $\frac{\partial\boldsymbol{\Gamma}(\boldsymbol{\mu}; \boldsymbol{\theta}_0)}{\partial\mu_i}$, if and only if $\boldsymbol{\phi}(\boldsymbol{\Gamma}(\boldsymbol{\mu}; \boldsymbol{\theta}_0))$ is positive semi-definite. Consequently, there exists a real $p \times p$ matrix $\mathbf{L}(\boldsymbol{\mu}; \boldsymbol{\theta}_0)$ such that $\mathbf{L}(\boldsymbol{\mu}; \boldsymbol{\theta}_0)\mathbf{L}(\boldsymbol{\mu}; \boldsymbol{\theta}_0)\T = \boldsymbol{\phi}(\boldsymbol{\Gamma}(\boldsymbol{\mu}; \boldsymbol{\theta}_0))$.

\subsection{Proof of Lemma \ref{lemma:differentiability}}

\paragraph{Proof of \ref{prop:differentiability_a}}

We have by \eqref{eq:approximate_covariance} and  \eqref{eq:Vmat} that 
\[
\boldsymbol{\Gamma}(\boldsymbol{\mu}; \boldsymbol{\theta}_0) = 
\begin{cases}
    \sum_{i \in \mathcal D} \mu_i^{-1}\mathbf{H}(\boldsymbol \theta_0)^{-1}\boldsymbol{\psi}_i(\boldsymbol{\theta}_0)\boldsymbol{\psi}_i(\boldsymbol{\theta}_0)\T \mathbf{H}(\boldsymbol \theta_0)^{-1}, & \text{for PO-WR or MULTI designs, and} \\
    \sum_{i \in \mathcal D} (\mu_i^{-1} - 1)\mathbf{H}(\boldsymbol \theta_0)^{-1}\boldsymbol{\psi}_i(\boldsymbol{\theta}_0)\boldsymbol{\psi}_i(\boldsymbol{\theta}_0)\T \mathbf{H}(\boldsymbol \theta_0)^{-1}, & \text{for PO-WOR}.
\end{cases}
\]
Taking the derivative with respect to $\mu_i$, we obtain
\[
\frac{\partial \boldsymbol{\Gamma}(\boldsymbol{\mu}; \boldsymbol{\theta}_0)}{\partial \mu_i} = -\mu_i^{-2}\mathbf{H}(\boldsymbol{\theta}_0)^{-1}\boldsymbol{\psi}_i(\boldsymbol{\theta}_0)\boldsymbol{\psi}_i(\boldsymbol{\theta}_0)\T\mathbf{H}(\boldsymbol{\theta}_0)^{-1}.
\]

\paragraph{Proof of \ref{prop:differentiability_b}–\ref{prop:differentiability_e}} Follows by the following rules from matrix differential calculus \citep{Petersen2012}: 

\begin{enumerate}[label = \alph*)]
    \setcounter{enumi}{1}
    \item $\frac{\partial \log \det(\mathbf{U})}{\partial U} = \mathbf{U}^{-1}$, provided that $\mathbf{U}$ is of full rank.
    \item $\frac{\partial \lambda_{\max}(\mathbf{U})}{\partial U} = \mathbf{v}\mathbf{v}\T$, where $\mathbf{v}$ is an eigenvector pertaining to the maximal eigenvalue of $\mathbf{U}$, provided that $\mathbf{v}$ is unique. 
    \item $\frac{\partial \tr(\mathbf{U}\boldsymbol{A})}{\partial U} = \boldsymbol{A}\T$.
    \item $\frac{\partial \tr(\mathbf{U}^q)}{\partial \mathbf{U}} = q(\mathbf{U}^{q-1})\T$, so that $\frac{\partial \tr(\mathbf{U}^q)^{1/q}}{\partial U} = 
    \frac{1}{q}\tr(\mathbf{U}^q)^{1/q-1}\frac{\partial \tr(\mathbf{U}^p)}{\partial \mathbf{U}} 
    = \tr(\mathbf{U}^q)^{1/q-1}(\mathbf{U}^{q-1})\T$, provided that $\mathbf{U}$ is of full rank. The final result follows by symmetry of $\boldsymbol{\Gamma}(\boldsymbol{\mu}; \boldsymbol{\theta}_0)$.
\end{enumerate}


\subsection{Proof of Lemma \ref{lemma:partial_derivatives}}

Combining the results of Lemma \ref{lemma:chain_rule} and \ref{lemma:differentiability}, we observe for PO-WR, PO-WOR and MULTI designs that the partial derivative of $\Phi(\boldsymbol{\Gamma}(\boldsymbol{\mu}; \boldsymbol{\theta}_0))$ with respect to $\mu_i$, whenever it exists, is given by
\begin{align*}
    \frac{\partial \Phi(\boldsymbol{\Gamma}(\boldsymbol{\mu}; \boldsymbol{\theta}_0))}{\partial \mu_i} 
    & = 
    -\tr\left(\boldsymbol{\phi}(\boldsymbol{\Gamma}(\boldsymbol{\mu}; \boldsymbol{\theta}_0))\mu_i^{-2}\mathbf{H}(\boldsymbol{\theta}_0)^{-1}\boldsymbol{\psi}_i(\boldsymbol{\theta}_0)\boldsymbol{\psi}_i(\boldsymbol{\theta}_0)\T\mathbf{H}(\boldsymbol{\theta}_0)^{-1}\right) \\
    & = 
    -\mu_i^{-2}
    \tr(\boldsymbol{\psi}_i(\boldsymbol{\theta}_0)\T
    \mathbf{H}(\boldsymbol{\theta}_0)^{-1}
    \mathbf{L}(\boldsymbol{\mu}; \boldsymbol{\theta}_0)\mathbf{L}(\boldsymbol{\mu}; \boldsymbol{\theta}_0)\T
    \mathbf{H}(\boldsymbol{\theta}_0)^{-1}
    \boldsymbol{\psi}_i(\boldsymbol{\theta}_0)) \\
    & = -\mu_i^{-2}
    \bigr\rvert\bigr\rvert
    \mathbf{L}(\boldsymbol{\mu}; \boldsymbol{\theta}_0)\T
    \mathbf{H}(\boldsymbol{\theta}_0)^{-1}
    \boldsymbol{\psi}_i(\boldsymbol{\theta}_0)
    \bigr\rvert\bigr\rvert_2^2.
\end{align*}
The second equality follows from the cyclic property of the trace and definition of $\mathbf{L}(\boldsymbol{\mu}, \boldsymbol{\theta}_0)$, and the third by noting that the expression within the parentheses is a scalar and equals the squared Euclidean norm of the vector $\mathbf{L}(\boldsymbol{\mu}; \boldsymbol{\theta}_0)\T\mathbf{H}(\boldsymbol{\theta}_0)^{-1}\boldsymbol{\psi}_i(\boldsymbol{\theta}_0)$.

\subsection{Proof of Lemma \ref{lemma:d_taylor_expansion}}

By a second order Taylor expansion around $\boldsymbol{\theta}_0$, we have that
\begin{align*}
d(\hat{\boldsymbol{\theta}}_{\!\boldsymbol{\mu}}) 
& = 
d(\boldsymbol{\theta}_0) + 
\nabla d(\boldsymbol{\theta})\T
\bigr\rvert_{\boldsymbol{\theta}=\boldsymbol{\theta}_0}
(\hat{\boldsymbol{\theta}}_{\!\boldsymbol{\mu}} - \boldsymbol{\theta}_0) 
+ 
\frac{1}{2}(\hat{\boldsymbol{\theta}}_{\!\boldsymbol{\mu}} - \boldsymbol{\theta}_0)\T
\mathbf{H}_d(\boldsymbol{\theta}_0)
(\hat{\boldsymbol{\theta}}_{\!\boldsymbol{\mu}} - \boldsymbol{\theta}_0) + o_p(||(\hat{\boldsymbol{\theta}}_{\!\boldsymbol{\mu}} - \boldsymbol{\theta}_0)||_2^2) ,
\end{align*}
where the first two terms, by definition of $d(\boldsymbol{\theta})$, are zero, and $\mathbf{H}_d(\boldsymbol{\theta}_0) = \frac{\partial^2 d(\boldsymbol{\theta})}
{\partial\boldsymbol{\theta}\partial \boldsymbol{\theta}\T}$. By the assumptions on $\hat{\boldsymbol{\theta}}_{\!\boldsymbol{\mu}}$, we have that $\hat{\boldsymbol{\theta}}_{\!\boldsymbol{\mu}} - \boldsymbol{\theta}_0 = o_p(n^{-1/2})$ and $\E[|\hat{\boldsymbol{\theta}}_{\!\boldsymbol{\mu}} - \boldsymbol{\theta}_0|^{2+\delta}]<\infty$ (elementwise) for some $\delta>0$. By bounded convergence, this implies for the remainder that $\E[o_p(||(\hat{\boldsymbol{\theta}}_{\!\boldsymbol{\mu}} - \boldsymbol{\theta}_0)||_2^2)] = o(n^{-1})$. We have further that
\begin{align*}
    \E\left[(\hat{\boldsymbol{\theta}}_{\!\boldsymbol{\mu}} - \boldsymbol{\theta}_0)\T
    \mathbf{H}_d(\boldsymbol{\theta}_0)
    (\hat{\boldsymbol{\theta}}_{\!\boldsymbol{\mu}} - \boldsymbol{\theta}_0)\right]
    & = 
    \tr 
    \left(
    \mathbf{H}_d(\boldsymbol{\theta}_0)
    \COV(\hat{\boldsymbol{\theta}}_{\!\boldsymbol{\mu}} - \boldsymbol{\theta}_0)
    \right) +
    \E[\hat{\boldsymbol{\theta}}_{\!\boldsymbol{\mu}} - \boldsymbol{\theta}_0]\T
    \mathbf{H}_d(\boldsymbol{\theta}_0)
    \E[\hat{\boldsymbol{\theta}}_{\!\boldsymbol{\mu}} - \boldsymbol{\theta}_0] \\
    & = \tr\left(
    \boldsymbol{\Gamma}(\boldsymbol{\mu}; \boldsymbol{\theta}_0) 
    \mathbf{H}_d(\boldsymbol{\theta}_0)
    \right) + o(n^{-1}),
\end{align*}
where the first equality follows from properties of quadratic forms \citep{Mathai1992}, and the second by assumptions \eqref{eq:unbiased}--\eqref{eq:approximate_covariance} on $\hat{\boldsymbol{\theta}}_{\!\boldsymbol{\mu}}$ and the cyclic property of the trace.

\subsection{Proof of Proposition \ref{prop:optimality_conditions}}

First we note that the matrix $\mathbf{L}(\boldsymbol{\mu}^*; \boldsymbol{\theta}_0)$ exists by Lemma \ref{lemma:chain_rule} whenever the objective function is differentiable at $\boldsymbol{\mu}^*$. Hence, the coefficients $c_i$ are positive, the square roots $\sqrt{c_i}$ are real, and the optimality conditions \eqref{eq:mui} and \eqref{eq:powor_a}–\eqref{eq:powor_c} well-defined.

\paragraph{Proof of \ref{prop:optimality_conditions_PO-WR}}

Consider the function $\Phi(\boldsymbol{\Gamma}(\boldsymbol{\mu}; \boldsymbol{\theta}_0))$ subject to the constraints $\sum_{i\in \mathcal D} \mu_i = n$, and $\mu_i > 0$ for all $i \in \mathcal D$. By the Lagrange multiplier method \citep{Boyd2004}, the constrained stationary points of $\Phi(\boldsymbol{\Gamma}(\boldsymbol{\mu}; \boldsymbol{\theta}_0))$ are obtained as the stationary points of the Lagrangian
\[
\Lambda(\boldsymbol \mu, \alpha) = \Phi(\boldsymbol{\Gamma}(\boldsymbol{\mu}; \boldsymbol{\theta}_0)) +\alpha g(\boldsymbol \mu) , \quad  g(\boldsymbol \mu) = \sum_{i\in \mathcal D} \mu_i - n.
\]
Taking the derivatives with respect to $\boldsymbol{\mu}$ and $\alpha$, we obtain the system of equations
\begin{align*}
    \nabla \Lambda(\boldsymbol{\mu}, \alpha) = \mathbf{0} \quad
    \Leftrightarrow \quad
    \begin{cases} 
    g(\boldsymbol{\mu}) = 0 \\ 
    -\nabla_{\!\boldsymbol{\mu}} \Phi(\boldsymbol{\Gamma}(\boldsymbol{\mu}; \boldsymbol{\theta}_0)) = \alpha \nabla g(\boldsymbol{\mu}) .
    \end{cases}
\end{align*}
Now, $\frac{\partial\Phi(\boldsymbol{\Gamma}(\boldsymbol{\mu}; \boldsymbol{\theta}_0))}{\partial \mu_i} = -c_i/\mu_i^2$ by Lemma \ref{lemma:partial_derivatives} and definition of $c_i$, and $\frac{\partial g(\boldsymbol{\mu})}{\partial \mu_i} = 1$. A stationary point therefore satisfies the system of equations $\alpha = c_1/\mu^{2}_1 = \ldots = c_N/\mu^{2}_N$ for all $i \in \mathcal D$. For $\boldsymbol{\mu}^*$ to be $\Phi$-optimal we must have $\mu_i^* \propto \sqrt{c_i}$, $\mu_i^* > 0$, and $\sum_{i \in \mathcal D} \mu_i = n$, and hence
\[
    \mu_i^* = n\frac{\sqrt{c_i}}{\sum_{j \in \mathcal D} \sqrt{c_j}} \text{ for all } i \in \mathcal D. 
\]

\paragraph{Proof of \ref{prop:optimality_conditions_PO-WOR}}

Consider the function $\Phi(\boldsymbol{\Gamma}(\boldsymbol{\mu}; \boldsymbol{\theta}_0))$ subject to the constraints $\sum_{i\in \mathcal D} \mu_i = n$ and $0 < \mu_i \le 1$ for all $i \in \mathcal D$. Also consider the Lagrangian
\[
\Lambda(\boldsymbol{\mu}, \alpha, \boldsymbol{\beta}) = \Phi(\boldsymbol{\Gamma}(\boldsymbol{\mu}; \boldsymbol{\theta}_0)) + 
\alpha g(\boldsymbol{\mu}) +
\sum_{i \in \mathcal D} \beta_i h_{i}(\boldsymbol{\mu}),
\]
where $g(\boldsymbol{\mu}) = \sum_{i \in \mathcal D} \mu_i - n$ and $h_i(\boldsymbol{\mu}) = \mu_i - 1$. The constrained stationary points of $\Phi(\boldsymbol{\Gamma}(\boldsymbol{\mu}; \boldsymbol{\theta}_0))$ are characterised as the solutions to the Karush-Kuhn-Tucker conditions \citep{Boyd2004}:
\begin{itemize}
    \item \textit{Stationarity}: $-\nabla_{\!\boldsymbol{\mu}} \Phi(\boldsymbol{\Gamma}(\boldsymbol{\mu}; \boldsymbol{\theta}_0)) = \alpha \nabla g(\boldsymbol{\mu}) + \sum_{i \in \mathcal D} \beta_i \nabla h_{i}(\boldsymbol{\mu})$.
    \item \textit{Primal feasibility}: $g(\boldsymbol{\mu}) = 0$, and $h_i(\boldsymbol{\mu}) \le 0$ for all $i \in \mathcal D$.
    \item \textit{Dual feasibility}: $\beta_i \ge 0$ for all $i \in \mathcal D$.
    \item \textit{Complementary slackness}: $\beta_i h_i(\boldsymbol{\mu}) = 0$ for all $i \in \mathcal D$.
\end{itemize}
First note that $\frac{\partial \Phi(\boldsymbol{\Gamma}(\boldsymbol{\mu}; \boldsymbol{\theta}_0))}{\partial \mu_i} = -c_i/\mu_i^2$  by Lemma \ref{lemma:partial_derivatives} and definition of $c_i$, $\frac{\partial g(\boldsymbol{\mu})}{\partial \mu_i} = 1$, and $\frac{\partial h_i(\boldsymbol{\mu})}{\partial \mu_j} = 1$ if $i = j$ and $0$ otherwise. Consider a sampling scheme $\boldsymbol{\mu} \in \mathcal M_n$ and let $\mathcal E = \{i \in \mathcal D: \mu_i = 1$\} and $n_{\mathcal E} = |\mathcal E|$. For the Karush-Kuhn-Tucker conditions to be satisfied, we must have that 
\begin{enumerate}[label = \roman*)]
    \item $\mu_i \le 1$ and $\sum_{i \in \mathcal D\setminus \mathcal E} \mu_i = n - n_{\mathcal E}$, by the primal feasibility condition,
    \label{list:primal_feasibility}
    \item $\beta_i = 0$ if $\mu_i < 1$, by the complementary slackness condition,
    \item $c_i / \mu_i^2 = \alpha + \beta_i$ by the stationarity condition, which by the above implies that 
    \begin{equation*} 
    c_i = 
    \begin{cases}
        \alpha + \beta_i & \text{if } \mu_i = 1 \\
        \alpha\mu_i^2 & \text{if } \mu_i < 1
    \end{cases}
    \quad \Leftrightarrow \quad 
    \begin{cases}
        \beta_i = c_i - \alpha & \text{if } \mu_i = 1 \\
        \alpha = c_i/\mu_i^2 & \text{if } \mu_i < 1,
    \end{cases}  
    \end{equation*}
    \label{list:stationarity}
    \item  $\beta_i \ge 0$ by the dual feasibility condition, which by the above implies that $\beta_i = c_i - \alpha = c_i - c_j/\mu_j^2 \ge 0$ for $i \in \mathcal E$ and $j \in \mathcal D \setminus \mathcal E$.
    \label{list:dual_feasibility}
\end{enumerate}
The condition \eqref{eq:powor_a} follows from \ref{list:primal_feasibility}, \eqref{eq:powor_b} from \ref{list:primal_feasibility} and \ref{list:stationarity}, and \eqref{eq:powor_c} from \ref{list:dual_feasibility}.



\subsection{Proof of Proposition \ref{prop:convex}}

Note first that the domain $\mathcal M_n$ of $\boldsymbol{\mu}$ is convex. The results hence follow from the second derivative test by showing that the Hessian matrix of $\Phi(\boldsymbol{\Gamma}(\boldsymbol{\mu}; \boldsymbol{\theta}_0))$ is positive semi-definite on $\mathcal M_n$.

\paragraph{Proof of \ref{prop:convex_a}}
We have by Lemma \ref{lemma:partial_derivatives} that
\[
    \frac{\partial \Phi(\boldsymbol{\Gamma}(\boldsymbol{\mu}; \boldsymbol{\theta}_0))}{\partial \mu_i} = -\mu_i^{-2}
    \bigr\rvert\bigr\rvert\mathbf{L}(\boldsymbol{\mu}; \boldsymbol{\theta}_0)\T
    \mathbf{H}(\boldsymbol{\theta}_0)^{-1}
    \boldsymbol{\psi}_i(\boldsymbol{\theta}_0)\bigr\rvert\bigr\rvert^2_{2}.
\]
For the L-optimality criterion the matrix $\mathbf{L}(\boldsymbol{\mu}; \boldsymbol{\theta}_0) = \mathbf{L}$ does not depend on $\boldsymbol{\mu}$. The second-order partial derivatives are given by
\[
    \frac{\partial^2 \Phi(\boldsymbol{\Gamma}(\boldsymbol{\mu}; \boldsymbol{\theta}_0))}{\partial \mu_i\partial \mu_j} = 
    \begin{cases}
    2\mu_i^{-3}\bigr\rvert\bigr\rvert\mathbf{L}\T
    \mathbf{H}(\boldsymbol{\theta}_0)^{-1}
    \boldsymbol{\psi}_i(\boldsymbol{\theta}_0)\bigr\rvert\bigr\rvert^2_{2} \ge 0 & \text{if }i = j, \\
    0 & \text{if }i \ne j.
    \end{cases}
\]
This matrix is diagonal with non-negative entries for all $\mu_i>0$, and hence positive semi-definite on $\mathcal M_n$.

\paragraph{Proof of \ref{prop:convex_b}}
We show that $\det\left(\boldsymbol{\Gamma}(\boldsymbol{\mu};\boldsymbol{\theta}_0)\right)$ is log-convex in $\boldsymbol{\mu}$, i.e., that $\log\det\left(\boldsymbol{\Gamma}(\boldsymbol{\mu};\boldsymbol{\theta}_0)\right)$ is convex.

First note that $\log\det\left(\boldsymbol{\Gamma}(\boldsymbol{\mu};\boldsymbol{\theta}_0)\right) = \log\det\left(\mathbf{V}(\boldsymbol{\mu}; \boldsymbol{\theta}_0))\right) -2\log\det\left(\mathbf{H}(\boldsymbol{\theta}_0)\right)$, where $\mathbf{V}(\boldsymbol{\mu}; \boldsymbol{\theta}_0)$ is given by \eqref{eq:Vmat} and $\mathbf{H}(\boldsymbol{\theta}_0)$ does not depend on $\boldsymbol{\mu}$. Thus, it suffices to show that $\log\det\left(\mathbf{V}(\boldsymbol{\mu}; \boldsymbol{\theta}_0)\right)$ is convex in $\boldsymbol{\mu}$. We obtain the desired result by showing that the Hessian of $\Phi(\mathbf{V}(\boldsymbol{\mu}; \boldsymbol{\theta}_0))$ can be decomposed as the Hadamard product between two positive semi-definite matrices, and hence is positive semi-definite \citep{Horn1990}.

Consider first a PO-WR or multinomial sampling design. The partial derivatives of $\Phi(\mathbf{V}(\boldsymbol{\mu}; \boldsymbol{\theta}_0))$ are given by
\begin{align*}
    \frac{\partial \Phi(\mathbf{V}(\boldsymbol{\mu}; \boldsymbol{\theta}_0))}{\partial \mu_i} 
    & = 
    -\mu_i^{-2}
    \boldsymbol{\psi}_i(\boldsymbol{\theta}_0)\T
    \mathbf{V}(\boldsymbol{\mu}; \boldsymbol{\theta}_0)^{-1}
    \boldsymbol{\psi}_i(\boldsymbol{\theta}_0), \\
    \frac{\partial^2 \Phi(\mathbf{V}(\boldsymbol{\mu}; \boldsymbol{\theta}_0))}{\partial \mu_i^2} 
    & = 2\mu_i^{-3}
    \boldsymbol{\psi}_i(\boldsymbol{\theta}_0)\T
    \mathbf{V}(\boldsymbol{\mu}; \boldsymbol{\theta}_0)^{-1}
    \boldsymbol{\psi}_i(\boldsymbol{\theta}_0) - 
    \mu_i^{-4}
    (\boldsymbol{\psi}_i(\boldsymbol{\theta}_0)\T
    \mathbf{V}(\boldsymbol{\mu}; \boldsymbol{\theta}_0)^{-1}
    \boldsymbol{\psi}_i(\boldsymbol{\theta}_0))^2, \\
    \frac{\partial^2 \Phi(\mathbf{V}(\boldsymbol{\mu}; \boldsymbol{\theta}_0))}{\partial \mu_i\partial \mu_j}
    & = - (\mu_i\mu_j)^{-2}
    (\boldsymbol{\psi}_i(\boldsymbol{\theta}_0)\T
    \mathbf{V}(\boldsymbol{\mu}; \boldsymbol{\theta}_0)^{-1}
    \boldsymbol{\psi}_j(\boldsymbol{\theta}_0))^2, \quad i \ne j.
\end{align*}
These results follow in analogy with the proof of Lemma \ref{lemma:partial_derivatives} by the chain rule \eqref{eq:chain_rule} and the following rules for matrix differentiation \citep{Petersen2012}: 
\[
\frac{\partial \mathbf{a}\T\mathbf{X}\mathbf{a}}{\partial \mathbf{X}} = \mathbf{a}\mathbf{a}\T, \quad 
\frac{\partial \log \det(\mathbf{Y})}{\partial x} = \tr\left(\mathbf{Y}^{-1}\frac{\partial \mathbf{Y}}{\partial x} \right), \quad \text{and} \quad 
\frac{\partial \mathbf{Y}^{-1}}{\partial x} = - \mathbf{Y}^{-1} \frac{\partial \mathbf{Y}}{\partial x}\mathbf{Y}^{-1}.
\]
Let $\boldsymbol{u}_i = \boldsymbol{\psi}_i(\boldsymbol{\theta}_0)/\sqrt{\mu_i}$ and $\mathbf{U}$ be the matrix with rows $\boldsymbol{u}_i\T$. Also, let $\mathbf{A} = \mathbf{U}(\mathbf{U}\T\mathbf{U})^{-1}\mathbf{U}\T$ and $a_{ij}$ the elements of $\mathbf{A}$. We note the following:
\begin{itemize}
    \item $\mathbf{U}\T\mathbf{U} = \sum_{i \in \mathcal D} 
    \mu_i^{-1}\boldsymbol{\psi}_i(\boldsymbol{\theta}_0)\boldsymbol{\psi}_i(\boldsymbol{\theta}_0)\T = 
    \mathbf{V}(\boldsymbol{\mu}; \boldsymbol{\theta}_0)$,
    \item $a_{ij} = 
    \boldsymbol{u}_i\T (\mathbf{U}\T\mathbf{U})^{-1} \boldsymbol{u}_j = 
    (\mu_i\mu_j)^{-1/2}\boldsymbol{\psi}_i(\boldsymbol{\theta}_0)\T \mathbf{V}(\boldsymbol{\mu}; \boldsymbol{\theta}_0)^{-1}
    \boldsymbol{\psi}_j(\boldsymbol{\theta}_0)$,
    \item $\mathbf{A}$ is an idempotent matrix, i.e., $\mathbf{A}^2 = \mathbf{A}$, which implies that $a_{ii} = \sum_{j} a_{ij}^2$, 
    \item $a_{ii} = \boldsymbol{u}_i\T (\mathbf{U}\T\mathbf{U})^{-1} \boldsymbol{u}_i = \boldsymbol{u}_i\T \mathbf{V}(\boldsymbol{\mu}; \boldsymbol{\theta}_0)^{-1}\boldsymbol{u}_i > 0$, since $\mathbf{V}(\boldsymbol{\mu}; \boldsymbol{\theta}_0)$ by assumption is positive definite.
\end{itemize}
We may now write
\[
\frac{\partial^2 \Phi(\mathbf{V}(\boldsymbol{\mu}; \boldsymbol{\theta}_0))}{\partial \mu_i\partial \mu_j} = 
\begin{cases}
    \mu_i^{-2}(2 a_{ii} - a_{ii}^2) & \text{if }i = j, \\
    (\mu_i\mu_j)^{-1} a_{ij}^2 & \text{if } i \ne j.
\end{cases}
\]
We recognise the Hessian matrix $\frac{\partial^2 \Phi(\mathbf{V}(\boldsymbol{\mu}; \boldsymbol{\theta}_0))}{\partial \boldsymbol{\mu}\partial \boldsymbol{\mu}\T}$ as the Hadamard product $\mathbf{M} \otimes \mathbf{B}$ of a rank-one matrix 
$\mathbf{M} = \mathbf{m}\mathbf{m}\T$ with $\mathbf{m} = (\mu_1^{-1}, \ldots, \mu_N^{-1})\T$, and a symmetric matrix $\mathbf{B}$ with entries
\[
b_{ij} = 
\begin{cases}
    2 a_{ii} - a_{ii}^2 & \text{if }i = j, \\
    a_{ij}^2 & \text{if } i \ne j.
\end{cases}
\]
The matrix $\mathbf{M}$ has eigenvalues $\mathbf{m}\T\mathbf{m}$ and $0$, and hence is positive semi-definite. The matrix $\mathbf{B}$ is diagonally dominant with positive entries, since $a_{ii} = \sum_{j}a_{ij}^2, b_{ii} = 2 a_{ii} - a_{ii}^2 = a_{ii} + \sum_{j\ne i}a_{ij}^2$, and $a_{ii}>0$ implies
\[
b_{ii} > a_{ii} > 0, \quad \text{and }
b_{ii} > \sum_{j\ne i}a_{ij}^2 = \sum_{j\ne i} b_{ij}.
\]
Hence, $\mathbf{B}$ is positive definite \citep{Horn1990}. It follows that the Hessian matrix $\frac{\partial^2 \Phi(\mathbf{V}(\boldsymbol{\mu}; \boldsymbol{\theta}_0))}{\partial \boldsymbol{\mu}\partial \boldsymbol{\mu}\T}$ is positive semi-definite on $\mathcal M_n$ for PO-WR and multinomial sampling designs.

It remains to prove convexity for PO-WOR. First note that the function $\log\det\left(\mathbf{V}(\boldsymbol{\mu}; \boldsymbol{\theta}_0)\right)$, by assumptions on $\mathbf{V}(\boldsymbol{\mu}; \boldsymbol{\theta}_0)$, is differentiable and continuous on $\mathcal M_n$. It suffices, by continuity, to prove that the Hessian is positive semi-definite on the interior of $\mathcal M_n$. Consider therefore a point $\boldsymbol{\mu} \in \mathcal M_n$ such that $\mu_i < 1$ for all $i$. Let $\boldsymbol{u}_i = \boldsymbol{\psi}_i(\boldsymbol{\theta}_0)\T\sqrt{1 - \mu_i}/\sqrt{\mu_i}$ and $\mathbf{U}$ be the matrix with rows $\boldsymbol{u}_i\T$. Also let $\mathbf{A} = \mathbf{U}(\mathbf{U}\T\mathbf{U})^{-1}\mathbf{U}\T$ and $a_{ij}$ the elements of $\mathbf{A}$. Similar to above, we may now write
\[
\frac{\partial^2 \Phi(\mathbf{V}(\boldsymbol{\mu}; \boldsymbol{\theta}_0))}{\partial \mu_i\partial \mu_j} = 
\begin{cases}
    \mu_i^{-2}(1-\mu_i)^{-1}(2 a_{ii} - a_{ii}^2) & \text{if }i = j, \\
    (\mu_i\mu_j)^{-1}(1-\mu_i)^{-1/2}(1-\mu_j)^{-1/2}a_{ij}^2 & \text{if } i \ne j.
\end{cases}
\]
We recognise the Hessian matrix $\frac{\partial^2 \Phi(\mathbf{V}(\boldsymbol{\mu}; \boldsymbol{\theta}_0))}{\partial \boldsymbol{\mu}\partial \boldsymbol{\mu}\T}$ as the Hadamard product $\mathbf{M} \otimes \mathbf{B}$ of a rank-one matrix 
$\mathbf{M} = \mathbf{m}\mathbf{m}\T$ with $\mathbf{m} = (\mu_1^{-1}(1-\mu_1)^{-1/2}, \ldots, \mu_N^{-1}(1-\mu_N)^{-1/2})\T$, and a symmetric matrix $\mathbf{B}$ with entries
\[
b_{ij} = 
\begin{cases}
    2 a_{ii} - a_{ii}^2 & \text{if }i = j, \\
    a_{ij}^2 & \text{if } i \ne j.
\end{cases}
\]
The remainder of the proof follows in complete analogy with the proof for PO-WR and multinomial sampling designs.

\subsection{Proof of Proposition \ref{prop:equivalence}}

\paragraph{Proof of \ref{prop:equivalence_a}}

First note that the Hessian $\mathbf{H}_d(\boldsymbol{\theta})$ is positive semi-definite at $\boldsymbol{\theta} = \boldsymbol{\theta}_0$, since $\boldsymbol{\theta}_0$ is the global minimiser of $d(\boldsymbol{\theta})$. Hence, there exists a matrix $\mathbf{L}$ such that $\mathbf{L}\mathbf{L}\T = \mathbf{H}_d(\boldsymbol{\theta})$. $\mathbf{L}$ is non-zero since $\mathbf{H}_d(\boldsymbol{\theta})$, by assumption, is non-zero. The \textit{d}-optimal sampling scheme $\boldsymbol{\mu}^*$ is defined as the minimiser of the function $
\tr(\boldsymbol{\Gamma}(\boldsymbol{\mu}; \boldsymbol{\theta}_0)\mathbf{H}_d(\boldsymbol{\theta}_0)) = 
\tr(\boldsymbol{\Gamma}(\boldsymbol{\mu}; \boldsymbol{\theta}_0)\mathbf{L}\mathbf{L}\T),
$
which by definition is equivalent to L-optimality with respect to a matrix $\mathbf{L}$ such that $\mathbf{L}\mathbf{L}\T = \mathbf{H}_d(\boldsymbol{\theta}).$

\paragraph{Proof of \ref{prop:equivalence_b}}

Assume that $\boldsymbol{\mu}^*$ is the minimser of $\Phi(\boldsymbol{\Gamma}(\boldsymbol{\mu}; \boldsymbol{\theta}_0))$ and let $d(\boldsymbol{\theta}) = \frac{1}{2}||\mathbf{L}(\boldsymbol{\mu}^*; \boldsymbol{\theta}_0)\T(\boldsymbol{\theta} - \boldsymbol{\theta}_0)||_2^2$ with Hessian matrix $\mathbf{H}_d(\boldsymbol{\theta}_0) = \mathbf{L}(\boldsymbol{\mu}^*; \boldsymbol{\theta}_0)\mathbf{L}(\boldsymbol{\mu}^*; \boldsymbol{\theta}_0)\T$. According to Proposition \ref{prop:optimality_conditions}, the $\Phi$-optimal sampling scheme $\boldsymbol{\mu}^*$ must satisfy the optimality conditions \eqref{eq:mui} or \eqref{eq:powor_a}--\eqref{eq:powor_c} with 
\begin{align*}
c_i & = ||\mathbf{L}(\boldsymbol{\mu}^*; \boldsymbol{\theta}_0)\T\mathbf{H}(\boldsymbol{\theta}_0)^{-1}\boldsymbol{\psi}_i(\boldsymbol{\theta}_0)||_2^2 \\
& =  \boldsymbol{\psi}_i(\boldsymbol{\theta}_0)\T
    \mathbf{H}(\boldsymbol{\theta}_0)^{-1}
    \mathbf{L}(\boldsymbol{\mu}^*; \boldsymbol{\theta}_0)
    \mathbf{L}(\boldsymbol{\mu}^*; \boldsymbol{\theta}_0)\T
    \mathbf{H}(\boldsymbol{\theta}_0)^{-1}
    \boldsymbol{\psi}_i(\boldsymbol{\theta}_0) \\
& = \boldsymbol{\psi}_i(\boldsymbol{\theta}_0)\T
    \mathbf{H}(\boldsymbol{\theta}_0)^{-1}
    \mathbf{H}_d(\boldsymbol{\theta}_0)
    \mathbf{H}(\boldsymbol{\theta}_0)^{-1}
    \boldsymbol{\psi}_i(\boldsymbol{\theta}_0).
\end{align*}
This is identical to the optimality conditions for the \textit{d}-optimality criterion. Moreover, the \textit{d}-optimality criterion is convex in $\boldsymbol{\mu}$ by Proposition \ref{prop:convex}\ref{prop:convex_a} and \ref{prop:equivalence}\ref{prop:equivalence_a}, so $\boldsymbol{\mu}^*$ must be the global minimiser for the \textit{d}-optimality criterion. Now, minimising $\tr(\boldsymbol{\Gamma}(\boldsymbol{\mu}; \boldsymbol{\theta}_0)\mathbf{H}_d(\boldsymbol{\theta}_0))$ is equivalent to minimising $\tr(k\boldsymbol{\Gamma}(\boldsymbol{\mu}; \boldsymbol{\theta}_0)\mathbf{H}_d(\boldsymbol{\theta}_0))$ for any constant $k>0$, so $\boldsymbol{\mu}^*$ is also \textit{d}-optimal with respect to the distance function $d(\boldsymbol{\theta}) = ||\mathbf{L}(\boldsymbol{\mu}^*; \boldsymbol{\theta}_0)\T(\boldsymbol{\theta} - \boldsymbol{\theta}_0)||_2^2$.

\subsection{Proof of Proposition \ref{prop:MD_ER_KL}}

The results follow from Proposition \ref{prop:equivalence}\ref{prop:equivalence_a} since the Hessian matrices of $d_{\mathrm{ER}}(\boldsymbol{\theta})$, $d_{\boldsymbol{\Sigma}}(\boldsymbol{\theta})$, and $d_{\mathrm{KL}}$ are given by 
\begin{align*}
    \frac{\partial^2 d_{\mathrm{ER}}(\boldsymbol{\theta})}{\partial\boldsymbol{\theta}\partial\boldsymbol{\theta}\T} 
    & = \frac{\partial^2 \ell_0(\boldsymbol{\theta})}{\partial\boldsymbol{\theta}\partial\boldsymbol{\theta}\T} 
    = \mathbf{H}(\boldsymbol{\theta}), \\
    \frac{\partial^2 d_{\boldsymbol{\Sigma}}(\boldsymbol{\theta})}{\partial\boldsymbol{\theta}\partial\boldsymbol{\theta}\T} 
    &  = 
    \frac{\partial^2}{\partial\boldsymbol{\theta}\partial\boldsymbol{\theta}\T}
    (\boldsymbol{\theta} - \boldsymbol{\theta}_0)
    \T\boldsymbol{\Sigma}^{-1}(\boldsymbol{\theta} - \boldsymbol{\theta}_0)
    = \boldsymbol{\Sigma}^{-1}, \text{ and} \\
    \frac{\partial^2 d_{\mathrm{KL}}(\boldsymbol{\theta})}{\partial\boldsymbol{\theta}\partial\boldsymbol{\theta}\T} 
    & = \frac{\partial^2}{\partial\boldsymbol{\theta}\partial\boldsymbol{\theta}\T} 
    \sum_{i \in \mathcal D} \int_{\mathcal Y} 
    \log \frac{f_{\boldsymbol{\theta}_0}\!(\boldsymbol{y}|\boldsymbol{x}_i)}{f_{\boldsymbol{\theta}}(\boldsymbol{y}|\boldsymbol{x}_i)}dF_{\boldsymbol{\theta}_0}\!(\boldsymbol{y}|\boldsymbol{x}_i) \\
    & = -\frac{\partial^2}{\partial\boldsymbol{\theta}\partial\boldsymbol{\theta}\T} 
    \sum_{i \in \mathcal D} \int_{\mathcal Y}  \log f_{\boldsymbol{\theta}}(\boldsymbol{y}|\boldsymbol{x}_i)dF_{\boldsymbol{\theta}_0}\!(\boldsymbol{y}|\boldsymbol{x}_i) \\
    & = -\sum_{i \in \mathcal D} \int_{\mathcal Y}  
    \frac{\partial^2\log f_{\boldsymbol{\theta}}(\boldsymbol{y}|\boldsymbol{x}_i)}{\partial\boldsymbol{\theta}\partial\boldsymbol{\theta}\T}  dF_{\boldsymbol{\theta}_0}\!(\boldsymbol{y}|\boldsymbol{x}_i) \\
    & = \E_{\boldsymbol{y}\sim f_{\boldsymbol{\theta}_0}\!(\boldsymbol{y}|\boldsymbol{x})}\!\left[-\sum_{i \in \mathcal D} \frac{\partial^2}{\partial\boldsymbol{\theta}\partial\boldsymbol{\theta}\T}\log f_{\boldsymbol{\theta}}(\boldsymbol{y}|\boldsymbol{x}_i)\right] \\
    & = \E_{\boldsymbol{y}\sim f_{\boldsymbol{\theta}_0}\!(\boldsymbol{y}|\boldsymbol{x})}[\mathbf{H}(\boldsymbol{\theta})] = \widetilde{\mathbf{H}}(\boldsymbol{\theta}). 
\end{align*}
For the Hessian of the Kullback-Leibler distance we have used the Leibniz integral rule to change the order of integration and differentiation \citep[cf.][]{Kullback1951}.

\subsection{Proof of Proposition \ref{prop:invariance}}

Consider a one-to-one differentiable mapping $\boldsymbol{g}: \boldsymbol{\theta} \mapsto \boldsymbol{\eta}$. Denote by
$
\ell^*_{0}(\boldsymbol{\eta}) = \sum_{i\in \mathcal D} \ell_i(\boldsymbol{g}^{-1}(\boldsymbol{\eta}))
$
the induced empirical risk, with the minimiser $\boldsymbol{\eta}_0 = \boldsymbol{g}(\boldsymbol{\theta}_0)$. By the chain rule, the Hessian matrix of $\ell^*_{0}(\boldsymbol{\eta})$ at $\boldsymbol{\eta} = \boldsymbol{\eta}_0$ is given by
\begin{align*}
    \mathbf{H}_{\ell_0^*}(\boldsymbol{\eta}_0) = \frac{\partial^2  \ell^*_{0}(\boldsymbol{\eta})}{\partial \boldsymbol{\eta}\partial\boldsymbol{\eta}\T}\biggr\rvert_{\boldsymbol{\eta} = \boldsymbol{\eta}_0}
    & = \mathbf{J}_{\boldsymbol{g}^{-1}}(\boldsymbol{g}^{-1}(\boldsymbol{\eta}_0))\T\mathbf{H}_{\ell_0}(\boldsymbol{g}^{-1}(\boldsymbol{\eta}_0))\mathbf{J}_{\boldsymbol{g}^{-1}}(\boldsymbol{g}^{-1}(\boldsymbol{\eta}_0)) \\
    & = \mathbf{J}_{\boldsymbol{g}}(\boldsymbol{\theta}_0)\iT\mathbf{H}(\boldsymbol{\theta}_0)\mathbf{J}_{\boldsymbol{g}}(\boldsymbol{\theta}_0)^{-1}. 
\end{align*}
Here we have also used the fact that $\nabla_{\boldsymbol{\theta}}\ell_0(\boldsymbol{\theta})\bigr\rvert_{\boldsymbol{\theta} = \boldsymbol{\theta}_0} = \boldsymbol{0}$, by definition of $\boldsymbol{\theta}_0$ as the minimiser of $\ell_0(\boldsymbol{\theta})$.

Now assume that $\boldsymbol{\mu}^*$ and $\tilde{\boldsymbol{\mu}}^*$ are $d_{\mathrm{ER}}$-optimal for $\boldsymbol{\theta}_0$ and $\boldsymbol{\eta}_0$, respectively. By the latter we mean that $\tilde{\boldsymbol{\mu}}^*$ minimises the expected distance of $\hat{\boldsymbol{\eta}}_{\!\boldsymbol{\mu}} = \boldsymbol{g}(\hat{\boldsymbol{\theta}}_{\!\boldsymbol{\mu}})$ from $\boldsymbol{\eta}_0 = \boldsymbol{g}(\boldsymbol{\theta}_0)$ with respect to the induced empirical risk distance $d_{\mathrm{ER}}^*(\boldsymbol{\eta}) = \ell_0^*(\boldsymbol{\eta}) - \ell_0^*(\boldsymbol{\eta}_0)$. By Proposition \ref{prop:MD_ER_KL}, $\boldsymbol{\mu}^*$ is L-optimal with respect to a matrix $\mathbf{L}$ such that $\mathbf{L}\mathbf{L}\T = \mathbf{H}(\boldsymbol{\theta}_0)$. Similarly, $\tilde{\boldsymbol{\mu}}^*$ is L-optimal with respect to a matrix $\widetilde{\mathbf{L}}$ such that 
\begin{equation}
    \label{eq:ltilde}
\widetilde{\mathbf{L}}\widetilde{\mathbf{L}}\T = \mathbf{H}_{\ell_0^*}(\boldsymbol{\eta}_0) = \mathbf{J}_{\boldsymbol{g}}(\boldsymbol{\theta}_0)\iT\mathbf{H}(\boldsymbol{\theta}_0)\mathbf{J}_{\boldsymbol{g}}(\boldsymbol{\theta}_0)^{-1} = \mathbf{J}_{\boldsymbol{g}}(\boldsymbol{\theta}_0)\iT
\mathbf{L}
\mathbf{L}\T
\mathbf{J}_{\boldsymbol{g}}(\boldsymbol{\theta}_0)^{-1}.    
\end{equation}
Now, $\boldsymbol{\mu}^*$ is the minimiser of the function
\begin{align*}
    \tr(\boldsymbol{\Gamma}(\boldsymbol{\mu}; \boldsymbol{\theta}_0)\mathbf{L}\mathbf{L}\T)
    & = 
    \tr(
    \mathbf{J}_{\boldsymbol{g}}(\boldsymbol{\theta}_0)^{-1}
    \mathbf{J}_{\boldsymbol{g}}(\boldsymbol{\theta}_0)
    \boldsymbol{\Gamma}(\boldsymbol{\mu}; \boldsymbol{\theta}_0)
    \mathbf{J}_{\boldsymbol{g}}(\boldsymbol{\theta}_0)\T
    \mathbf{J}_{\boldsymbol{g}}(\boldsymbol{\theta}_0)\iT
    \mathbf{L}\mathbf{L}\T)
    \\
    & = 
    \tr(
    \boldsymbol{\Gamma}_{\!\boldsymbol{g}}(\boldsymbol{\mu}; \boldsymbol{\theta}_0)
    \mathbf{J}_{\boldsymbol{g}}(\boldsymbol{\theta}_0)\iT
    \mathbf{L}
    \mathbf{L}\T
    \mathbf{J}_{\boldsymbol{g}}(\boldsymbol{\theta}_0)^{-1}
    )    
    \\
    & = 
    \tr(
    \boldsymbol{\Gamma}_{\!\boldsymbol{g}}(\boldsymbol{\mu}; \boldsymbol{\theta}_0)
    \widetilde{\mathbf{L}}
    \widetilde{\mathbf{L}}\T
    ).
\end{align*}
The first equality follows by inserting the identity matrix $\mathbf{I}_{p \times p} = \mathbf{J}_{\boldsymbol{g}}(\boldsymbol{\theta}_0)^{-1}
\mathbf{J}_{\boldsymbol{g}}(\boldsymbol{\theta}_0) = \mathbf{J}_{\boldsymbol{g}}(\boldsymbol{\theta}_0)\T \mathbf{J}_{\boldsymbol{g}}(\boldsymbol{\theta}_0)\iT$ twice, the second equality by \eqref{eq:cov_eta} and the cyclic property of the trace, and the third equality by \eqref{eq:ltilde}. But $\tilde{\boldsymbol{\mu}}^*$ is also a minimiser of $\tr(\boldsymbol{\Gamma}_{\!\boldsymbol{g}}(\boldsymbol{\mu}; \boldsymbol{\theta}_0)
\widetilde{\mathbf{L}}\widetilde{\mathbf{L}}\T)$. Since the L-optimality criterion is convex in $\boldsymbol{\mu}$, the optimum is unique and we must have $\boldsymbol{\mu}^* = \tilde{\boldsymbol{\mu}}^*$. Hence, the $d_{\mathrm{ER}}$-optimality criterion is invariant under the re-parameterisation $\boldsymbol{g}: \boldsymbol{\theta} \mapsto \boldsymbol{\eta}$.

The results for $d_{S}$- and $d_{\mathrm{KL}}$-optimality follow analogously.

\clearpage

\section{Supplementary Figures}
\label{sec:supplementary_figures}
\renewcommand{\thefigure}{S\arabic{figure}}

\FloatBarrier

\begin{figure}[!htb]
    \centering
    \includegraphics[scale = 1]{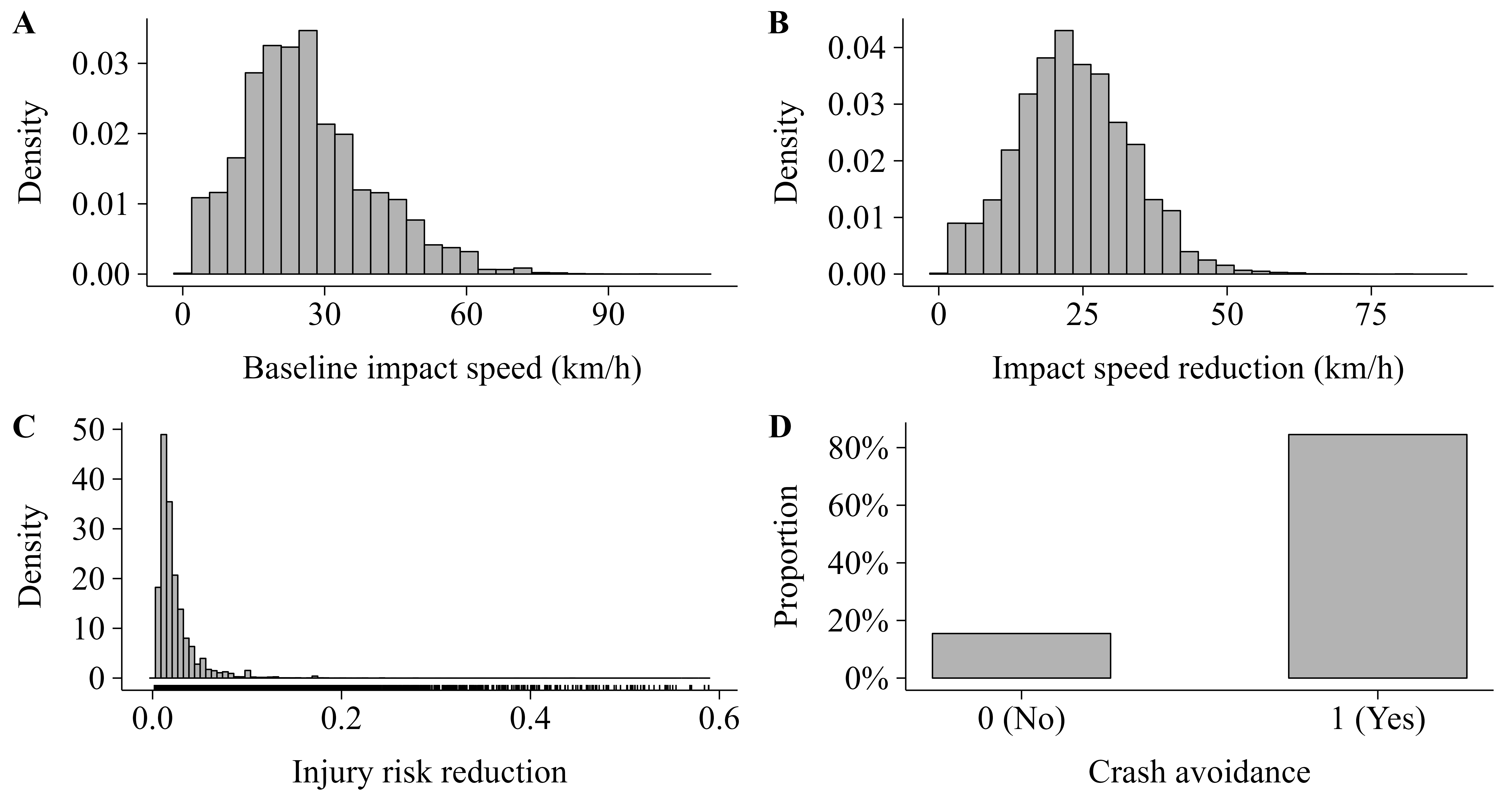}
    \caption{Characteristics of the vehicle safety assessment dataset considered in Section \ref{sec:application}. \textbf{A}: Impact speed distribution under a baseline manual driving scenario. \textbf{B–D}: Distribution of the impact speed reduction, injury risk reduction, and crash avoidance rate, with an automatic emergency system compared to the baseline manual driving scenario.}
    \label{fig:descriptives}
\end{figure}

\begin{figure}[!htb]
    \centering
    \includegraphics[scale = 1]{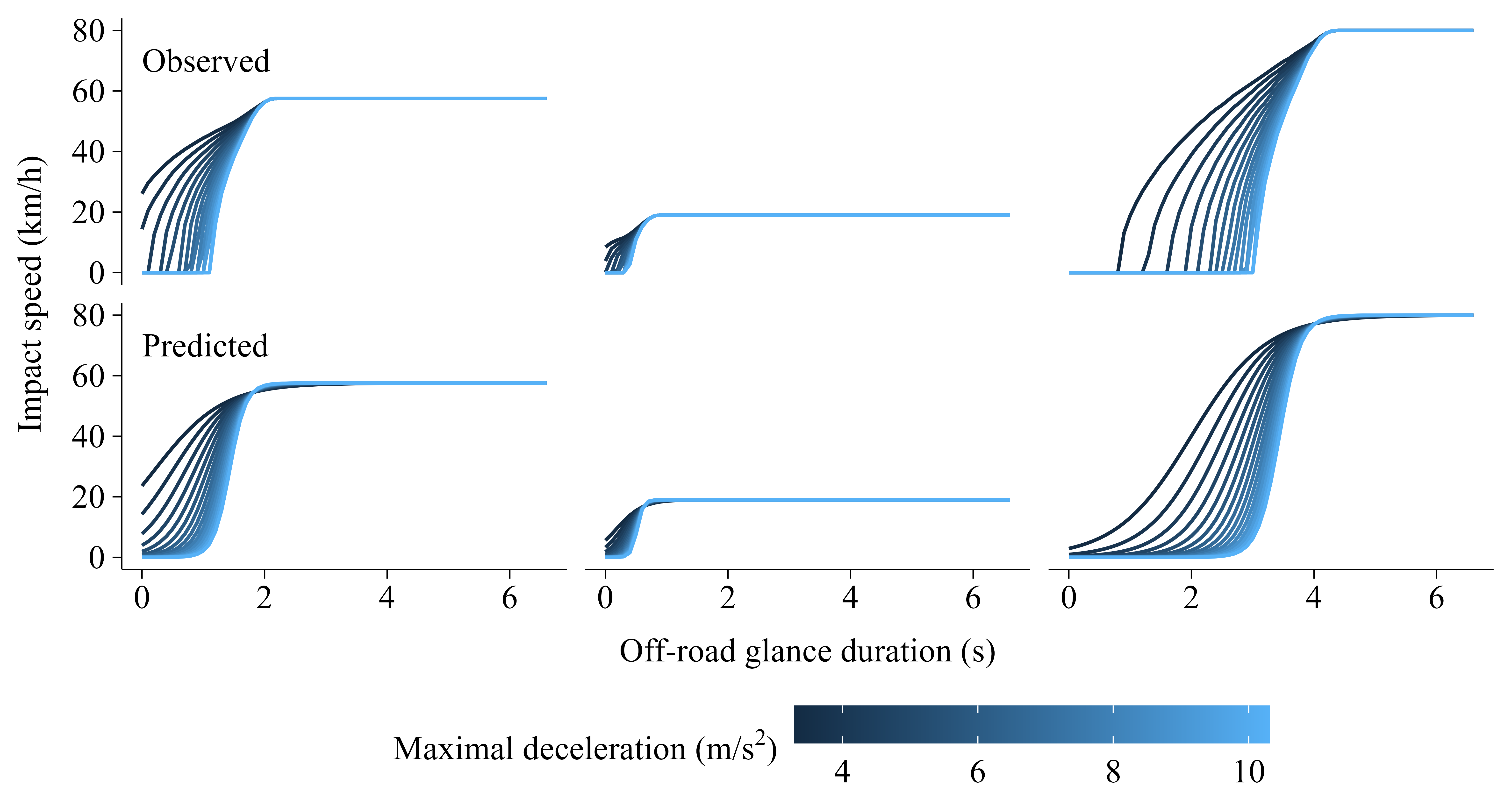}
    \caption{Impact speed response surface as a function of off-road glance duration and maximal deceleration during braking for counterfactual variations of three reconstructed rear-end crashes. \textbf{Top panel}: Observed impact speed. \textbf{Bottom panel}: Predicted impact speed using the quasi-binomial logistic regression model \eqref{eq:qblr}. The response values have been mapped from the model range $[0, 1]$ to the original range $[0, y_{\mathrm{max},k}]$, where $y_{\mathrm{max},k}$ is the maximal possible impact speed for the variations generated from case $k$, $k =1, \ldots, 44$.}
    \label{fig:response_surface}
\end{figure}

\end{document}